\magnification=\magstephalf
\input amstex
\documentstyle{amsppt}
\input xy
\xyoption{all}
\UseComputerModernTips

\topmatter
\title
Algebraic $K$-theory of the fraction field of topological $K$-theory
\endtitle
\rightheadtext{$K$-theory of the fraction field of $K$-theory}
%% \date
%% November 25th 2009
%% \enddate
\author
Christian Ausoni and John Rognes
\endauthor
\address Institute of Mathematics, University of Bonn, DE-53115 Bonn,
	Germany \endaddress
\email ausoni\@math.uni-bonn.de \endemail
\address Department of Mathematics, University of Oslo, NO--0316 Oslo,
	Norway \endaddress
\email rognes\@math.uio.no \endemail
\abstract
We compute the algebraic $K$-theory modulo~$p$ and~$v_1$ of the
$S$-algebra $\ell/p = k(1)$, using topological cyclic homology.  We use
this to compute the homotopy cofiber of a transfer map $K(L/p) \to
K(L_p)$, which we interpret as the algebraic $K$-theory of the ``fraction
field'' of the $p$-complete Adams summand of topological $K$-theory.
The results suggest that there is an arithmetic duality theorem for this
fraction field, much like Tate--Poitou duality for $p$-adic fields.
\endabstract
\endtopmatter

\define\ADer{\operatorname{ADer}}

\define\cC{\Cal C}

\define\cO{\Cal O}
\define\cok{\operatorname{cok}}
\define\dlog{\operatorname{dlog}}
\define\dlp{{\textstyle dp \over \textstyle p}}
\define\dlv{{\textstyle dv_1 \over \textstyle v_1}}
\define\Ext{\operatorname{Ext}}
\define\F{\Bbb F}
\define\ff{\text{\it ff\/}}
\define\Gal{\operatorname{Gal}}
\define\G{\Bbb G}
\redefine\H{\Bbb H}
\define\holim{\operatornamewithlimits{holim}}
\define\Q{\Bbb Q}
\define\Rlim{\mathop{{\lim}^1}}
\define\Spec{\operatorname{Spec}}
\define\Spnc{\operatorname{Spnc}}

\define\W{\Bbb W}
\define\Z{\Bbb Z}
\redefine\phi{\varphi}

\document

\head 1. Introduction \endhead

In this paper we continue the investigation of the algebraic $K$-theory
of topological $K$-theory initiated in~\cite{AuR02}. There we computed
the mod~$p$ and~$v_1$ homotopy groups of $K(\ell_p)$ at primes
$p\geq5$, denoted by $V(1)_* K(\ell_p)$ or $K/(p,v_1)_*(\ell_p)$,
where $\ell_p$ is the Adams summand of the $p$-complete connective
complex $K$-theory spectrum $ku_p$, with coefficients $\pi_*(\ell_p)
= \ell_{p*} = \Z_p[v_1]$, $|v_1|=2p-2$. Denoting by $P(v_2)$ the
polynomial algebra over $\F_p$ generated by $v_2$, $|v_2|=2p^2-2$, we
showed that $V(1)_* K(\ell_p)$ is essentially a free $P(v_2)$-module on
$4p+4$ generators. This is reminiscent, up to a chromatic shift of one
in the sense of stable homotopy theory, of the structure of the mod~$p$
algebraic $K$-theory $V(0)_* K(\Z_p) = K/(p)_*(\Z_p)$ of the $p$-adic
integers, which is a free $P(v_1)$-module on $p+3$ generators.

The structure of the  algebraic $K$-theory of $\Z_p$ can be conceptually
explained in terms of localization, Galois descent and motivic truncation.
First, there is a localization sequence
$$
K(\Z/p)\to K(\Z_p) \to K(\Q_p)
$$
relating the algebraic $K$-theory of $\Z_p$ to that of its fraction
field $\Q_p = \ff(\Z_p)$ and its residue field $\Z/p$.  Second, the
$v_1$-periodic mod~$p$ algebraic $K$-theory of $\Q_p$ is the target of
a Galois descent spectral sequence~\cite{Th85}
$$
E^2_{s,t} = H^{-s}_{Gal}(\Q_p; \F_p(t/2))
\Longrightarrow V(0)_{s+t} K(\Q_p) [v_1^{-1}] \,.
$$
The coefficient module $\F_p(t/2)$ can be interpreted as $V(0)_t
K(\bar\Q_p) [v_1^{-1}]$ by Suslin's Theorem~\cite{Su84}, and this spectral
sequence can be thought of as the continuous homotopy fixed-point spectral
sequence in mod~$p$ homotopy associated to the action of the absolute
Galois group of $\Q_p$ on $K(\bar\Q_p)$.  Third, the unlocalized mod~$p$
algebraic $K$-theory of $\Q_p$ is the target of a motivic spectral
sequence \cite{BL:ss}, \cite{FS02}
$$
E^2_{s,t} = H^{-s}_{mot}(\Q_p; \F_p(t/2))
\Longrightarrow V(0)_{s+t} K(\Q_p) \,,
$$
where
$$
H^r_{mot}(\Q_p; \F_p(i)) \cong \cases
H^r_{Gal}(\Q_p; \F_p(i)) & \text{for $r \le i$ \,,} \\
0 & \text{for $r>i$ \,.}
\endcases
$$
In other words, the mod~$p$ motivic cohomology of $\Q_p$ is a specific
truncation of its Galois cohomology.  In terms of Bloch's higher Chow
groups \cite{Bl86}, the vanishing of $H^r_{mot}(\Q_p; \F_p(i))$ for $r>i$
follows from the fact that there are no varieties of negative dimension
relative to $\Spec(\Q_p)$.  The identification with Galois cohomology
for $r\le i$ is part of the Beilinson conjectures \cite{BMS87}.

Our aim here is to conceptually understand $V(1)_* K(\ell_p)$ in the
same way as we understand $V(0)_* K(\Z_p)$, using the Galois theory for
commutative $S$-algebras developed by the second author~\cite{Rog08}.
The first task is therefore to identify the fraction field $\ff(\ell_p)$
of $\ell_p$.  There is a localization sequence
$$
K(\Z_p)\to K(\ell_p) \to K(L_p) \,,
$$
the existence of which was conjectured by the second author and
established by Blumberg and Mandell~\cite{BM08}, relating the algebraic
$K$-theory of $\ell_p$ to that of the periodic $p$-complete Adams
summand $L_p$.  The fact that $p$ is not invertible in $\pi_*(L_p) =
L_{p*} = \Z_p[v_1^{\pm1}]$, and the presence of the nontrivial residue
$S$-algebra $L/p = K(1)$ of $L_p$, give some indications as to why $L_p$
should not qualify as the fraction field we are after.

A naive guess would be to take the Bousfield localization $L_p[p^{-1}]$
of $L_p$ away from $p$, inverting $p$ in $\pi_*(L_p)$.  However,
$L_p[p^{-1}]$ is an $H\Q_p$-algebra, and we cannot possibly recover the
$v_2$-periodic mod~$p$ and~$v_1$ homotopy of $K(\ell_p)$ or $K(L_p)$
from that of $K(L_p[p^{-1}])$, which is $v_2$-torsion.  See Remark~3.12.
Therefore $L_p[p^{-1}]$ does not qualify as the desired fraction field
of $\ell_p$.  Instead, we have reasons to expect that if $\ff(\ell_p)$
is the ``correct'' fraction field, by which we mean that it will suit
the purpose of understanding algebraic $K$-theory by Galois descent, then
its algebraic $K$-theory $K(\ff(\ell_p))$ ought to fit in the following
$3\times3$-diagram of localization cofiber sequences, see Definition~3.9.
$$
\xymatrix{
K(\Z/p) \ar[r]^{i_*} \ar[d]^{\pi_*} &
K(\Z_p)\ar[r]^{j^*} \ar[d]^{\pi_*} &
K(\Q_p) \ar[d]^{\pi_*} \\
K(\ell/p) \ar[r]^{i_*} \ar[d]^{\rho^*} &
K(\ell_p) \ar[r]^{j^*} \ar[d]^{\rho^*} &
K(p^{-1}\ell_p) \ar[d]^{\rho^*} \\
K(L/p) \ar[r]^{i_*} &
K(L_p) \ar[r]^{j^*} &
K(\ff(\ell_p))\,.
}
\tag{1.1}
$$
We emphasize that for now $p^{-1}\ell_p$ and $\ff(\ell_p)$ are just formal
symbols, not the same as the Bousfield localizations $\ell_p[p^{-1}]$
and $L_p[p^{-1}]$.  The $V(1)$-homotopy groups of $K(\Z/p)$, $K(\Z_p)$
and $K(\ell_p)$ are known by~\cite{Qu72}, \cite{HM97} and~\cite{AuR02}.
It therefore remains to compute $V(1)_* K(\ell/p)$ and the
transfer maps $i_*$ and $\pi_*$ in order to be in
a position to evaluate the $V(1)$-homotopy of the iterated cofiber
$K(\ff(\ell_p))$. This is our main result, corresponding to Theorem~8.10
in the body of the paper.

\proclaim{Theorem 1.2}
There is an isomorphism of $E(\lambda_1, \lambda_2) \otimes P(v_2)$-modules
$$
\align
V(1)_* K(\ell/p) &\cong P(v_2) \otimes E(\bar\epsilon_1) \otimes
	\F_p\{1, \partial\lambda_2, \lambda_2, \partial v_2\} \\
&\qquad \oplus P(v_2) \otimes E(\dlog v_1) \otimes \F_p\{t^d
	v_2 \mid 0<d<p^2-p, p\nmid d\} \\
&\qquad \oplus P(v_2) \otimes E(\bar\epsilon_1) \otimes \F_p\{t^{dp}
	\lambda_2 \mid 0<d<p\} \,.
\endalign
$$
Here $|\lambda_1| = |\bar\epsilon_1| = 2p-1$, $|\lambda_2| = 2p^2-1$,
$|v_2| = 2p^2-2$, $|\dlog v_1| = 1$, $|\partial| = -1$ and $|t| = -2$.
This is a free $P(v_2)$-module of rank~$(2p^2-2p+8)$ and of zero Euler
characteristic, where $p\ge5$ is assumed.
\endproclaim

As explained in Section~2, the spectrum $\ell/p$ admits infinitely
many structures of an associative $\ell_p$-algebra, none of which are
commutative.  However, all of these $\ell_p$-algebras have equivalent
underlying $S$-algebra structures \cite{An:un}, and indeed it turns out
that the additive structure of $V(1)_* K(\ell/p)$ is independent of the
choice of $\ell_p$-algebra structure. The classes $\lambda_1$, $\lambda_2$
and~$v_2$ are present in $V(1)_* K(\ell_p)$, and in the statement
above we refer to the $V(1)_* K(\ell_p)$-module structure of $V(1)_*
K(\ell/p)$.  We prove this theorem by means of the cyclotomic trace
map to topological cyclic homology $TC(\ell/p)$.  On the way we evaluate
$V(1)_* THH(\ell/p)$, where $THH$ denotes topological Hochschild homology,
as well as $V(1)_* TC(\ell/p)$, see Proposition~5.4 and Theorem~8.8.

Defining $K(\ff(\ell_p))$ as the iterated cofiber in diagram~(1.1), we
can then evaluate its $V(1)$-homotopy.  The formulation given
in Theorem~9.4 readily implies the following statement.

\proclaim{Theorem 1.3}
Let $p\ge5$ be a prime. 
There is an isomorphism of $P(v_2^{\pm1})$-modules
$$
V(1)_* K(\ff(\ell_p)) [v_2^{-1}] \cong P(v_2^{\pm1}) \otimes \Lambda_*
$$
where
$$
\align
\Lambda_* &\cong E(\partial v_2, \dlog p, \dlog v_1) \\
  &\qquad \oplus E(\dlog v_1) \otimes \F_p\{ t^d \lambda_1 \mid 0<d<p\} \\
  &\qquad \oplus E(\dlog v_1) \otimes \F_p\{ t^d v_2 \dlog p \mid
  	0<d<p^2-p, p \nmid d\} \\
  &\qquad \oplus E(\dlog p) \otimes \F_p\{t^{dp} \lambda_2 \mid 0<d<p\}
\endalign
$$
is displayed in Figure~10.3.  Here $|\dlog p\,| = 1$, and the degrees
of the other classes are as in Theorem~1.2.  The localization homomorphism
$$
V(1)_* K(\ff(\ell_p)) \to V(1)_* K(\ff(\ell_p)) [v_2^{-1}]
$$
is an isomorphism in degrees $* \ge 2p$.
\endproclaim

This proves that the iterated cofiber $K(\ff(\ell_p))$ cannot be a
$K(\Q_p)$-module, since $V(1)_* K(\Q_p)$ is a torsion
$P(v_2)$-module.

In the final part of this paper, we conjecturally interpret the
computation of Theorem~1.3 in terms of Galois descent.  Indeed, the
second author conjectured that if $\Omega_1$ is a separable closure of
$\ff(\ell_p)$, then there is a homotopy equivalence
$$
L_{K(2)} K(\Omega_1) \simeq E_2 \,.
$$
Here $E_2$ is Morava's $E$-theory~\cite{GH04} with coefficients
$\pi_* E_2 = \W(\F_{p^2}) [[u_1]] [u^{\pm1}]$, and $L_{K(2)}$ denotes
Bousfield localization with respect to the Morava $K$-theory $K(2)$
with coefficients $\pi_* K(2) = \F_p[v_2^{\pm1}]$. The $v_2$-periodic
$V(1)$-homotopy groups of $K(\Omega_1)$ will then be given by
$$
V(1)_* K(\Omega_1) [v_2^{-1}] \cong \F_{p^2}[u^{\pm1}].
$$
We therefore conjecture to have a corresponding Galois descent spectral
sequence
$$
E^2_{s,t} = H^{-s}_{Gal}(\ff(\ell_p); \F_{p^2}(t/2))
\Longrightarrow V(1)_{s+t} K(\ff(\ell_p)) [v_2^{-1}]  \,,
$$
see Conjecture~10.5.  Starting with $V(1)_* K(\ff(\ell_p)) [v_2^{-1}]$
and working backwards, we give a conjectural description of the Galois
cohomology of $\ff(\ell_p)$ with coefficients in $V(1)_* K(\Omega_1)
[v_2^{-1}]$, see Figure~10.3.  This fits very well with the example of
the Galois cohomology of $\Q_p$ with coefficients in $V(0)_* K(\bar\Q_p)
[v_1^{-1}]$, with the difference that the absolute Galois group of
$\ff(\ell_p)$ has cohomological dimension $3$ instead of~$2$.  Also, there
appears to be a perfect arithmetic duality pairing in the conjectural
Galois cohomology of~$\ff(\ell_p)$, analogous to Tate--Poitou duality for
local number fields, which is not present in the analogous cohomology
rings for $\ell_p$ or $L_p$.  This indicates that $\ff(\ell_p)$ ought
to be a form of $S$-algebraic two-local field, mixing three different
residue characteristics.  As argued in section~3 of this paper, we expect
very similar results to hold for $\ell_p$ replaced by $ku_p$.

\medskip

The paper is organized as follows.  In section~2 we fix our
notations, show that $\ell/p$ admits the structure of an associative
$\ell$-algebra, and give a similar discussion for $ku/p$ and the periodic
versions $L/p$ and $KU/p$. We
give derived algebro-geometric and log-geometric interpretations
of these constructions.  In section~3, we discuss localization
sequences in algebraic $K$-theory.
Section~4 contains the computation of
the mod $p$ homology of $THH(\ell/p)$, and in section~5 we evaluate
its $V(1)$-homotopy. Sections~6, 7 and~8 deal with the computation of
$TC(\ell/p)$ and $K(\ell/p)$ in $V(1)$-homotopy.  In section~9 the various
computations are assembled in a diagram of localization sequences,
to evaluate $V(1)_* K(\ff(\ell_p))$.  In section~10, the structure of
$V(1)_* K(\ff(\ell_p))$ is interpreted in terms of Galois descent, yielding
a conjectural description of the Galois cohomology of $\ff(\ell_p)$.

\head 2. Base change squares of $S$-algebras \endhead

We fix some notations.  Let $p$ be a prime, even or odd for now.  Write
$X_{(p)}$ and $X_p$ for the $p$-localization and the $p$-completion,
respectively, of any spectrum or abelian group $X$.  Let $ku$ and $KU$
be the connective and the periodic complex $K$-theory spectra, with
homotopy rings $ku_* = \Z[u]$ and $KU_* = \Z[u^{\pm1}]$, where
$|u|=2$.  Let $\ell = BP\langle1\rangle$ and $L = E(1)$ be the
$p$-local Adams summands, with $\ell_* = \Z_{(p)}[v_1]$ and $L_* =
\Z_{(p)}[v_1^{\pm1}]$, where $|v_1| = 2p-2$.  The inclusion $\ell \to
ku_{(p)}$ maps $v_1$ to $u^{p-1}$.  Alternate notations in the
$p$-complete cases are $KU_p = E_1$ and $L_p = \widehat{E(1)}$.  These
ring spectra are all commutative $S$-algebras, in the sense that each
admits a unique $E_\infty$ ring spectrum structure.  See \cite{BR05}
for proofs of uniqueness in the periodic cases.

Let $ku/p$ and $KU/p$ be the connective and periodic
mod~$p$ complex $K$-theory spectra, with coefficients $(ku/p)_* =
\Z/p[u]$ and $(KU/p)_* = \Z/p[u^{\pm1}]$.  These are $2$-periodic
versions of the Morava $K$-theory spectra $\ell/p = k(1)$ and $L/p =
K(1)$, with $(\ell/p)_* = \Z/p[v_1]$ and $(L/p)_* = \Z/p[v_1^{\pm1}]$.
Each of these can be constructed as the cofiber of the multiplication
by $p$ map, as a module over the corresponding commutative
$S$-algebra.  For example, there is a cofiber sequence of $ku$-modules
$ku @>p>> ku @>i>> ku/p$.

Let $HR$ be the Eilenberg--Mac\,Lane spectrum of a ring $R$.  When $R$
is associative, $HR$ admits a unique associative $S$-algebra structure,
and when $R$ is commutative, $HR$ admits a unique commutative
$S$-algebra structure.  The zeroth Postnikov section defines unique
maps of commutative $S$-algebras $\pi \: ku \to H\Z$ and $\pi \: \ell
\to H\Z_{(p)}$, which can be followed by unique commutative $S$-algebra
maps to $H\Z/p$.

The $ku$-module spectrum $ku/p$ does not admit the structure of a
commutative $ku$-algebra.  It cannot even be an $E_2$ or $H_2$ ring
spectrum, since the homomorphism induced in mod~$p$ homology by the
resulting map $\pi \: ku/p \to H\Z/p$ of $H_2$ ring spectra would not
commute with the homology operation $Q^1(\bar\tau_0) = \bar\tau_1$ in
the target $H_*(H\Z/p; \F_p)$ \cite{BMMS86, III.2.3}.  Similar remarks
apply for $KU/p$, $\ell/p$ and $L/p$.  Associative algebra structures,
or $A_\infty$ ring spectrum structures, are easier to come by.  The
following result is a direct application of the methods of \cite{La01,
\S\S9--11}.  We adapt the notation of \cite{BJ02, \S3} to provide some
details in our case.

\proclaim{Proposition 2.1}
The $ku$-module spectrum $ku/p$ admits the structure of an associative
$ku$-algebra, but the structure is not unique.  Similar statements hold
for $KU/p$ as a $KU$-algebra, $\ell/p$ as an $\ell$-algebra and $L/p$
as an $L$-algebra.
\endproclaim

\demo{Proof}
We construct $ku/p$ as the (homotopy) limit of its Postnikov tower of
associative $ku$-algebras $P^{2m-2} = ku/(p, u^m)$, with coefficient rings
$ku/(p, u^m)_* = ku_*/(p, u^m)$ for $m\ge1$.  To start the induction,
$P^0 = H\Z/p$ is a $ku$-algebra via $i \circ \pi \: ku \to H\Z \to H\Z/p$.
Assume inductively for $m\ge1$ that $P = P^{2m-2}$ has been constructed.
We will define $P^{2m}$ by a (homotopy) pullback diagram
$$
\xymatrix{
P^{2m} \ar[r] \ar[d] & P \ar[d]^{in_1} \\
P \ar[r]^-{d} & P \vee \Sigma^{2m+1} H\Z/p
}
$$
in the category of associative $ku$-algebras.  Here
$$
d \in \ADer_{ku}^{2m+1}(P, H\Z/p)
\cong THH_{ku}^{2m+2}(P, H\Z/p)
$$
is an associative $ku$-algebra derivation of $P$ with values in
$\Sigma^{2m+1} H\Z/p$, and the group of such can be identified with the
indicated topological Hochschild cohomology group of $P$ over
$ku$.  We recall that these are the homotopy groups (cohomologically
graded) of the function spectrum $F_{P \wedge_{ku} P^{op}}(P, H\Z/p)$.
The composite map $pr_2 \circ d \: P \to \Sigma^{2m+1}
H\Z/p$ of $ku$-modules, where $pr_2$ projects onto the second wedge
summand, is restricted to equal the $ku$-module Postnikov $k$-invariant
of $ku/p$ in
$$
H^{2m+1}_{ku}(P; \Z/p) = \pi_0 F_{ku}(P, \Sigma^{2m+1} H\Z/p) \,.
$$
We compute that $\pi_*(P \wedge_{ku} P^{op}) = ku_*/(p, u^m) \otimes
E(\tau_0, \tau_{1,m})$, where $|\tau_0| = 1$, $|\tau_{1,m}| = 2m+1$ and
$E(-)$ denotes the exterior algebra on the given generators.  (For $p=2$,
the use of the opposite product is essential here \cite{An08, \S3}.)
The function spectrum description of topological Hochschild cohomology
leads to the spectral sequence
$$
\align
E_2^{**} &= \Ext^{**}_{\pi_*(P \wedge_{ku} P^{op})}(\pi_*(P), \Z/p) \\
&\cong \Z/p[y_0, y_{1,m}] \\
&\Longrightarrow THH^*_{ku}(P, H\Z/p) \,,
\endalign
$$
where $y_0$ and $y_{1,m}$ have cohomological bidegrees~$(1,1)$
and~$(1,2m+1)$, respectively.  The spectral sequence collapses at $E_2 =
E_\infty$, since it is concentrated in even total degrees.  In particular,
$$
\ADer^{2m+1}_{ku}(P, H\Z/p) \cong \F_p\{y_{1,m}, y_0^{m+1}\} \,.
$$
Additively, $H^{2m+1}_{ku}(P; \Z/p) \cong \F_p\{Q_{1,m}\}$ is generated
by a class dual to $\tau_{1,m}$, which is the image of $y_{1,m}$ under
left composition with $pr_2$.  It equals the $ku$-module $k$-invariant
of $ku/p$.  Thus there are precisely $p$ choices $d = y_{1,m} +
\alpha y_0^{m+1}$, with $\alpha \in \F_p$, for how to extend any given
associative $ku$-algebra structure on $P = P^{2m-2}$ to one on $P^{2m}
= ku/(p, u^{m+1})$.  In the limit, we find that there are an uncountable
number of associative $ku$-algebra structures on $ku/p = \holim_m P^{2m}$,
each indexed by the sequence of choices $\alpha \in \F_p$.

The periodic spectrum $KU/p$ can be obtained from $ku/p$ by Bousfield
$KU$-localization in the category of $ku$-modules \cite{EKMM97,
VIII.4}, which makes it an associative $KU$-algebra.  The
classification of periodic $S$-algebra structures is the same as in the connective
case, since the original $ku$-algebra structure on $ku/p$ can be
recovered from that on $KU/p$ by a functorial passage to the connective
cover.  To construct $\ell/p$ as an associative $\ell$-algebra, or
$L/p$ as an associative $L$-algebra, replace $u$ by $v_1$ in these
arguments.
\qed
\enddemo

By varying the ground $S$-algebra, we obtain the same conclusions about
$ku/p$ as a $ku_{(p)}$-algebra or $ku_p$-algebra, and about $\ell/p$ as
an $\ell_p$-algebra.  Similarly, for each $\nu\ge1$ we can realize the
$ku$-module spectrum $ku/p^\nu$, with $(ku/p^\nu)_* = \Z/p^\nu[u]$, as
an associative $ku$-algebra, and with a little care the $ku$-algebra
structures may be so chosen that the $p$-adic tower
$$
\dots \to ku/p^{\nu+1} \to ku/p^\nu \to \dots \to ku/p
\tag 2.2
$$
is one of associative $ku$-algebras.  Uncountably many choices of
structures remain, though.  Again, the same remarks apply over $\ell$,
where the $p$-adic tower $\{\ell/p^\nu\}_\nu$ can be chosen to be one
of associative $\ell$-algebras.

For each choice of $ku$-algebra structure on $ku/p$, the zeroth
Postnikov section $\pi \: ku/p \to H\Z/p$ is a $ku$-algebra map, with
the unique $ku$-algebra structure on the target.  Hence there is
a commutative square of associative $ku$-algebras
$$
\xymatrix{
ku \ar[r]^i \ar[d]^{\pi} & ku/p \ar[d]^{\pi} \\
H\Z \ar[r]^i & H\Z/p
}
\tag 2.3
$$
and similarly in the $p$-local and $p$-complete cases.  In view of the
weak equivalence $H\Z \wedge_{ku} ku/p \simeq H\Z/p$, this square
expresses the associative $H\Z$-algebra $H\Z/p$ as the base change of
the associative $ku$-algebra $ku/p$ along $\pi \: ku \to H\Z$.
Likewise, there is a commutative square of associative $\ell$-algebras
$$
\xymatrix{
\ell \ar[r]^i \ar[d]^{\pi} & \ell/p \ar[d]^{\pi} \\
H\Z_{(p)} \ar[r]^i & H\Z/p
}
\tag 2.4
$$
that expresses $H\Z/p$ as the base change of $\ell/p$ along $\ell \to
H\Z_{(p)}$, and similarly in the $p$-complete case.  By omission of
structure, these squares are also diagrams of $S$-algebras and
$S$-algebra maps.

\remark{Remark 2.5}
We might think algebro-geometrically about these diagrams.  In other
words, we would like to view the $S$-algebras in diagram~(2.3) as the
rings of functions on a diagram
$$
\xymatrix{
\Spec(\Z/p) \ar[r]^-i \ar[d]^{\pi} & \Spec(\Z) \ar[d]^{\pi} \\
\Spnc(ku/p) \ar[r]^-i & \Spec(ku)
}
$$
of affine derived schemes, in the sense of Jacob Lurie, or affine brave
new schemes in the sense of To{\"e}n and Vezzosi.  Here $\Spec(ku)$ has
the same underlying prime ideal space as $\Spec(\pi_0(ku)) = \Spec(\Z)$,
but it is ringed in commutative $S$-algebras, rather than in ordinary
commutative rings.  Over the open subset $\Spec(\Z[f^{-1}])$, where $f$
is a finite product of primes, the ring of functions is the commutative
$S$-algebra $ku[f^{-1}]$ with coefficient ring $ku_*[f^{-1}]$.  However,
since $ku/p$ is not a commutative $S$-algebra, we do not know how to
make good sense of its structure space.  At best, $i \: \Spnc(ku/p)
\to \Spec(ku)$ can be though of as a map from a non-commutative derived
scheme to the derived scheme $\Spec(ku)$.  The base change assertion
above then says that $\Spec(\Z/p)$ is the fiber of $\Spnc(ku/p)$ over
the closed subspace $\Spec(\Z)$ in $\Spec(ku)$, which in turn happens
to be an honest (commutative) scheme.
\endremark

\remark{Remark 2.6}
At this point we recall how $\Spec$ of the fraction field $F = \ff(R)$
of an integral domain $R$ is obtained geometrically by cutting away all
the proper closed subspaces of $\Spec(R)$.  The remaining open subset
represents the generic point of $\Spec(R)$.  For any localization $T$
of $R$, with $R \subset T \subset F$, we also have $F = \ff(T)$, so $F$
is simultaneously the fraction field of $R$ and of $T$.  In the case
of $\Spec(ku)$, with $R = ku$, the analogous (sub-)spaces appear to be
$\Spec(\Z)$ (included by $\pi$) and the $\Spnc(ku/p)$ (mapped by $i$)
for all primes $p$, so $\Spec$ of the fraction field $\ff(ku)$ should
consist of the complement in $\Spec(ku)$ of all of these spaces.
In the $p$-complete case, this means cutting out the two spaces
$\Spec(\Z_p)$ and $\Spnc(ku/p)$ from $\Spec(ku_p)$, meeting precisely
along their intersection $\Spec(\Z/p)$.  The remainder should be
$\Spec(\ff(ku_p))$.  We expect that the complement of $\Spec(\Z_p)$ in
$\Spec(ku_p)$ will play the role of $\Spec(KU_p)$, and that the complement
of $\Spec(\Z/p)$ in $\Spnc(ku/p)$ will play the role of $\Spnc(KU/p)$,
so the $S$-algebraic fraction field of $ku_p$ is also realized as the
$S$-algebraic fraction field of $KU_p$, denoted $\ff(ku_p) = \ff(KU_p)
= p^{-1}KU_p$.  Similarly, $\Spec$ of the fraction field
$$
\ff(\ell_p) = \ff(L_p) = p^{-1}L_p
$$
of $\ell_p$ or $L_p$ would be the complement of $\Spec(\Z_p)$ and
$\Spnc(\ell/p)$ in $\Spec(\ell_p)$, meeting at $\Spec(\Z/p)$.

Based on the algebraic $K$-theory computations to follow, we have strong
reasons to believe that the complement $\Spec(p^{-1}ku)$ of $\Spnc(ku/p)$
in $\Spec(ku)$ is not equal to the derived algebraic scheme associated to
the commutative $ku$-algebra $ku[p^{-1}]$ that is obtained by (algebraic)
Bousfield localization away from $p$.  See Remark~3.12 for the calculational
justification of our expectation.  Similarly, in the $p$-complete case
the complement $\Spec(p^{-1}ku_p)$ of $\Spnc(ku/p)$ in $\Spec(ku_p)$
is not equal to $\Spec(ku_p[p^{-1}])$, and in the periodic case the
complement $\Spec(p^{-1}KU_p) = \Spec(\ff(ku_p))$ is different from
$\Spec(KU_p[p^{-1}])$.  In fact, we do not expect that $\Spec(p^{-1}ku_p)$
and $\Spec(\ff(ku_p))$ are affine derived schemes in the sense that
they are equal to $\Spec(B)$ for a commutative $S$-algebra~$B$.  Still,
they may be derived stacks, in a sense approximated by the following
remark.
\endremark

\remark{Remark 2.7}
Since $\Spnc(ku/p)$ is not actually a closed subscheme of
$\Spec(ku_p)$, but only maps to the complement of the open
subscheme $\Spec(ku_p[p^{-1}])$, we expect that a derived scheme
$\Spec(p^{-1} ku_p)$ that is complementary to $\Spnc(ku/p)$, is
a derived algebro-geometric object which is intermediate between
$\Spec(ku_p[p^{-1}])$ and $\Spec(ku_p)$.  In forthcoming work, by
the second author together with Clark Barwick and Steffen Sagave, we
hope to give good sense to this intermediate object as a logarithmic
$S$-algebra $(ku_p, M)$, where $M$ is a pointed $E_\infty$ space.
See also Remark~3.3.
This makes $(ku_p, M)$ a logarithmic $S$-algebra and $\Spec(ku_p, M)$
an affine derived logarithmic scheme.  Following Martin Olsson's work
in the classical algebraic case \cite{Ol03}, this derived logarithmic
scheme can be faithfully represented by the derived stack representing
the $\infty$-category of derived logarithmic schemes over it.  Hence we
expect that this moduli stack of $\Spec(ku_p, M)$ will provide the desired
derived stack interpretation of $\Spec(p^{-1}ku_p)$.  Related topological
logarithmic structures on $ku_p$ will then provide useful interpretations
for $\Spec(KU_p)$ and $\Spec(\ff(ku_p))$.  See \cite{Rog09} for a set
of foundations of topological logarithmic geometry, and \cite{Rog:ltc}
for the construction of topological cyclic homology of a logarithmic
$S$-algebra.
\endremark

\head 3. Localization squares in algebraic $K$-theory \endhead

We are interested in the algebraic $K$-theory spectra of the
$S$-algebras considered in the previous section.  For an $S$-algebra
$B$, we adapt \cite{EKMM97, VI} and let $\cC_B$ be the category of
finite cell $B$-modules and $B$-module maps.  This is a category with
cofibrations and weak equivalences in the sense of \cite{Wa85, \S1},
with cofibrations the maps that are isomorphic to the inclusion of a
subcomplex, and weak equivalences the homotopy equivalences.
The algebraic
$K$-theory of~$B$, denoted $K(B)$, is defined to be the algebraic
$K$-theory of this category with cofibrations and weak equivalences.
When $B = HR$ is an Eilenberg--Mac\,Lane spectrum, it is known that $K(HR)
\simeq K(R)$ recovers Quillen's algebraic $K$-theory of the ring $R$.

If $\phi \: A \to B$ is a map of $S$-algebras, the exact functor $B
\wedge_A (-) \: \cC_A \to \cC_B$ induces a map $\phi^* \: K(A) \to K(B)$,
making $K(-)$ a covariant functor in the $S$-algebra.  We shall follow
the variance conventions of algebraic geometry, thinking of $K(-)$
as a contravariant functor in a (not yet generally defined) geometric
object $\Spec(-)$ associated to the $S$-algebra, so that $\phi^*$
is derived from the inverse image functor on sheaves.  This is done
to avoid confusion with the transfer maps $\phi_* \: K(B) \to K(A)$,
derived from the direct image functor on sheaves.

Suppose that $\phi \: A \to B$ makes $B$ a finite cell $A$-module.
Then each finite cell $B$-module $M$ can be viewed as an $A$-module along
$\phi$, and each $B$-cell attachment needed to build $M$ can be achieved
by attaching finitely many $A$-cells, one for each $A$-cell in $B$.
The resulting exact functor $\cC_B \to \cC_A$ induces the transfer map
$\phi_* \:  K(B) \to K(A)$.

When $A$ is a commutative $S$-algebra, the smash product $(-) \wedge_A
(-)$ induces a pairing on $\cC_A$ that makes $K(A)$ a commutative
$S$-algebra.  When $A$ is central in $B$, so that $B$ is an
$A$-algebra, the smash product makes $K(B)$ a $K(A)$-module, and the
maps $\phi^* \: K(A) \to K(B)$ and $\phi_* \: K(B) \to K(A)$ (when
defined) are $K(A)$-module maps.  In other words, the transfer is
a module map over its target.  When $B$ is commutative, then $K(B)$ is
likewise a commutative $S$-algebra, and $\phi^*$ is a commutative
$S$-algebra map.  The $K(A)$-linearity of $\phi_*$ can then be
expressed by the projection/reciprocity formula $\phi_*(\phi^*(x)
\wedge_B y) = x \wedge_A \phi_*(y)$.

We apply this to the maps $i \: ku \to ku/p$ and $\pi \: ku \to H\Z$.
Here $ku/p = ku \cup_p e^1$ and $H\Z = ku \cup_u e^3$ in the category
of $ku$-modules, meaning that $ku/p$ is obtained from $ku$ by attaching
a $1$-cell along $p$, and that $H\Z$ is obtained from $ku$ by attaching
a $3$-cell along $u$.   The transfer maps $i_*$ and $\pi_*$ and their
variants then define the upper left hand square in the following
commutative diagram, where all rows and columns are cofiber sequences.
$$
\xymatrix{
K(\Z/p) \ar[r]^{i_*} \ar[d]^{\pi_*} &
K(\Z) \ar[r]^{j^*} \ar[d]^{\pi_*} &
K(\Z[p^{-1}]) \ar[d]^{\pi_*} \\
K(ku/p) \ar[r]^{i_*} \ar[d]^{\rho^*} &
K(ku) \ar[r]^{j^*} \ar[d]^{\rho^*} &
K(p^{-1}ku) \ar[d]^{\rho^*} \\
K(KU/p) \ar[r]^{i_*} &
K(KU) \ar[r]^{j^*} &
K(p^{-1}KU)
}
\tag 3.1
$$
The upper row is a cofiber sequence by Quillen's localization theorem
for Dedekind domains \cite{Qu73}.  
Andrew Blumberg and Mike Mandell \cite{BM08} have proved that the
middle column is a cofiber sequence, as conjectured by
the second author.  The result is analogous to that of Quillen, in that
$ku/u \simeq H\Z$ and $ku[u^{-1}] \simeq KU$.  
The same argument shows that the left hand column
is a cofiber sequence.

We complete the $3 \times 3$ diagram of cofiber sequences by making the
following ad hoc definition.

\definition{Definition 3.2}
Let $K(p^{-1}ku)$ and $K(p^{-1}KU)$ denote the (homotopy) cofibers of the
transfer maps $i_* \: K(ku/p) \to K(ku)$ and $i_* \: K(KU/p) \to
K(KU)$, respectively.
\enddefinition

\remark{Remark 3.3}
The theory of logarithmic structures on $S$-algebras, see \cite{Rog09},
is designed to model intermediate objects between a commutative
$S$-algebra $B$, like $ku$, and its localizations, like $ku[p^{-1}]$.
A logarithmic structure on $B$ is given by a commutative $S$-algebra map
$S[M] \to B$ from the spherical monoid ring on a ``commutative monoid''
$M$, where the latter term should be interpreted in a sufficiently
flexible category to also encompass $E_\infty$ spaces.  In joint work with
Steffen Sagave, the second author has defined a category of logarithmic
modules over a logarithmic $S$-algebra $(B, M)$, and one may define
the logarithmic $K$-theory $K(B, M)$ as the algebraic $K$-theory of the
category of suitably finite objects in that module category.
It still remains to be seen if one can realize $K(p^{-1} ku)$ as $K(ku,
M)$ for a suitable logarithmic $S$-algebra $(ku, M)$.
\endremark

\remark{Remark 3.4}
There is a similar $3 \times 3$ diagram to~(3.1), in which the left
hand column is replaced by the sum of corresponding columns for all
primes $p$, and $i_*$ is replaced by the sum of the transfer maps.  By
Quillen's localization sequence the upper right hand term can be
written as $K(\Q)$.  We expect that the lower right hand term then
plays the role of the algebraic $K$-theory of the $S$-algebraic
fraction field of $ku$, denoted $K(\ff(ku))$.  The hypothetical
fraction field has a valuation with uniformizer $u$, and the middle
right hand term $K(\cO)$ in the $3 \times 3$ diagram plays the role of
the algebraic $K$-theory of an $S$-algebraic valuation ring $\cO$ of
$\ff(ku)$.  However, we do not expect that this fraction field
$\ff(ku)$ is the rationalization $L_0 KU = KU\Q$ of the $S$-algebra
$KU$.  This expectation is again supported by the algebraic $K$-theory
computations that follow.
\endremark

\medskip

We shall use the topological cyclic homology $TC(-)$ and the cyclotomic
trace map $trc \: K(B) \to TC(B)$ of \cite{BHM93} to compute the
algebraic $K$-theory spectra and maps in the diagram~(3.1) above, after
making two modifications.  Firstly, we will replace $ku$ by its
$p$-completion $ku_p$, and secondly we will pass to the Adams summand
$\ell_p$ of $ku_p$.  We will now motivate these two changes, and discuss
their effect on the algebraic $K$-theory diagram.

First, the cyclotomic trace map is very close to a $p$-complete
equivalence when applied to connective $p$-complete $S$-algebras $B$
with $\pi_0(B) = \Z/p$ or $\Z_p$, or a finite extension of these.  More
precisely, $K(B)_p$ is the connective cover of $TC(B)_p$ in these cases
\cite{HM97}.  Therefore we shall pass to the $p$-complete situation,
replacing $ku$ and $H\Z$ by $ku_p$ and $H\Z_p$, respectively.  In view
of \cite{Du97}, the iterated (homotopy) fiber of the square
$$
\xymatrix{
K(ku)_p \ar[r]^{\kappa^*} \ar[d]^{\pi^*} & K(ku_p)_p \ar[d]^{\pi^*} \\
K(\Z)_p \ar[r]^{\kappa^*} & K(\Z_p)_p \,,
}
$$
where $\kappa \: ku \to ku_p$ and $\kappa \: H\Z \to H\Z_p$ denote the
completion maps, is homotopy equivalent to the iterated fiber
of the corresponding square of topological cyclic homology spectra, and
the latter is trivial, since $TC(B)_p \simeq TC(B_p)_p$.  Therefore the
fiber of $K(ku)_p \to K(ku_p)_p$ is homotopy equivalent to the
fiber of $K(\Z)_p \to K(\Z_p)_p$, so the change created by passing to
the $p$-complete case is as well understood for $ku$ and $ku_p$ as in
the number ring case.

In fact, our calculations shall concentrate on the $v_2$-periodic behavior
of the algebraic $K$-theory of the $S$-algebras related to topological
$K$-theory, rather than the subsequent issues of divisibility by $v_0 =
p$ and $v_1 = u^{p-1}$.  We achieve this by working with homotopy with
mod~$p$ and~$v_1$ coefficients, i.e., with the $V(1)$-homotopy functor
$X \mapsto V(1)_* X = \pi_*(V(1) \wedge X)$, where $V(1) = S/(p, v_1)$
is the Smith--Toda complex \cite{Sm70} defined for $p$ odd as the mapping
cone of the Adams self-map $v_1 \: \Sigma^{2p-2} V(0) \to V(0)$ of the mod~$p$
Moore spectrum $V(0) = S/p = S \cup_p e^1$.
Hence there is a cofiber sequence
$$
\Sigma^{2p-2} V(0) @>v_1>> V(0) @>i_1>> V(1) @>j_1>> \Sigma^{2p-1} V(0) \,.
$$
The Smith--Toda complex $V(1)$ is a homotopy commutative and associative
ring spectrum for $p \ge 5$ \cite{Ok84}, and its $BP$-homology satisfies
$BP_* V(1) \cong BP_*/(p, v_1)$.  In this sense $V(1)$ generalizes the
mod~$p$ Moore spectrum representing mod~$p$ homotopy, and represents
mod~$p$ and~$v_1$ homotopy.  There is an element $v_2 \in \pi_{2p^2-2}
V(1)$, with image $v_2$ in $BP_* V(1)$.  So $V(1)$ is a finite complex of
type~$2$ in the sense of \cite{HoSm98}, with telescope $V(1)[v_2^{-1}]$ the
spectrum that represents $v_2$-periodic homotopy, $V(1)_*(X)[v_2^{-1}]$.
The composite map $\beta_{1,1} = i_1 j_1 \: V(1) \to \Sigma^{2p-1} V(1)$
defines the primary $v_1$-Bockstein homomorphism acting naturally on
$V(1)_* X$.

It is known by \cite{BM94} and \cite{BM95} that $V(1)_* K(\Z_p)$ is a
finite $\Z/p$-module.  Furthermore, the announced proof of the Bloch--Kato
conjecture by Vladimir Voevodsky and Markus Rost also implies the
Lichtenbaum--Quillen conjecture on the algebraic $K$-theory of rings of
integers in number fields, which in turn implies that $V(1)_* K(\Z)$ is
a finite $\Z/p$-module.  So the common fiber of the maps $K(\Z)_p \to
K(\Z_p)_p$ and $K(ku)_p \to K(ku_p)_p$ has (totally) finite
$V(1)$-homotopy.  Therefore, replacing $K(\Z)$ by $K(\Z_p)$ and $K(ku)$
by $K(ku_p)$ in diagram~(3.1), and still defining the lower row and the
right hand column so as to obtain a $3 \times 3$ square of
cofibrations, only changes the $V(1)$-homotopy in the diagram by a
finite $\Z/p$-module.  In this case the upper right hand corner is
$K(\Q_p)$, since $\Z_p[p^{-1}] = \Q_p$, and the lower right hand corner
$p^{-1}KU_p$ plays the role of the $S$-algebraic fraction field of $ku_p$,
denoted $\ff(ku_p)$.  It follows that the $v_2$-periodic homotopy,
given by $V(1)_* K(ku) [v_2^{-1}]$, etc., does not change at all:
$$
V(1)_* K(p^{-1}KU) [v_2^{-1}] \cong
V(1)_* K(\ff(ku_p)) [v_2^{-1}] \,.
\tag 3.5
$$

Along the same lines, we might have started with the version of
diagram~(3.1) discussed in Remark~3.4, inverting all primes in $\Z$
rather than just $p$.  This adds a copy of the cofiber sequence
$K(\Z/r) @>\pi_*>> K(ku/r) @>>> K(KU/r)$ to the left hand column, for
each prime $r \ne p$.  Here $\pi^* \: K(ku/r)_p \to K(\Z/r)_p$ is
easily seen to be an equivalence, using the $BGL(-)^+$-model for
algebraic $K$-theory, and from the $F\psi^r$-model for $K(\Z/r)$ from
\cite{Qu72} it is clear that $V(1)_* K(\Z/r)$ is finite and
concentrated in degrees $0 \le * < 2p-2$.  Therefore, replacing the
right hand column by the cofiber sequence $K(\Q) @>\pi_*>> K(\cO) \to
K(\ff(ku))$, where $\cO = \cO_{\ff(ku)}$ is the valuation ring of the
fraction field of $ku$, only changes the $V(1)$-homotopy in a finite
range of degrees, and the $v_2$-periodic homotopy is unchanged:
$$
V(1)_* K(p^{-1}KU) [v_2^{-1}] \cong V(1)_* K(\ff(ku)) [v_2^{-1}] \,.
\tag 3.6
$$

Second, the Hurewicz image of the Bott element $u$ is nilpotent in the
mod~$p$ homology of the spectrum $ku$, since $u^{p-1} = v_1$ has Adams
filtration~$1$, and this significantly complicates the algebra involved
in computing the topological Hochschild homology $THH(-)$ and the
topological cyclic homology $TC(-)$ of $ku$, compared to the
computations for its Adams summand $\ell$.  Compare \cite{MS93},
\cite{AuR02} and \cite{Au05}.  Therefore, for the detailed
calculations in this paper we shall restrict attention to
the Adams summands $\ell_p$ and $\ell/p$ of $ku_p$ and $ku/p$,
respectively.  The map $L_p \to KU_p$ is a $\Delta$-Galois extension of
commutative $S$-algebras, in the sense of \cite{Rog08, \S4.1}, where
$\Delta = (\Z/p)^\times \cong \Z/(p-1)$ acts on $KU_p$ through the
$p$-adic Adams operations fixing $L_p$.  Generalizing the
Lichtenbaum--Quillen conjecture, it is therefore to be expected that
the map
$$
\Phi \: K(L_p) \to K(KU_p)^{h\Delta}
$$
is (close to) an equivalence after $p$-adic completion, and similarly
for the fraction fields.  Indeed, the first expectation is known to be
correct:  Consider the map of horizontal cofiber sequences
$$
\xymatrix{
K(\Z_p)_p \ar[r]^{\pi_*} \ar[d]^{p-1} &
K(\ell_p)_p \ar[r] \ar[d]^{\phi} &
K(L_p)_p \ar[d]^{\Phi} \\
K(\Z_p)^{h\Delta}_p \ar[r]^{\pi_*^{h\Delta}} &
K(ku_p)^{h\Delta}_p \ar[r] &
K(KU_p)^{h\Delta}_p \,,
}
\tag 3.7
$$
where the lower row is obtained from the Blumberg--Mandell localization
sequence by taking $\Delta$-homotopy fixed points ($\Delta$ acts
trivially on $\Z_p = \pi_0(ku_p)$), and the upper row is its analogue
for the Adams summand.  The left hand vertical map multiplies by the
ramification index $(p-1)$, since the uniformizer $v_1$ in $L_p$ maps
to that power of the uniformizer $u$ in $KU_p$.  More directly, the
upper map $\pi_*$ takes the unit class of $H\Z_p$ to that of the
complex $v_1 \: \Sigma^{2p-2} \ell_p \to \ell_p$, which extends along
$\ell_p \to ku_p$ to the class of $v_1 = u^{p-1} \: \Sigma^{2p-2} ku_p
\to ku_p$.  This is equivalent to $(p-1)$ times the class of $u \:
\Sigma^2 ku_p \to ku_p$, which is the image of the unit class of
$H\Z_p$ under the lower map $\pi_*^{h\Delta}$.  Either way, the left
hand vertical map is a $p$-adic equivalence.  Thus the middle vertical
map $\phi$ is an equivalence if and only if the right hand vertical map
$\Phi$ is one.  These equivalent conclusions are indeed precisely
verified by the calculations of the first author in \cite{Au05},
showing that $V(1)_* K(\ell_p) \cong V(1)_* K(ku_p)^{\Delta}$, from
which the corresponding assertion for $p$-adic integral homotopy groups
follows by a Bockstein argument.

Conversely, the later calculations of \cite{Au:tcku} show that there is
an element $b \in V(1)_{2p+2} K(ku)$ with $b^{p-1} = -v_2$, related to
$a_{(0,1)} \in K(2)_{2p+2} K(\Z/p, 2)$ \cite{RaWi80, \S9} via the
composite map
$$
K(\Z/p,2) \to K(\Z,3) \to BGL_1(ku) \to K(ku)
$$
and the unit map $V(1) \to V(1) \wedge E(2) = K(2)$.  Furthermore,
there is an isomorphism
$$
P(b^{\pm1}) \otimes_{P(v_2^{\pm1})} V(1)_* K(\ell_p) [v_2^{-1}]
	\cong V(1)_* K(ku_p) [v_2^{-1}]
$$
where $P(x^{\pm1})$ denotes the Laurent polynomial algebra over $\F_p$
generated by $x$.  The left hand side can also be written as the algebraic
extension
$$
V(1)_* K(\ell_p)[v_2^{-1}, b]/(b^{p-1}+v_2) \,.
$$
Here $\Delta$ acts faithfully on $b$, fixing $b^{p-1} = -v_2$, so this
both shows how the $v_2$-periodic algebraic $K$-theory of $ku_p$ is
generated by that of $\ell_p$ and the element~$b$, and conversely how
the $v_2$-periodic algebraic $K$-theory of $\ell_p$ can be recovered
from that of $ku_p$ as the fixed points of the Galois action.  In view
of the cofiber sequences~(3.7) above, the corresponding formulas also
hold with $L_p$ and $KU_p$ in place of $\ell_p$ and $ku_p$.

We are therefore confident that $ku$-based calculations of $K(ku/p)$,
$K(ku_p)$ and $K(\ff(ku_p))$, like the $\ell$-based calculations of
$K(\ell/p)$, $K(\ell_p)$ and $K(\ff(\ell_p))$ to be presented here, will
verify the, for now conjectural, formula
$$
P(b^{\pm1}) \otimes_{P(v_2^{\pm1})} V(1)_* K(\ff(\ell_p)) [v_2^{-1}]
	\cong V(1)_* K(\ff(ku_p)) [v_2^{-1}]
\tag 3.8
$$
and show that $V(1)_* K(\ff(\ell_p)) [v_2^{-1}]$ can be recovered as the
$\Delta$-invariants of the right hand side.  Thus we expect these
more complicated calculations to give qualitatively the same answers.
One reason to perform such more complicated calculations, other than
to confirm our expectations, would be to precisely identify the image
of $V(1)_* K(\ff(ku_p))$ in $V(1)_* K(\ff(ku_p)) [v_2^{-1}]$, and to check
that the map inverting $v_2$ is very close to being injective, but we will
not address this matter in the present paper.

With this we are done with our justifications for concentrating on the
following framework.  We collect the relevant definitions in one place,
for the reader's convenience.

\definition{Definition 3.9}
Let $\ell_p$ be the $p$-complete connective Adams summand, which is
a well-defined commutative $S$-algebra, with $\pi_* \ell_p =
\Z_p[v_1]$.  For each choice of associative $\ell_p$-algebra structure
on $\ell/p$, the $p$-complete form
$$
\xymatrix{
\ell_p \ar[r]^i \ar[d]^{\pi} & \ell/p \ar[d]^{\pi} \\
H\Z_p \ar[r]^i & H\Z/p
}
$$
of diagram~(2.4) is a base change square of associative $\ell_p$-algebras.
Each map $i \: \ell_p \to \ell/p$, $i \: H\Z_p \to H\Z/p$, $\pi \: \ell_p
\to H\Z_p$ and $\pi \: \ell/p \to H\Z/p$ makes the target a finite 2-cell
complex over the source.  The respective transfer maps in algebraic
$K$-theory define the upper left hand square in the diagram~(3.10) below.

The upper row, and the left hand and middle columns, are cofiber
sequences by the algebraic $K$-theory localization sequences of Quillen
and Blumberg--Mandell, respectively.  We define $K(p^{-1}\ell_p)$ and
$K(p^{-1}L_p)$ to be the cofibers of the transfer maps $i_*$
in the middle and lower rows, respectively.  This makes $K(p^{-1}L_p)$ the
iterated cofiber of the upper left hand square.  Furthermore,
we write $\ff(\ell_p) = \ff(L_p)$ for the symbol $p^{-1}L_p$, and
consider it as an $S$-algebraic fraction field of~$\ell_p$.  Likewise,
we consider the symbol $p^{-1}\ell_p$ as an $S$-algebraic valuation
ring for a valuation on $\ff(\ell_p)$ with uniformizer $v_1$.  Then the
following diagram is a $3 \times 3$ square of cofiber sequences.
$$
\xymatrix{
K(\Z/p) \ar[r]^{i_*} \ar[d]^{\pi_*} &
K(\Z_p) \ar[r]^{j^*} \ar[d]^{\pi_*} &
K(\Q_p) \ar[d]^{\pi_*} \\
K(\ell/p) \ar[r]^{i_*} \ar[d]^{\rho^*} &
K(\ell_p) \ar[r]^{j^*} \ar[d]^{\rho^*} &
K(p^{-1}\ell_p) \ar[d]^{\rho^*} \\
K(L/p) \ar[r]^{i_*} &
K(L_p) \ar[r]^{j^*} &
K(\ff(\ell_p))
}
\tag 3.10
$$
\enddefinition

\remark{Remark 3.11}
We gave the existence of a non-trivial $S$-algebra homomorphism $L_p
\to L/p$ as one argument for why $L_p$ is not an $S$-algebraic field.
However, this is in an essential way a non-commutative phenomenon.
It follows from \cite{BMMS86, III.4.1} that any commutative $S$-algebra with
$p=0$ in $\pi_0$ is an $H\Z/p$-module, hence any commutative $L_p$-algebra
with $p=0$ in $\pi_0$ is trivial, since $H\Z/p \wedge L_p \simeq *$.
Algebraic $K$-theory is a functor of associative $S$-algebras, not
just of commutative $S$-algebras, so this should not be perceived as an
insurmountable obstacle, but it still illustrates some of the subtleties
involved in thinking algebro-geometrically about these $S$-algebras.
\endremark

\remark{Remark 3.12}
It is time to explain why we do not think of $K(p^{-1}ku)$ as the
algebraic $K$-theory of the $S$-algebra $ku[p^{-1}]$, etc.  For primes
$p\ge5$ we computed $V(1)_* K(\ell_p)$ in \cite{AuR02}, and we shall compute
$V(1)_* K(\ell/p)$ in Theorem~8.10 below.  The results are different,
also after inverting $v_2$, so $V(1)_* K(p^{-1}\ell_p) [v_2^{-1}]$ is
certainly nonzero.  This is a direct summand of $V(1)_* K(p^{-1}ku_p)
[v_2^{-1}]$, which is isomorphic to $V(1)_* K(p^{-1}ku) [v_2^{-1}]$,
$V(1)_* K(\ff(ku_p)) [v_2^{-1}]$ and $V(1)_* K(\ff(ku)) [v_2^{-1}]$,
so all of these are also nonzero.  On the other hand, $ku_p$ is an
$S_p$-algebra, so the localization $ku_p[p^{-1}]$ is an $S_p[p^{-1}]
\simeq H\Q_p$-algebra.  Therefore $V(1)_* K(ku_p[p^{-1}])$ is a module
over $V(1)_* K(\Q_p)$, which is finite, so $V(1)_* K(ku_p[p^{-1}])
[v_2^{-1}]$ is a module over $V(1)_* K(\Q_p) [v_2^{-1}]$, which is zero.
Hence $V(1)_* K(ku_p[p^{-1}]) [v_2^{-1}] = 0$.  Therefore $p^{-1}ku_p
\not\simeq ku_p[p^{-1}]$ and $\ff(ku_p) \not\simeq KU_p[p^{-1}]$.
\endremark

\head 4. Topological Hochschild homology \endhead

In this and the following sections, we shall compute the $V(1)$-homotopy
of the topological Hochschild homology $THH(-)$ and topological cyclic
homology $TC(-)$ of the $S$-algebras in the upper left hand square in
diagram~(3.10), for primes $p\ge5$.  Passing to connective covers, this
also computes the $V(1)$-homotopy of the algebraic $K$-theory spectra
appearing in that square, which by the projection formula for the transfer
maps allows us to compute the $V(1)$-homotopy of the remaining entries
in the diagram, including $K(p^{-1}L_p) = K(\ff(\ell_p))$.  With these
coefficients, or more generally, after $p$-adic completion, the functors
$THH$ and $TC$ are insensitive to $p$-completion in the argument, so we
shall simplify the notation slightly by working with the associative
$S$-algebras $\ell$ and $H\Z_{(p)}$ in place of $\ell_p$ and $H\Z_p$.
For ordinary rings $R$
we almost always shorten notations like $THH(HR)$ to $THH(R)$.

The computations follow the strategy of \cite{Bo:zzp}, \cite{BM94},
\cite{BM95} and \cite{HM97} for $H\Z/p$ and $H\Z$, and of \cite{MS93}
and \cite{AuR02} for $\ell$.  See also \cite{AnR05, \S\S4--7} for further
discussion of the $THH$-part of such computations.  In this section
we shall compute the mod~$p$ homology of the topological Hochschild
homology of $\ell/p$ as a module over the corresponding homology for
$\ell$, for any odd prime~$p$.

We write $E(x) = \F_p[x]/(x^2)$ for the exterior algebra, $P(x) = \F_p[x]$
for the polynomial algebra and $P(x^{\pm1}) = \F_p[x, x^{-1}]$ for the
Laurent polynomial algebra on one generator $x$, and similarly for a
list of generators.  We will also write $\Gamma(x) = \F_p\{ \gamma_i(x)
\mid i \ge 0\}$ for the divided power algebra, with $\gamma_i(x) \cdot
\gamma_j(x) = (i,j) \gamma_{i+j}(x)$, where $(i,j) = (i+j)!/i!j!$ is the
binomial coefficient.  We use the obvious abbreviations $\gamma_0(x) =
1$ and $\gamma_1(x) = x$.  Finally, we write $P_h(x) = \F_p[x]/(x^h)$ for
the truncated polynomial algebra of height $h$, and recall the isomorphism
$\Gamma(x) \cong P_p(\gamma_{p^e}(x) \mid e\ge0)$ in characteristic~$p$.

We write $H_*(-)$ for homology with mod~$p$ coefficients.  It takes
values in $A_*$-comodules, where $A_*$ is the dual Steenrod algebra
\cite{Mi58}.  Explicitly (for $p$ odd),
$$
A_* = P(\bar\xi_k \mid k\ge1) \otimes E(\bar\tau_k \mid k\ge0)
$$
with coproduct
$$
\psi(\bar\xi_k) = \sum_{i+j=k} \bar\xi_i \otimes \bar\xi_j^{p^i}
$$
and
$$
\psi(\bar\tau_k) = 1 \otimes \bar\tau_k + \sum_{i+j=k}  \bar\tau_i
\otimes \bar\xi_j^{p^i} \,.
$$
Here $\bar\xi_0 = 1$, $\bar\xi_k = \chi(\xi_k)$ has degree $2(p^k-1)$
and $\bar\tau_k = \chi(\tau_k)$ has degree $2p^k-1$, where $\chi$ is
the canonical conjugation \cite{MM65}.  Then the zeroth
Postnikov sections induce identifications
$$
\align
H_*(H\Z_{(p)}) &=  P(\bar\xi_k \mid k\ge1) \otimes E(\bar\tau_k \mid k\ge1) \\
H_*(\ell) &=  P(\bar\xi_k \mid k\ge1) \otimes E(\bar\tau_k \mid k\ge2) \\
H_*(\ell/p) &=  P(\bar\xi_k \mid k\ge1) \otimes E(\bar\tau_0, \bar\tau_k
\mid k\ge2)
\endalign
$$
as $A_*$-comodule subalgebras of $H_*(H\Z/p) = A_*$.
We often make use of the following $A_*$-comodule coactions
$$
\align
\nu(\bar\tau_0) &= 1 \otimes \bar\tau_0 + \bar\tau_0 \otimes 1 \\
\nu(\bar\xi_1) &= 1 \otimes \bar\xi_1 + \bar\xi_1 \otimes 1 \\
\nu(\bar\tau_1) &= 1 \otimes \bar\tau_1 + \bar\tau_0 \otimes \bar\xi_1
+ \bar\tau_1 \otimes 1 \\
\nu(\bar\xi_2) &= 1 \otimes \bar\xi_2 + \bar\xi_1 \otimes \bar\xi_1^p +
\bar\xi_2 \otimes 1 \\
\nu(\bar\tau_2) &= 1 \otimes \bar\tau_2 + \bar\tau_0 \otimes \bar\xi_2
+ \bar\tau_1 \otimes \bar\xi_1^p + \bar\tau_2 \otimes 1 \,.
\endalign
$$

The B{\"o}kstedt spectral sequences
$$
E^2_{**}(B) = HH_*(H_*(B)) \Longrightarrow H_*(THH(B))
$$
for the commutative $S$-algebras $B = H\Z/p$, $H\Z_{(p)}$ and $\ell$
begin
$$
\align
E^2_{**}(\Z/p) &= A_* \otimes E(\sigma\bar\xi_k \mid k\ge1) \otimes
\Gamma(\sigma\bar\tau_k \mid k\ge0) \\
E^2_{**}(\Z_{(p)}) &= H_*(H\Z_{(p)}) \otimes E(\sigma\bar\xi_k \mid
k\ge1) \otimes \Gamma(\sigma\bar\tau_k \mid k\ge1) \\
E^2_{**}(\ell) &= H_*(\ell) \otimes E(\sigma\bar\xi_k \mid k\ge1)
\otimes \Gamma(\sigma\bar\tau_k \mid k\ge2) \,.
\endalign
$$
They are (graded) commutative $A_*$-comodule algebra spectral sequences,
and there are differentials
$$
d^{p-1}(\gamma_j\sigma\bar\tau_k) = \sigma\bar\xi_{k+1}
\cdot \gamma_{j-p}\sigma\bar\tau_k
$$
for $j\ge p$ and $k\ge0$, see \cite{Bo:zzp}, \cite{Hu96}
or \cite{Au05, 4.3}, leaving
$$
\align
E^\infty_{**}(\Z/p) &= A_* \otimes P_p(\sigma\bar\tau_k \mid k\ge0) \\
E^\infty_{**}(\Z_{(p)}) &= H_*(H\Z_{(p)}) \otimes E(\sigma\bar\xi_1)
\otimes P_p(\sigma\bar\tau_k \mid k\ge1) \\
E^\infty_{**}(\ell) &= H_*(\ell) \otimes E(\sigma\bar\xi_1,
\sigma\bar\xi_2) \otimes P_p(\sigma\bar\tau_k \mid k\ge2) \,.
\endalign
$$
The inclusion of $0$-simplices $\eta \: B \to THH(B)$ is split for
commutative $B$ by the augmentation $\epsilon \: THH(B) \to B$.  Thus
there are unique representatives in B{\"o}kstedt filtration~$1$, with
zero augmentation, for each of the classes $\sigma x$.  They correspond
to $1 \otimes x - x \otimes 1$ in the Hochschild complex, or just $1
\otimes x$ in the normalized Hochschild complex.  There are
multiplicative extensions $(\sigma\bar\tau_k)^p = \sigma\bar\tau_{k+1}$
for $k\ge0$, so
$$
\aligned
H_*(THH(\Z/p)) &= A_* \otimes P(\sigma\bar\tau_0) \\
H_*(THH(\Z_{(p)})) &= H_*(H\Z_{(p)}) \otimes E(\sigma\bar\xi_1) \otimes
P(\sigma\bar\tau_1) \\
H_*(THH(\ell)) &= H_*(\ell) \otimes E(\sigma\bar\xi_1, \sigma\bar\xi_2)
\otimes P(\sigma\bar\tau_2)
\endaligned
\tag 4.1
$$
as $A_*$-comodule algebras.  The $A_*$-comodule coactions are given by
$$
\aligned
\nu(\sigma\bar\tau_0) &= 1 \otimes \sigma\bar\tau_0 \\
\nu(\sigma\bar\xi_1) &= 1 \otimes \sigma\bar\xi_1 \\
\nu(\sigma\bar\tau_1) &= 1 \otimes \sigma\bar\tau_1  + \bar\tau_0
\otimes \sigma\bar\xi_1 \\
\nu(\sigma\bar\xi_2) &= 1 \otimes \sigma\bar\xi_2 \\
\nu(\sigma\bar\tau_2) &= 1 \otimes \sigma\bar\tau_2  + \bar\tau_0
\otimes \sigma\bar\xi_2 \,.
\endaligned
\tag 4.2
$$
The map $\pi^* \: THH(\ell) \to THH(\Z_{(p)})$ takes $\sigma\bar\xi_2$
to $0$ and $\sigma\bar\tau_2$ to $(\sigma\bar\tau_1)^p$.  The map $i^*
\: THH(\Z_{(p)}) \to THH(\Z/p)$ takes $\sigma\bar\xi_1$ to $0$ and
$\sigma\bar\tau_1$ to $(\sigma\bar\tau_0)^p$.

The B{\"o}kstedt spectral sequence for the associative $S$-algebra
$B = \ell/p$ begins
$$
E^2_{**}(\ell/p) = H_*(\ell/p) \otimes E(\sigma\bar\xi_k \mid k\ge1)
\otimes \Gamma(\sigma\bar\tau_0, \sigma\bar\tau_k \mid k\ge2) \,.
$$
It is an $A_*$-comodule module spectral sequence over the B{\"o}kstedt
spectral sequence for $\ell$, since the $\ell$-algebra multiplication
$\ell \wedge \ell/p \to \ell/p$ is a map of associative $S$-algebras.
However, it is not itself an algebra spectral sequence, since the
product on $\ell/p$ is not commutative enough to induce a natural
product structure on $THH(\ell/p)$.  Nonetheless, we will use the
algebra structure present at the $E^2$-term to help in naming classes.

The map $\pi \: \ell/p \to H\Z/p$ induces an injection of B{\"o}kstedt
spectral sequence $E^2$-terms, so there are differentials generated
algebraically by
$$
d^{p-1}(\gamma_j\sigma\bar\tau_k) = \sigma\bar\xi_{k+1}
\cdot \gamma_{j-p}\sigma\bar\tau_k
$$
for $j\ge p$, $k=0$ or $k\ge2$, leaving
$$
E^\infty_{**}(\ell/p) = H_*(\ell/p) \otimes E(\sigma\bar\xi_2) \otimes
P_p(\sigma\bar\tau_0, \sigma\bar\tau_k \mid k\ge2)
\tag 4.3
$$
as an $A_*$-comodule module over $E^\infty_{**}(\ell)$.  We need to
resolve the $A_*$-comodule and $H_*(THH(\ell))$-module extensions in
order to obtain $H_*(THH(\ell/p))$.  This is achieved in Lemma~4.6 below.

The map $\pi^* \: E^\infty_{**}(\ell/p) \to E^\infty_{**}(\Z/p)$ is an
isomorphism in total degrees $\le (2p-2)$ and injective in total degrees
$\le (2p^2-2)$.  The first class in the kernel is $\sigma\bar\xi_2$.
Hence there are unique classes
$$
1 \ ,\  \bar\tau_0 \ ,\  \sigma\bar\tau_0 \ ,\  \bar\tau_0
\sigma\bar\tau_0 \ ,\  \dots \ ,\ (\sigma\bar\tau_0)^{p-1}
$$
in degrees $0 \le * \le 2p-2$ of $H_*(THH(\ell/p))$, mapping to classes
with the same names in $H_*(THH(\Z/p))$.  More concisely, these are the
monomials $\bar\tau_0^\delta (\sigma\bar\tau_0)^i$ for $0\le\delta\le1$
and $0\le i\le p-1$, except that the degree~$(2p-1)$ case $(\delta,i) =
(1,p-1)$ is omitted.  The $A_*$-comodule coaction on these classes is
given by the same formulas in $H_*(THH(\ell/p))$ as in $H_*(THH(\Z/p))$,
cf.~(4.2).

There is also a class $\bar\xi_1$ in degree~$(2p-2)$ of
$H_*(THH(\ell/p))$
mapping to a class with the same name, and same $A_*$-coaction, in
$H_*(THH(\Z/p))$.

In degree~$(2p-1)$ there is a map $\pi^*$ of extensions from
$$
0 \to \F_p\{\bar\xi_1 \bar\tau_0\} \to H_{2p-1}(THH(\ell/p)) \to
\F_p\{\bar\tau_0 (\sigma\bar\tau_0)^{p-1}\} \to 0
$$
to
$$
0 \to \F_p\{\bar\tau_1, \bar\xi_1 \bar\tau_0\} \to H_{2p-1}(THH(\Z/p))
\to \F_p\{\bar\tau_0 (\sigma\bar\tau_0)^{p-1}\} \to 0 \,.
$$
The latter extension is canonically split by the augmentation
$\epsilon \: THH(\Z/p) \to H\Z/p$, which uses the commutativity
of the $S$-algebra $H\Z/p$.

In degree~$2p$, the map $\pi^*$ goes from
$$
H_{2p}(THH(\ell/p)) = \F_p\{\bar\xi_1 \sigma\bar\tau_0\}
$$
to
$$
0 \to \F_p\{\bar\tau_0 \bar\tau_1\} \to H_{2p}(THH(\Z/p)) \to
\F_p\{\sigma\bar\tau_1, \bar\xi_1 \sigma\bar\tau_0\} \to 0 \,.
$$
Again the latter extension is canonically split.

\proclaim{Lemma 4.4}
There is a unique class $y$ in $H_{2p-1}(THH(\ell/p))$ that is represented by
$\bar\tau_0 (\sigma\bar\tau_0)^{p-1}$ in $E^\infty_{p-1,p}(\ell/p)$
and that maps to $\bar\tau_0 (\sigma\bar\tau_0)^{p-1} -
\bar\tau_1$ in $H_*(THH(\Z/p))$.
\endproclaim

\demo{Proof}
This follows by naturality of the suspension operator $\sigma$ and
the multiplicative relation $(\sigma\bar\tau_0)^p = \sigma\bar\tau_1$
in $H_*(THH(\Z/p))$.  A class $y$ in $H_{2p-1}(THH(\ell/p))$ represented
by $\bar\tau_0 (\sigma\bar\tau_0)^{p-1}$ is determined modulo $\bar\xi_1
\bar\tau_0$.  Its image in $H_{2p-1}(THH(\Z/p))$ thus has the form $\alpha
\bar\tau_1 + \bar\tau_0 (\sigma\bar\tau_0)^{p-1}$ modulo $\bar\xi_1
\bar\tau_0$, for some $\alpha \in \F_p$.  The suspension $\sigma y$
lies in $H_{2p}(THH(\ell/p)) = \F_p\{\bar\xi_1 \sigma\bar\tau_0\}$, so
its image in $H_{2p}(THH(\Z/p))$ is $0$ modulo $\bar\tau_0 \bar\tau_1$
and $\bar\xi_1 \sigma\bar\tau_0$.  It is also the suspension of
$\alpha \bar\tau_1 + \bar\tau_0 (\sigma\bar\tau_0)^{p-1}$ modulo
$\bar\xi_1 \bar\tau_0$, which equals $\sigma(\alpha \bar\tau_1) +
(\sigma\bar\tau_0)^p = (\alpha+1) \sigma\bar\tau_1$.  In particular, the
coefficient $(\alpha + 1)$ of $\sigma\bar\tau_1$ is $0$, so $\alpha=-1$.
\qed
\enddemo

\remark{Remark 4.5}
For $p=2$ this can alternatively be read off from the explicit form
\cite{Wu91} of the commutator for the product $\mu$ in $\ell/p$.
The coequalizer $C$ of the two maps
$$
\xymatrix{
\ell/p \wedge \ell/p \ar@<0.6ex>[r]^-{\mu} \ar@<-0.6ex>[r]_-{\mu\tau} & \ell/p
}
$$
maps to (the $1$-skeleton of) $THH(\ell/p)$.  The commutator $\mu -
\mu\tau$ factors as
$$
\ell/p \wedge \ell/p @>\beta\wedge\beta>> \Sigma \ell/p \wedge \Sigma \ell/p
@>\mu>> \Sigma^2 \ell/p @>v_1>> \ell/p
$$
where $\beta$ is the mod~$p$ Bockstein associated to the cofiber
sequence $\ell @>p>> \ell @>i>> \ell/p$ and the cofiber of $v_1$ is
$H\Z/p$.  We get a map of cofiber sequences
$$
\xymatrix{
\ell/p \wedge \ell/p \ar[r]^-{\mu - \mu\tau} \ar[d]_{\mu(\beta \wedge
	\beta)} & \ell/p \ar[r] \ar@{=}[d] & C \ar[d] \\
\Sigma^2 \ell/p \ar[r]^-{v_1} & \ell/p \ar[r] & H\Z/p \,,
}
$$
so there is a class in $H_3(C)$ that maps to $\bar\xi_1
\otimes \bar\xi_1$ in $H_2(\ell/p \wedge \ell/p)$ and to $\bar\xi_1
\sigma\bar\xi_1$ in $H_3(THH(\ell/p))$, which also maps to $\bar\xi_2$
in the cofiber of $v_1$, i.e., whose $A_*$-coaction contains the term
$\bar\xi_2 \otimes 1$.  (The classes $\bar\tau_0$ and $\bar\tau_1$ go
by the names $\bar\xi_1$ and $\bar\xi_2$ at $p=2$.)

For odd primes there is a similar interpretation of how the
non-com\-mutativity of the product on $\ell/p$ provides an obstruction
to splitting off the $0$-simplices from the $(p-1)$-skeleton of
$THH(\ell/p)$, where the cyclic permutation of the $p$ factors in the
$(p-1)$-simplex $\bar\tau_0 (\sigma\bar\tau_0)^{p-1}$, represented by
the Hochschild cycle $\bar\tau_0 \otimes \dots \otimes \bar\tau_0$,
plays a similar role to the twist map $\tau$ above.
\endremark

\medskip

Let
$$
H_*(THH(\ell))/(\sigma\bar\xi_1) \cong
  H_*(\ell) \otimes E(\sigma\bar\xi_2) \otimes P(\sigma\bar\tau_2)
$$
denote the quotient algebra of $H_*(THH(\ell))$ by the ideal generated
by $\sigma\bar\xi_1$.

\proclaim{Lemma 4.6}
There is an isomorphism of $H_*(THH(\ell))$-modules
$$
H_*(THH(\ell/p)) \cong H_*(THH(\ell))/(\sigma\bar\xi_1)
  \otimes \F_p\{1, \bar\tau_0, \sigma\bar\tau_0, \bar\tau_0
  \sigma\bar\tau_0, \dots, (\sigma\bar\tau_0)^{p-1}, y\} \,.
$$
Hence $H_*(THH(\ell/p))$ is a free module of rank~$2p$ over
$H_*(THH(\ell))/(\sigma\bar\xi_1)$, generated by classes
$$
1 \ ,\  \bar\tau_0 \ ,\  \sigma\bar\tau_0 \ ,\  \bar\tau_0
\sigma\bar\tau_0 \ ,\  \dots \ ,\ (\sigma\bar\tau_0)^{p-1} \ ,\  y
$$
in degrees $0$ through~$2p-1$.  These generators are represented in
$E^\infty_{**}(\ell/p)$ by the classes
$$
1 \ ,\  \bar\tau_0 \ ,\  \sigma\bar\tau_0 \ ,\  \bar\tau_0
\sigma\bar\tau_0 \ ,\  \dots \ ,\ (\sigma\bar\tau_0)^{p-1} \ ,\
\bar\tau_0 (\sigma\bar\tau_0)^{p-1} \,,
$$
and map under $\pi^*$ to classes with the same names in $H_*(THH(\Z/p))$,
except for~$y$, which maps to
$$
\bar\tau_0 (\sigma\bar\tau_0)^{p-1} - \bar\tau_1 \,.
$$
The $A_*$-comodule coactions are given by
$$
\nu((\sigma\bar\tau_0)^i) = 1 \otimes (\sigma\bar\tau_0)^i
$$
for $0\le i\le p-1$,
$$
\nu(\bar\tau_0 (\sigma\bar\tau_0)^i) = 1 \otimes \bar\tau_0
  (\sigma\bar\tau_0)^i + \bar\tau_0 \otimes (\sigma\bar\tau_0)^i
$$
for $0 \le i \le p-2$, and
$$
\nu(y) = 1 \otimes y + \bar\tau_0 \otimes (\sigma\bar\tau_0)^{p-1} -
  \bar\tau_0 \otimes \bar\xi_1 - \bar\tau_1 \otimes 1 \,.
$$
\endproclaim

\demo{Proof}
$H_*(\ell/p)$ is freely generated as a module over $H_*(\ell)$ by $1$
and $\bar\tau_0$, and the classes $\sigma\bar\xi_2$ and
$\sigma\bar\tau_2$ in $H_*(THH(\ell))$ induce multiplication by the
same symbols in $E^\infty_{**}(\ell/p)$, as given in~(4.3).  This
generates all of $E^\infty_{**}(\ell/p)$ from the $2p$ classes
$\bar\tau_0^\delta (\sigma\bar\tau_0)^i$ for $0\le\delta\le1$
and $0\le i\le p-1$.

We claim that multiplication by $\sigma\bar\xi_1$ acts trivially on
$H_*(THH(\ell/p))$.  It suffices to verify this on the module
generators $\bar\tau_0^\delta (\sigma\bar\tau_0)^i$, for which the
product with $\sigma\bar\xi_1$ remains in the range of degrees where
the map to $H_*(THH(\Z/p))$ is injective.  The action of
$\sigma\bar\xi_1$ is trivial on $H_*(THH(\Z/p))$, since
$d^{p-1}(\gamma_p \sigma\bar\tau_0) = \sigma\bar\xi_1$ and
$\epsilon(\sigma\bar\xi_1) = 0$, from which the claim follows.

The $A_*$-comodule coaction on each module generator, including $y$, is
determined by that on its image under $\pi^*$.
In the latter case, the thing to check is that
$$
\align
(1 \otimes \pi^*)(\nu(y)) &= \nu(\pi^*(y)) =
	\nu(\bar\tau_0(\sigma\bar\tau_0)^{p-1}-\bar\tau_1) \\
&= 1 \otimes \bar\tau_0(\sigma\bar\tau_0)^{p-1} + \bar\tau_0 \otimes
	(\sigma\bar\tau_0)^{p-1} - 1 \otimes \bar\tau_1 - \bar\tau_0
	\otimes \bar\xi_1 - \bar\tau_1 \otimes 1
\endalign
$$
equals
$$
(1 \otimes \pi^*)(1 \otimes y + \bar\tau_0 \otimes
	(\sigma\bar\tau_0)^{p-1} - \bar\tau_0 \otimes \bar\xi_1 -
	\bar\tau_1 \otimes 1) \,.
$$
\qed
\enddemo

We note that these results do not visibly depend on the particular
choice of $\ell$-algebra structure on $\ell/p$.

\head 5. Passage to $V(1)$-homotopy \endhead

For $p\ge5$ the Smith--Toda complex $V(1)$ is a homotopy commutative
ring spectrum.  In this section we compute $V(1)_* THH(\ell/p)$ as a
module over $V(1)_* THH(\ell)$, for any prime $p\ge5$.  The unique ring
spectrum map from $V(1)$ to $H\Z/p$ induces the identification
$$
H_*(V(1)) = E(\tau_0, \tau_1)
$$
(no conjugations) as $A_*$-comodule subalgebras of $A_*$.  Here
$$
\align
\nu(\tau_0) &= 1 \otimes \tau_0 + \tau_0 \otimes 1 \\
\nu(\tau_1) &= 1 \otimes \tau_1 + \xi_1 \otimes \tau_0 + \tau_1 \otimes 1
\,.
\endalign
$$
For each $\ell$-algebra $B$, $V(1) \wedge THH(B)$ is a module spectrum
over $V(1) \wedge THH(\ell)$ and thus over $V(1) \wedge \ell \simeq
H\Z/p$, so $H_*(V(1) \wedge THH(B))$ is a sum of copies of $A_*$ as an
$A_*$-comodule.  In particular, $V(1)_* THH(B) = \pi_*(V(1) \wedge
THH(B))$ is identified with the subgroup of $A_*$-comodule primitives
in
$$
H_*(V(1) \wedge THH(B)) \cong H_*(V(1)) \otimes H_*(THH(B))
$$
with the diagonal $A_*$-comodule coaction.  We write $v \wedge x$ for
the image of $v \otimes x$ under this identification, with $v \in
H_*(V(1))$ and $x \in H_*(THH(B))$.  Let
$$
\aligned
\epsilon_0 &= 1 \wedge \bar\tau_0 + \tau_0 \wedge 1 \\
\epsilon_1 &= 1 \wedge \bar\tau_1 + \tau_0 \wedge \bar\xi_1 + \tau_1
\wedge 1 \\
\lambda_1 &= 1 \wedge \sigma\bar\xi_1 \\
\lambda_2 &= 1 \wedge \sigma\bar\xi_2 \\
\mu_0 &= 1 \wedge \sigma\bar\tau_0 \\
\mu_1 &= 1 \wedge \sigma\bar\tau_1 + \tau_0 \wedge \sigma\bar\xi_1 \\
\mu_2 &= 1 \wedge \sigma\bar\tau_2 + \tau_0 \wedge \sigma\bar\xi_2 \,.
\endaligned
\tag 5.1
$$
These are all $A_*$-comodule primitive, where defined.  By a dimension
count,
$$
\aligned
V(1)_* THH(\Z/p) &= E(\epsilon_0, \epsilon_1) \otimes P(\mu_0) \\
V(1)_* THH(\Z_{(p)}) &= E(\epsilon_1) \otimes E(\lambda_1) \otimes P(\mu_1) \\
V(1)_* THH(\ell) &= E(\lambda_1, \lambda_2) \otimes P(\mu_2)
\endaligned
\tag 5.2
$$
as commutative $\F_p$-algebras.  The map $\pi \: \ell \to H\Z_{(p)}$
takes $\lambda_2$ to $0$ and $\mu_2$ to $\mu_1^p$.  The map $i \:
H\Z_{(p)} \to H\Z/p$ takes $\lambda_1$ to $0$ and $\mu_1$ to
$\mu_0^p$.  Note that $\mu_2 \in V(1)_{2p^2} THH(\ell)$ was
simply denoted $\mu$ in \cite{AuR02}.

In degrees $\le (2p-2)$ of $H_*(V(1) \wedge THH(\ell/p))$ the classes
$$
\mu_0^i := 1 \wedge (\sigma\bar\tau_0)^i
\tag 5.3.a
$$
for $0\le i\le p-1$ and
$$
\epsilon_0 \mu_0^i := 1 \wedge \bar\tau_0 (\sigma\bar\tau_0)^i + \tau_0
\wedge (\sigma\bar\tau_0)^i
\tag 5.3.b
$$
for $0\le i\le p-2$ are $A_*$-comodule primitive.  These map to the
classes $\epsilon_0^\delta \mu_0^i$ in $V(1)_* THH(\Z/p)$ for
$0\le\delta\le1$ and $0\le i\le p-1$, except that the degree bound
excludes the top case of $\epsilon_0\mu_0^{p-1}$.

In degree $(2p-1)$ of $H_*(V(1) \wedge THH(\ell/p))$ we have generators
$1 \wedge \bar\xi_1 \bar\tau_0$, $\tau_0 \wedge
(\sigma\bar\tau_0)^{p-1}$, $\tau_0 \wedge \bar\xi_1$, $\tau_1 \wedge 1$
and $1 \wedge y$.  These have coactions
$$
\align
\nu(1 \wedge \bar\xi_1 \bar\tau_0) &= 1 \otimes 1 \wedge \bar\xi_1
\bar\tau_0 + \bar\tau_0 \otimes 1 \wedge \bar\xi_1 + \bar\xi_1 \otimes
1 \wedge \bar\tau_0 + \bar\xi_1 \bar\tau_0 \otimes 1 \wedge 1 \\
\nu(\tau_0 \wedge (\sigma\bar\tau_0)^{p-1}) &= 1 \otimes \tau_0 \wedge 
(\sigma\bar\tau_0)^{p-1} + \tau_0 \otimes 1 \wedge (\sigma\bar\tau_0)^{p-1} \\
\nu(\tau_0 \wedge \bar\xi_1) &= 1 \otimes \tau_0 \wedge \bar\xi_1 +
\tau_0 \otimes 1 \wedge \bar\xi_1 + \bar\xi_1 \otimes \tau_0 \wedge 1
+ \bar\xi_1\tau_0 \otimes 1 \wedge 1 \\
\nu(\tau_1 \wedge 1) &= 1 \otimes \tau_1 \wedge 1 + \xi_1 \otimes \tau_0
\wedge 1 + \tau_1 \otimes 1 \wedge 1
\endalign
$$
and
$$
\nu(1 \wedge y) = 1 \otimes 1 \wedge y + \bar\tau_0 \otimes 1 \wedge
	(\sigma\bar\tau_0)^{p-1} - \bar\tau_0 \otimes 1 \wedge \bar\xi_1 -
	\bar\tau_1 \otimes 1 \wedge 1 \,.
$$
Hence the sum
$$
\bar\epsilon_1 := 1 \wedge y + \tau_0 \wedge
(\sigma\bar\tau_0)^{p-1} - \tau_0 \wedge \bar\xi_1 - \tau_1 \wedge 1
\tag 5.3.c
$$
is $A_*$-comodule primitive.  Its image under $\pi^*$ in $H_*(V(1)
\wedge THH(\Z/p))$
is
$$
\epsilon_0 \mu_0^{p-1} - \epsilon_1 = 1 \wedge \bar\tau_0
(\sigma\bar\tau_0)^{p-1} + \tau_0 \wedge (\sigma\bar\tau_0)^{p-1} - 1
\wedge \bar\tau_1 - \tau_0 \wedge \bar\xi_1 - \tau_1 \wedge 1 \,.
$$

Let
$$
V(1)_* THH(\ell)/(\lambda_1) \cong E(\lambda_2) \otimes P(\mu_2)
$$
be the quotient algebra of $V(1)_* THH(\ell)$ by the ideal generated
by $\lambda_1$.

\proclaim{Proposition 5.4}
There is an isomorphism of $V(1)_* THH(\ell)$-modules
$$
V(1)_* THH(\ell/p) = V(1)_* THH(\ell)/(\lambda_1) \otimes \F_p\{1,
\epsilon_0, \mu_0, \epsilon_0 \mu_0, \dots, \mu_0^{p-1}, \bar\epsilon_1\}
\,,
$$
where the classes $\mu_0^i$, $\epsilon_0 \mu_0^i$ and $\bar\epsilon_1$
are defined in~(5.3.abc) above.  Multiplication by $\lambda_1$ is $0$, so
this is a free module on the $2p$ generators
$$
1 \ ,\  \epsilon_0 \ ,\  \mu_0 \ ,\  \epsilon_0 \mu_0 \ ,\  \dots \ ,\
\mu_0^{p-1} \ ,\  \bar\epsilon_1
$$
over $V(1)_* THH(\ell)/(\lambda_1)$.  The map $\pi^*$ to $V(1)_* THH(\Z/p)$
takes $\epsilon_0^\delta \mu_0^i$ in degree $0 \le \delta + 2i \le 2p-2$
to $\epsilon_0^\delta \mu_0^i$, and takes $\bar\epsilon_1$ in degree
$(2p-1)$ to $\epsilon_0 \mu_0^{p-1} - \epsilon_1$.
\endproclaim

\demo{Proof}
Additively, this follows by another dimension count.  The
multiplication by $\lambda_1$ is $0$ for degree and filtration
reasons:  $\lambda_1$ has B{\"o}kstedt filtration~$1$ and cannot map to
$\bar\epsilon_1$ in B{\"o}kstedt filtration~$(p-1)$.  Similarly in
higher degrees.
\qed
\enddemo

We end this section with a suggestive conjectural calculation of the
topological Hochschild homology of the fraction field $\ff(\ell)
= p^{-1}L$, which may play the role of the de\,Rham complex over
$\Spec(\ff(\ell))$ in derived algebraic geometry.  The calculation is
not needed for the rest of the paper.  We work with the $p$-local $\ell$,
but could equally well have worked with the $p$-complete $\ell_p$.

Thus consider a $3 \times 3$ square of cofiber sequences
$$
\xymatrix{
THH(\Z/p) \ar[r]^{i_*} \ar[d]^{\pi_*} &
THH(\Z_{(p)}) \ar[r]^{j^*} \ar[d]^{\pi_*} &
THH(\Z_{(p)}|\Q) \ar[d]^{\pi_*} \\
THH(\ell/p) \ar[r]^{i_*} \ar[d]^{\rho^*} &
THH(\ell) \ar[r]^{j^*} \ar[d]^{\rho^*} &
THH(\ell|p^{-1}\ell) \ar[d]^{\rho^*} \\
THH(\ell/p|L/p) \ar[r]^{i_*} &
THH(\ell|L) \ar[r]^{j^*} &
THH(\ff(\ell))
}
\tag 5.5
$$
generated by the upper left hand square.  The transfer map $i_* \:
THH(\Z/p) \to THH(\Z_{(p)})$ was properly defined in \cite{HM03},
and $THH(\Z_{(p)}|\Q)$ is its cofiber.  To construct the remaining
transfer maps one may use the definition in \cite{BM:loc} of $THH$
in terms of spectral categories.  By analogy with the algebraic case,
we write $THH(\ell|p^{-1}\ell)$, $THH(\ell/p|L/p)$ and $THH(\ell|L)$
for three of the remaining cofibers in the diagram.  We might have
denoted the term in the lower right hand corner by something like
$THH(\ell|p^{-1}\ell,L|p^{-1}L)$, but for simplicity we abbreviate this
to $THH(\ff(\ell))$.

In the top row we get \cite{HM03, 2.4.1} a long exact sequence
in $V(1)$-homotopy
$$
\dots @>\partial>>
E(\epsilon_0, \epsilon_1) \otimes P(\mu_0)
@>i_*>> E(\epsilon_1, \lambda_1) \otimes P(\mu_1)
@>j^*>> E(\dlog p, \epsilon_1) \otimes P(\kappa_0)
@>\partial>> \dots
$$
where $i_*(\epsilon_0\mu_0^{p-1}) = \lambda_1$, $j^*(\mu_1) =
\kappa_0^p$, $\partial(\dlog p) = 1$ and $\partial(\kappa_0) =
\epsilon_0$, tensored with the identity on $E(\epsilon_1)$.

In the middle row we expect a long exact sequence
$$
\multline
\dots @>\partial>>
\F_p\{1, \epsilon_0, \mu_0, \epsilon_0 \mu_0, \dots, \mu_0^{p-1},
  \bar\epsilon_1\} \otimes E(\lambda_2) \otimes P(\mu_2) @>i_*>> \\
@>i_*>> E(\lambda_1, \lambda_2) \otimes P(\mu_2)
@>j^*>> E(\dlog p, \lambda_2) \otimes P_p(\kappa_0) \otimes P(\mu_2)
@>\partial>> \dots
\endmultline
$$
where $i_*(\bar\epsilon_1) = \lambda_1$, $\partial(\dlog p) = 1$ and
$\partial(\kappa_0) = \epsilon_0$, tensored with the identity on
$E(\lambda_2)$.

In the middle column we expect \cite{Au05, \S10} a long exact
sequence
$$
\dots @>\partial>>
E(\epsilon_1, \lambda_1) \otimes P(\mu_1)
@>\pi_*>> E(\lambda_1, \lambda_2) \otimes P(\mu_2)
@>\rho^*>> E(\dlog v_1, \lambda_1) \otimes P(\kappa_1)
@>\partial>> \dots
$$
where $\pi_*(\epsilon_1\mu_1^{p-1}) = \lambda_2$, $\rho^*(\mu_2) =
\kappa_1^p$, $\partial(\dlog v_1) = 1$ and $\partial(\kappa_1) =
\epsilon_1$, tensored with the identity on $E(\lambda_1)$.

In the right hand column we expect a long exact sequence
$$
\multline
\dots @>\partial>>
E(\dlog p, \epsilon_1) \otimes P(\kappa_0) @>\pi_*>> \\
@>\pi_*>> E(\dlog p, \lambda_2) \otimes P_p(\kappa_0) \otimes P(\mu_2)
@>\rho^*>> E(\dlog p, \dlog v_1) \otimes P(\kappa_0)
@>\partial>> \dots
\endmultline
$$
where $\pi_*(\epsilon_1\kappa_0^{p^2-p}) = \lambda_2$, $\rho^*(\mu_2) =
\kappa_0^{p^2}$, $\partial(\dlog v_1) = 1$ and $\partial(\kappa_0^p) =
\epsilon_1$, tensored with the identity on $E(\dlog p)$.  Note that this
$\rho^*$ is not multiplicative, since it takes the truncated polynomial
algebra $P_p(\kappa_0)$ into $P(\kappa_0)$.

This leads to the following formula.  When compared with \cite{HM03,
2.4.1}, it strongly suggests that $E(\dlog p, \dlog v_1)$ is the
$V(1)$-homotopy of a de\,Rham complex for $\ell$ with logarithmic
poles along $(p)$ and $(v_1)$.

\proclaim{Conjecture 5.6}
There is an isomorphism
$$
V(1)_* THH(\ff(\ell)) \cong E(\dlog p, \dlog v_1) \otimes P(\kappa_0)
$$
with $|\dlog p\,| = |\dlog v_1| = 1$ and $|\kappa_0| = 2$.
Here $THH(\ff(\ell))$
is defined as the iterated cofiber of the upper left hand square in
diagram~(5.5), $\dlog p$ is in the image from $\pi_1 THH(\ell|p^{-1}\ell)$,
with $\partial(\dlog p) = 1$ in $\pi_0 THH(\ell/p)$, $\dlog v_1$
is in the image from $\pi_1 THH(\ell|L)$, with $\partial(\dlog v_1) = 1$
in $\pi_0 THH(\Z_{(p)})$, and $\kappa_0$ satisfies $\kappa_0^{p^2} =
\mu_2$, with $\mu_2$ in the image from $V(1)_{2p^2} THH(\ell)$.
\endproclaim

\head 6. The $C_p$-Tate construction \endhead

Let $C = C_{p^n}$ denote the cyclic group of order $p^n$, considered as a
closed subgroup of the circle group $S^1$.  For each spectrum $X$
with $C$-action, $X_{hC} = EC_+ \wedge_C X$ and $X^{hC} = F(EC_+, X)^C$
denote its homotopy orbit and homotopy fixed point spectra, as usual.
We now write $X^{tC} = [\widetilde{EC} \wedge F(EC_+, X)]^C$ for the
$C$-Tate construction on~$X$, denoted $t_C(X)^C$ in \cite{GM95} and
$\hat\H(C, X)$ in \cite{AuR02}.  There are $C$-homotopy fixed point
and $C$-Tate spectral sequences in $V(1)$-homotopy for $X$, with
$$
E^2_{s,t}(C, X) = H_{gp}^{-s}(C; V(1)_t(X))
\Longrightarrow V(1)_{s+t}(X^{hC})
$$
and
$$
\hat E^2_{s,t}(C, X) = \hat H_{gp}^{-s}(C; V(1)_t(X))
\Longrightarrow V(1)_{s+t}(X^{tC}) \,.
$$
We write $H_{gp}^*(C_{p^n}; \F_p) = E(u_n) \otimes P(t)$ and $\hat
H_{gp}^*(C_{p^n}; \F_p) = E(u_n) \otimes P(t^{\pm1})$ with $u_n$ in degree
$1$ and $t$ in degree $2$.  So $u_n$, $t$ and $x \in V(1)_t(X)$ have
bidegree $(-1,0)$, $(-2,0)$ and $(0,t)$ in either spectral sequence,
respectively.  See \cite{HM03, \S4.3} for proofs of the multiplicative
properties of these spectral sequences.

We are principally interested in the case when $X = THH(B)$, with the
$S^1$-action given by the cyclic structure.  It is a cyclotomic
spectrum, in the sense of \cite{HM97}, leading to the
commutative diagram
$$
\xymatrix{
THH(B)_{hC_{p^n}} \ar[r]^N \ar@{=}[d] &
THH(B)^{C_{p^n}} \ar[r]^R \ar[d]^{\Gamma_n} &
THH(B)^{C_{p^{n-1}}} \ar[d]^{\hat\Gamma_n} \\
THH(B)_{hC_{p^n}} \ar[r]^{N^h} &
THH(B)^{hC_{p^n}} \ar[r]^{R^h} &
THH(B)^{tC_{p^n}}
}
$$
of horizontal cofiber sequences.  We abbreviate $\hat E^2_{**}(C,
THH(B))$ to $\hat E^2_{**}(C, B)$, etc.  When $B$ is a commutative
$S$-algebra, this is a commutative algebra spectral sequence, and when
$B$ is an associative $A$-algebra, with $A$ commutative, then $\hat
E^*(C, B)$ is a module spectral sequence over $\hat E^*(C,
A)$.  The map $R^h$ corresponds to the inclusion $E^2_{**}(C, B) \to
\hat E^2_{**}(C, B)$ from the second quadrant to the upper half-plane,
for connective $B$.

In this section we compute $V(1)_* THH(\ell/p)^{tC_p}$ by means of the
$C_p$-Tate spectral sequence in $V(1)$-homotopy for $THH(\ell/p)$.  In
Propositions~6.8 and~6.9 we show that the comparison map $\hat\Gamma_1
\: V(1)_* THH(\ell/p) \to V(1)_* THH(\ell/p)^{tC_p}$ is
$(2p-2)$-coconnected and can be identified with the algebraic
localization homomorphism that inverts $\mu_2$.

First we recall the structure of the $C_p$-Tate spectral sequence for
$THH(\Z/p)$, with $V(0)$- and $V(1)$-coefficients.  We have $V(0)_*
THH(\Z/p) = E(\epsilon_0) \otimes P(\mu_0)$, and with an obvious
notation the $E^2$-terms are
$$
\align
\hat E^2_{**}(C_p, \Z/p; V(0)) &= E(u_1) \otimes P(t^{\pm1}) \otimes
E(\epsilon_0) \otimes P(\mu_0) \\
\hat E^2_{**}(C_p, \Z/p) &=  E(u_1) \otimes P(t^{\pm1}) \otimes
E(\epsilon_0, \epsilon_1) \otimes P(\mu_0) \,.
\endalign
$$
In each $C$-Tate spectral sequence we have a first differential
$$
d^2(x) = t \cdot \sigma x \,,
$$
see e.g.~\cite{Rog98, 3.3}.  We easily deduce $\sigma\epsilon_0 = \mu_0$
and $\sigma\epsilon_1 = \mu_0^p$ from~(5.1), so
$$
\align
\hat E^3_{**}(C_p, \Z/p; V(0)) &= E(u_1) \otimes P(t^{\pm1}) \\
\hat E^3_{**}(C_p, \Z/p) &=  E(u_1) \otimes P(t^{\pm1}) \otimes
E(\epsilon_0 \mu_0^{p-1} - \epsilon_1) \,.
\endalign
$$
Thus the $V(0)$-homotopy spectral sequence collapses at $\hat E^3 =
\hat E^\infty$.  By naturality with respect to the map $i_1 \: V(0) \to
V(1)$, all the classes on the horizontal axis of $\hat E^3(C_p, \Z/p)$
are infinite cycles, so also the latter spectral sequence collapses at
$\hat E^3_{**}(C_p, \Z/p) = \hat E^\infty_{**}(C_p, \Z/p)$.

We know from \cite{HM97, Prop.~5.3} that the comparison map 
$$
\hat\Gamma_1 \:
V(0)_* THH(\Z/p) \to V(0)_* THH(\Z/p)^{tC_p}
$$ 
takes $\epsilon_0^\delta
\mu_0^i$ to $(u_1 t^{-1})^\delta t^{-i}$, for all $0\le\delta\le1$,
$i\ge0$.  In particular, the integral map $\hat\Gamma_1 \: \pi_*
THH(\Z/p) \to \pi_* THH(\Z/p)^{tC_p}$ is $(-2)$-coconnected, meaning
that it induces an injection in degree~$(-2)$ and an isomorphism in all
higher degrees.  From this we can deduce the following behavior of the
comparison map $\hat\Gamma_1$ in $V(1)$-homotopy.

\proclaim{Lemma 6.1}
The map
$$
\hat\Gamma_1 \: V(1)_* THH(\Z/p) \to V(1)_* THH(\Z/p)^{tC_p}
$$
takes the classes $\epsilon_0^\delta \mu_0^i$ from $V(0)_* THH(\Z/p)$,
for $0\le\delta\le1$ and $i\ge0$, to classes represented in $\hat
E^\infty_{**}(C_p, \Z/p)$ by $(u_1 t^{-1})^\delta t^{-i}$ (on the
horizontal axis).

Furthermore, it takes the class $\epsilon_0 \mu_0^{p-1} - \epsilon_1$
in degree $(2p-1)$ to a class represented by $\epsilon_0 \mu_0^{p-1} -
\epsilon_1$ (on the vertical axis).
\endproclaim

\demo{Proof}
The classes $\epsilon_0^\delta \mu_0^i$ are in the image from
$V(0)$-homotopy, and we recalled above that they are detected by
$(u_1 t^{-1})^\delta t^{-i}$ in the $V(0)$-homotopy $C_p$-Tate spectral
sequence for $THH(\Z/p)$.  By naturality along $i_1 \: V(0) \to V(1)$,
they are detected by the same (nonzero) classes in the $V(1)$-homotopy
spectral sequence $\hat E^\infty_{**}(C_p, \Z/p)$.

To find the representative for $\hat\Gamma_1(\epsilon_0 \mu_0^{p-1} -
\epsilon_1)$ in degree~$(2p-1)$, we appeal to the cyclotomic trace map
from algebraic $K$-theory, or more precisely, to the commutative
diagram
$$
\xymatrix{
& K(B) \ar[dr]^{tr} \ar[d]^{tr_1} \ar[dl]_{tr} \\
THH(B) & THH(B)^{C_p} \ar[r]^R \ar[d]^{\Gamma_1} \ar[l]_F & THH(B)
  \ar[d]^{\hat\Gamma_1} \\
& THH(B)^{hC_p} \ar[r]^{R^h} \ar[ul] & THH(B)^{tC_p} \,.
}
\tag 6.2
$$
The B{\"o}kstedt trace map $tr \: K(B) \to THH(B)$ admits a preferred
lift $tr_n$ through each fixed-point spectrum $THH(B)^{C_{p^n}}$, which
equalizes the iterated restriction and Frobenius maps $R^n$ and $F^n$
to $THH(B)$ \cite{BHM93}.  In particular, the circle action and the
$\sigma$-operator act trivially on classes in the image of $tr$.

In the case $B = H\Z/p$, we know that $K(\Z/p)_p \simeq H\Z_p$ after
$p$-adic completion, so $V(1)_* K(\Z/p) = E(\bar\epsilon_1)$, where
the $v_1$-Bockstein of $\bar\epsilon_1$ is $-1$.  The B{\"o}kstedt trace
image $tr(\bar\epsilon_1) \in V(1)_* THH(\Z/p)$ lies in $\F_p\{\epsilon_1,
\epsilon_0 \mu_0^{p-1}\}$, has $v_1$-Bockstein $tr(-1) = -1$ and suspends
by $\sigma$ to $0$.  Hence
$$
tr(\bar\epsilon_1) = \epsilon_0 \mu_0^{p-1} - \epsilon_1 \,.
$$
As we recalled above, the map $\hat\Gamma_1 \: \pi_* THH(\Z/p) \to
\pi_* THH(\Z/p)^{tC_p}$ is $(-2)$-coconnected, so the corresponding map
in $V(1)$-homotopy is at least $(2p-2)$-coconnected.  Thus it takes
$\epsilon_0 \mu_0^{p-1} - \epsilon_1$ to a nonzero class in $V(1)_*
THH(\Z/p)^{tC_p}$, represented somewhere in total degree~$(2p-1)$ of $\hat
E^\infty_{**}(C_p, \Z/p)$, in the lower right hand corner of the diagram.

Going down the middle of the diagram, we reach a class $(\Gamma_1
\circ tr_1)(\bar\epsilon_1)$, represented in total degree~$(2p-1)$
of the left half-plane $C_p$-homotopy fixed point spectral sequence
$E^\infty_{**}(C_p, \Z/p)$.  Its image under the edge homomorphism
to $V(1)_* THH(\Z/p)$ equals $(F \circ tr_1)(\bar\epsilon_1) =
tr(\bar\epsilon_1)$, hence $(\Gamma_1 \circ tr_1)(\bar\epsilon_1)$
is represented by $\epsilon_0 \mu_0^{p-1} - \epsilon_1$ in
$E^\infty_{0,2p-1}(C_p, \Z/p)$.  Its image under $R^h$ in the $C_p$-Tate
spectral sequence is the generator of $\hat E^\infty_{0,2p-1}(C_p, \Z/p)
= \F_p\{\epsilon_0 \mu_0^{p-1} - \epsilon_1\}$, hence that generator is
the $E^\infty$-representative of $\hat\Gamma_1(\epsilon_0 \mu_0^{p-1}
- \epsilon_1)$.
\qed
\enddemo

We can lift the algebraic $K$-theory class $\bar\epsilon_1$ to $\ell/p$.

\definition{Definition 6.3}
The map $\pi \: \ell/p \to H\Z/p$ is $(2p-2)$-connected, so it induces
a $(2p-1)$-connected map $V(1)_* K(\ell/p) \to V(1)_* K(\Z/p) =
E(\bar\epsilon_1)$, by \cite{BM94, 10.9}.  We can therefore choose a class
$$
\bar\epsilon_1^K \in V(1)_{2p-1} K(\ell/p)
$$
that maps to the generator $\bar\epsilon_1$ in $V(1)_* K(\Z/p)
\cong \Z/p$.
\enddefinition

\proclaim{Lemma 6.4}
The B{\"o}kstedt trace $tr \: V(1)_* K(\ell/p) \to V(1)_* THH(\ell/p)$
takes $\bar\epsilon_1^K$ to $\bar\epsilon_1$.
\endproclaim

\demo{Proof}
In the commutative square
$$
\xymatrix{
V(1)_* K(\ell/p) \ar[d]^{\pi^*} \ar[r]^-{tr} & V(1)_* THH(\ell/p)
  \ar[d]^{\pi^*} \\
V(1)_* K(\Z/p) \ar[r]^-{tr} & V(1)_* THH(\Z/p)
}
$$
the trace image $tr(\bar\epsilon_1^K)$ in $V(1)_* THH(\ell/p)$ must map
under $\pi^*$ to $tr(\bar\epsilon_1) = \epsilon_0 \mu_0^{p-1} - \epsilon_1$
in $V(1)_* THH(\Z/p)$, which by Proposition~5.4 characterizes it as
being equal to the class $\bar\epsilon_1$.
Hence $tr(\bar\epsilon_1^K) = \bar\epsilon_1$.
\qed
\enddemo

Next we turn to the $C_p$-Tate spectral sequence $\hat E^*(C_p,
\ell/p)$ in $V(1)$-homotopy for $THH(\ell/p)$.  
Its $E^2$-term is
$$
\hat E^2_{**}(C_p, \ell/p) = E(u_1) \otimes P(t^{\pm1}) \otimes
\F_p\{1, \epsilon_0, \mu_0, \epsilon_0 \mu_0, \dots, \mu_0^{p-1},
  \bar\epsilon_1\} \otimes E(\lambda_2) \otimes P(\mu_2) \,.
$$
We have $d^2(x) = t \cdot \sigma x$, where
$$
\sigma(\epsilon_0^\delta \mu_0^{i-1}) = \cases
\mu_0^i & \text{for $\delta=1$, $0<i<p$,} \\
0 & \text{otherwise}
\endcases
$$
is readily deduced from~(5.1), and $\sigma(\bar\epsilon_1) = 0$ since
$\bar\epsilon_1$ is in the image of $tr$.  Thus
$$
\hat E^3_{**}(C_p, \ell/p) = E(u_1) \otimes P(t^{\pm1}) \otimes
E(\bar\epsilon_1) \otimes E(\lambda_2) \otimes P(t\mu_2) \,.
\tag 6.5
$$
We prefer to use $t\mu_2$ rather than $\mu_2$ as a generator, since it
represents multiplication by $v_2$ in all module spectral sequences
over $E^*(S^1, \ell)$, by \cite{AuR02, 4.8}.

To proceed, we shall use that $\hat E^*(C_p ,\ell/p)$ is a module
over the spectral sequence for $THH(\ell)$.  We therefore recall the
structure of the latter spectral sequence, from \cite{AuR02, 5.5}.
It begins
$$
\hat E^2_{**}(C_p, \ell) =
E(u_1) \otimes P(t^{\pm1}) \otimes E(\lambda_1, \lambda_2)
\otimes P(\mu_2) \,.
$$
The classes $\lambda_1$, $\lambda_2$ and $t\mu_2$ are infinite
cycles, and the differentials
$$
\align
d^{2p}(t^{1-p}) &= t\lambda_1 \\
d^{2p^2}(t^{p-p^2}) &= t^p\lambda_2 \\
d^{2p^2+1}(u_1t^{-p^2}) &= t\mu_2
\endalign
$$
up to units in $\F_p$, which we will always suppress, leave the terms
$$
\align
\hat E^{2p+1}_{**}(C_p, \ell) &= E(u_1, \lambda_1, \lambda_2)
	\otimes P(t^{\pm p}, t\mu_2) \\
\hat E^{2p^2+1}_{**}(C_p, \ell) &= E(u_1, \lambda_1, \lambda_2)
	\otimes P(t^{\pm p^2}, t\mu_2) \\
\hat E^{2p^2+2}_{**}(C_p, \ell) &= E(\lambda_1, \lambda_2)
	\otimes P(t^{\pm p^2})
\endalign
$$
with $\hat E^{2p^2+2} = \hat E^\infty$, converging to $V(1)_*
THH(\ell)^{tC_p}$.  The comparison map $\hat\Gamma_1$ maps $\lambda_1$,
$\lambda_2$ and $\mu_2$ to $\lambda_1$, $\lambda_2$ and $t^{-p^2}$,
respectively, inducing the algebraic localization map and identification
$$
\hat\Gamma_1 \: V(1)_* THH(\ell) \to
V(1)_* THH(\ell) [\mu_2^{-1}] \cong V(1)_* THH(\ell)^{tC_p} \,.
$$

\proclaim{Lemma 6.6}
In $\hat E^*(C_p, \ell/p)$, the class
$u_1 t^{-p}$ supports the nonzero differential
$$
d^{2p^2}(u_1 t^{-p}) = u_1 t^{p^2-p} \lambda_2
$$
and does not survive to the $E^\infty$-term.
\endproclaim

\demo{Proof}
In $\hat E^*(C_p, \ell)$, there is such a nonzero differential, up to a
unit in $\F_p$, which we have already declared that we will suppress.
By naturality along $i \: \ell \to \ell/p$, it follows that there is
also such a differential in $\hat E^*(C_p, \ell/p)$.  It remains to
argue that the target is nonzero.  Considering the $E^3$-term in~(6.5),
the only possible source of a previous differential hitting $u_1
t^{p^2-p} \lambda_2$ is $\bar\epsilon_1$.  But $\bar\epsilon_1$ is in
an even column and $u_1 t^{p^2-p} \lambda_2$ is in an odd column.  By
naturality with respect to the Frobenius (group restriction) map from
the $S^1$-Tate spectral sequence to the $C_p$-Tate spectral sequence,
which takes $\hat E^2_{**}(S^1, B)$ isomorphically to the even columns
of $\hat E^2_{**}(C_p, B)$, any such differential from an even to an
odd column must be zero.
\qed
\enddemo

To determine the map $\hat\Gamma_1$ we use naturality with
respect to the map $\pi \: \ell/p \to H\Z/p$. 

\proclaim{Lemma 6.7}
The classes $1, \epsilon_0, \mu_0, \epsilon_0 \mu_0, \dots, \mu_0^{p-1}$
and $\bar\epsilon_1$ in $V(1)_* THH(\ell/p)$ map under $\hat\Gamma_1$
to classes in $V(1)_* THH(\ell/p)^{tC_p}$ that are represented in $\hat
E^\infty_{**}(C_p, \ell/p)$ by the permanent cycles $(u_1 t^{-1})^\delta
t^{-i}$ (on the horizontal axis) in degrees $\le (2p-2)$, and by the
permanent cycle $\bar\epsilon_1$ (on the vertical axis) in degree
$(2p-1)$.
\endproclaim

\demo{Proof}
In the commutative square
$$
\xymatrix{
V(1)_* THH(\ell/p) \ar[r]^-{\hat\Gamma_1} \ar[d]^{\pi^*} & V(1)_*
THH(\ell/p)^{tC_p} \ar[d]^{\pi^*} \\
V(1)_* THH(\Z/p) \ar[r]^-{\hat\Gamma_1} & V(1)_* THH(\Z/p)^{tC_p}
}
$$
the classes $\epsilon_0^\delta \mu_0^i$ in the upper left hand corner
map to classes in the lower right hand corner that are represented by
$(u_1 t^{-1})^\delta t^{-i}$ in degrees $\le (2p-2)$, and
$\bar\epsilon_1$ maps to $\epsilon_0 \mu_0^{p-1} - \epsilon_1$ in
degree $(2p-1)$.  This follows by combining Proposition~5.4 and
Lemma~6.1.

The first $(2p-1)$ of these are represented in maximal filtration (on
the horizontal axis), so their images in the upper right hand corner
must be represented by permanent cycles $(u_1 t^{-1})^\delta t^{-i}$ in
the Tate spectral sequence $\hat E^\infty_{**}(C_p, \ell/p)$.

The image of the last class, $\bar\epsilon_1$, in the upper right hand
corner could either be represented by $\bar\epsilon_1$ in
bidegree~$(0,2p-1)$ or by $u_1 t^{-p}$ in bidegree~$(2p-1,0)$.
However, the last class supports a differential $d^{2p^2}(u_1 t^{-p}) =
u_1 t^{p^2-p} \lambda_2$, by Lemma~6.6 above.  This only leaves the
other possibility, that $\hat\Gamma_1(\bar\epsilon_1)$ is represented
by $\bar\epsilon_1$ in $\hat E^\infty_{**}(C_p, \ell/p)$.
\qed
\enddemo

We proceed to determine the differential structure in $\hat E^*(C_p,
\ell/p)$, making use of the permanent cycles identified above.

\proclaim{Proposition 6.8}
The $C_p$-Tate spectral sequence in $V(1)$-homotopy for $THH(\ell/p)$
has
$$
\hat E^3_{**}(C_p, \ell/p) = E(u_1, \bar\epsilon_1, \lambda_2) \otimes
	P(t^{\pm1}, t\mu_2) \,.
$$
It has differentials generated by
$$
d^{2p^2-2p+2}(t^{p-p^2} \cdot t^{-i}\bar\epsilon_1) = t\mu_2
\cdot t^{-i} 
$$
for $0<i<p$, $d^{2p^2}(t^{p-p^2}) = t^p\lambda_2$ and
$d^{2p^2+1}(u_1 t^{-p^2}) = t\mu_2$.  The following terms are
$$
\align
\hat E^{2p^2-2p+3}_{**}(C_p, \ell/p) &= E(u_1, \lambda_2) \otimes
	\F_p\{t^{-i} \mid 0<i<p\} \otimes P(t^{\pm p}) \\
	&\qquad \oplus E(u_1, \bar\epsilon_1, \lambda_2) \otimes
	P(t^{\pm p}, t\mu_2) \\
\hat E^{2p^2+1}_{**}(C_p, \ell/p) &= E(u_1, \lambda_2) \otimes
	\F_p\{t^{-i} \mid 0<i<p\} \otimes P(t^{\pm p^2}) \\
	&\qquad \oplus E(u_1, \bar\epsilon_1, \lambda_2) \otimes
	P(t^{\pm p^2}, t\mu_2) \\
\hat E^{2p^2+2}_{**}(C_p, \ell/p) &= E(u_1, \lambda_2) \otimes
	\F_p\{t^{-i} \mid 0<i<p\} \otimes P(t^{\pm p^2}) \\
	&\qquad \oplus E(\bar\epsilon_1, \lambda_2) \otimes P(t^{\pm
	p^2}) \,.
\endalign
$$
The last term can be rewritten as
$$
\hat E^\infty(C_p, \ell/p) = \bigl( E(u_1) \otimes \F_p\{t^{-i} \mid
0<i<p\} \oplus E(\bar\epsilon_1) \bigr) \otimes E(\lambda_2) \otimes
P(t^{\pm p^2}) \,.
$$
\endproclaim

\demo{Proof}
We have already identified the $E^2$- and $E^3$-terms above.  The
$E^3$-term~(6.5) is generated over $\hat E^3(C_p, \ell)$ by (an
$\F_p$-basis for) $E(\bar\epsilon_1)$, so the next possible
differential is induced by $d^{2p}(t^{1-p}) = t\lambda_1$.  But
multiplication by $\lambda_1$ is trivial in $V(1)_* THH(\ell/p)$, by
Lemma~5.4, so $\hat E^3(C_p, \ell/p) = \hat E^{2p+1}(C_p, \ell/p)$.
This term is generated over $\hat E^{2p+1}(C_p, \ell)$ by $P_p(t^{-1})
\otimes E(\bar\epsilon_1)$.  Here $1, t^{-1}, \dots, t^{1-p}$ and
$\bar\epsilon_1$ are permanent cycles, by Lemma~6.7.  Any
$d^r$-differential before $d^{2p^2}$ must therefore originate on a
class $t^{-i} \bar\epsilon_1$ for $0<i<p$, and be of even length~$r$,
since these classes lie in even columns.  For bidegree reasons, the
first possibility is $r = 2p^2-2p+2$, so $\hat E^3(C_p, \ell/p) = \hat
E^{2p^2-2p+2}(C_p, \ell/p)$.

Multiplication by $v_2$ acts trivially on $V(1)_* THH(\ell)$ and
$V(1)_* THH(\ell)^{tC_p}$ for degree reasons, and therefore also on
$V(1)_* THH(\ell/p)$ and $V(1)_* THH(\ell/p)^{tC_p}$ by the module
structure.  The class $v_2$ maps to $t\mu_2$ in the $S^1$-Tate spectral
sequence for $\ell$, as recalled above, so multiplication by $v_2$ is
represented by multiplication by $t\mu_2$ in the $C_p$-Tate spectral
sequence for $\ell/p$.  Applied to the permanent cycles $(u_1 t^{-1})^\delta
t^{-i}$ in degrees $\le (2p-2)$, this implies that the products
$$
t\mu_2 \cdot (u_1 t^{-1})^\delta t^{-i}
$$
must be infinite cycles representing zero, i.e., they must be hit by
differentials.  In the cases $\delta=1$, $0\le i\le p-2$, these classes
in odd columns cannot be hit by differentials of odd length,
such as $d^{2p^2+1}$, so the only possibility is
$$
d^{2p^2-2p+2}(t^{p-p^2} \cdot (u_1 t^{-1})t^{-i} \bar\epsilon_1 )
	= t\mu_2 \cdot (u_1 t^{-1})t^{-i}
$$
for $0\le i\le p-2$.  By the module structure (consider multiplication
by $u_1$) it follows that
$$
d^{2p^2-2p+2}(t^{p-p^2} \cdot t^{-i} \bar\epsilon_1 )
	= t\mu_2 \cdot t^{-i}
$$
for $0<i<p$.  Hence we can compute from~(6.5) that
$$
\align
\hat E^{2p^2-2p+3}_{**}(C_p, \ell/p) &= E(u_1) \otimes P(t^{\pm p})
	\otimes \F_p\{t^{-i} \mid 0<i<p\} \otimes E(\lambda_2) \\
&\qquad \oplus E(u_1) \otimes P(t^{\pm p}) \otimes E(\bar\epsilon_1)
	\otimes E(\lambda_2) \otimes P(t\mu_2) \,.
\endalign
$$
This is generated over $\hat E^{2p+1}(C_p, \ell)$ by the permanent
cycles $1, t^{-1}, \dots, t^{1-p}$ and~$\bar\epsilon_1$, so the next
differential is induced by $d^{2p^2}(t^{p-p^2}) = t^p\lambda_2$.
This leaves
$$
\align
\hat E^{2p^2+1}_{**}(C_p, \ell/p) &= E(u_1) \otimes P(t^{\pm p^2})
	\otimes \F_p\{t^{-i} \mid 0<i<p\} \otimes E(\lambda_2) \\
&\qquad \oplus E(u_1) \otimes P(t^{\pm p^2}) \otimes E(\bar\epsilon_1)
	\otimes E(\lambda_2) \otimes P(t\mu_2) \,.
\endalign
$$
Finally, $d^{2p^2+1}(u_1 t^{-p^2}) = t\mu_2$ applies, and leaves
$$
\align
\hat E^{2p^2+2}_{**}(C_p, \ell/p) &= E(u_1) \otimes P(t^{\pm p^2})
	\otimes \F_p\{t^{-i} \mid 0<i<p\} \otimes E(\lambda_2) \\
&\qquad \oplus P(t^{\pm p^2}) \otimes E(\bar\epsilon_1) \otimes
	E(\lambda_2) \,.
\endalign
$$
For bidegree reasons, $\hat E^{2p^2+2} = \hat E^\infty$.
\qed
\enddemo

\proclaim{Proposition 6.9}
The comparison map $\hat\Gamma_1$ takes the classes $\epsilon_0^\delta
\mu_0^i$, $\bar\epsilon_1$, $\lambda_2$ and $\mu_2$ in $V(1)_*
THH(\ell/p)$ to classes in $V(1)_* THH(\ell/p)^{tC_p}$ represented by
$(u_1 t^{-1})^\delta t^{-i}$, $\bar\epsilon_1$, $\lambda_2$ and
$t^{-p^2}$ in $\hat E^\infty_{**}(C_p, \ell/p)$, respectively.
Thus
$$
V(1)_* THH(\ell/p)^{tC_p} \cong \F_p\{1, \epsilon_0, \mu_0, \epsilon_0
  \mu_0, \dots, \mu_0^{p-1}, \bar\epsilon_1\} \otimes E(\lambda_2) \otimes
  P(\mu_2^{\pm1})
$$
and $\hat\Gamma_1$ factors as the algebraic localization map and identification
$$
\hat\Gamma_1 \: V(1)_* THH(\ell/p) \to V(1)_* THH(\ell/p) [\mu_2^{-1}]
\cong V(1)_* THH(\ell/p)^{tC_p} \,.
$$
In particular, this map is $(2p-2)$-coconnected.
\endproclaim

\demo{Proof}
The action of the map $\hat\Gamma_1$ on the classes $1, \epsilon_0,
\mu_0, \epsilon_0 \mu_0, \dots, \mu_0^{p-1}$ and $\bar\epsilon_1$
was given in Lemma~6.7, and the action on the classes $\lambda_2$
and $\mu_2$ was already recalled from \cite{AuR02}.  The structure of
$V(1)_* THH(\ell/p)^{tC_p}$ is then immediate from the $E^\infty$-term
in Proposition~6.8.  The top class not in the image of $\hat\Gamma_1$
is $\bar\epsilon_1 \lambda_2 \mu_2^{-1}$, in degree~$(2p-2)$.
\qed
\enddemo

Recall that
$$
\align
TF(B) &= \holim_{n,F} THH(B)^{C_{p^n}} \\
TR(B) &= \holim_{n,R} THH(B)^{C_{p^n}}
\endalign
$$
are defined as the homotopy limits over the Frobenius and the restriction
maps $F, R \: THH(B)^{C_{p^n}} \to THH(B)^{C_{p^{n-1}}}$, respectively.

\proclaim{Corollary 6.10}
The comparison maps
$$
\align
\Gamma_n &\: THH(\ell/p)^{C_{p^n}} \to THH(\ell/p)^{hC_{p^n}} \\
\hat\Gamma_n &\: THH(\ell/p)^{C_{p^{n-1}}} \to THH(\ell/p)^{tC_{p^n}}
\endalign
$$
for $n\ge1$, and
$$
\align
\Gamma &\: TF(\ell/p) \to THH(\ell/p)^{hS^1} \\
\hat\Gamma &\: TF(\ell/p) \to THH(\ell/p)^{tS^1}
\endalign
$$
all induce $(2p-2)$-coconnected maps on $V(1)$-homotopy.
\endproclaim

\demo{Proof}
This follows from a theorem of Tsalidis \cite{Ts98} and Proposition~6.9
above, just like in \cite{AuR02, 5.7}.  See also \cite{BBLR:cf}.
\qed
\enddemo

\head 7. Higher fixed points \endhead

Let $n\ge1$.  Write $v_p(i)$ for the $p$-adic valuation of $i$.
Define a numerical function $\rho(-)$ by
$$
\align
\rho(2k-1) &= (p^{2k+1}+1)/(p+1) = p^{2k} - p^{2k-1} + \dots - p + 1 \\
\rho(2k) &= (p^{2k+2}-p^2)/(p^2-1) = p^{2k} + p^{2k-2} + \dots + p^2
\endalign
$$
for $k\ge0$, so $\rho(-1) = 1$ and $\rho(0) = 0$.  For even arguments,
$\rho(2k) = r(2k)$ as defined in \cite{AuR02, 2.5}.

In all of the following spectral sequences we know that $\lambda_2$,
$t\mu_2$ and $\bar\epsilon_1$ are infinite cycles.  For $\lambda_2$ and
$\bar\epsilon_1$ this follows from the $C_{p^n}$-fixed point analogue
of diagram~(6.2), by \cite{AuR02, 2.8} and Lemma~6.4.  For $t\mu_2$
it follows from \cite{AuR02, 4.8}, by naturality.

\proclaim{Theorem 7.1}
The $C_{p^n}$-Tate spectral sequence $\hat E^*(C_{p^n}, \ell/p)$ in
$V(1)$-homotopy for $THH(\ell/p)$ begins
$$
\hat E^2_{**}(C_{p^n}, \ell/p) = E(u_n, \lambda_2) \otimes \F_p\{1,
  \epsilon_0, \mu_0, \epsilon_0 \mu_0, \dots, \mu_0^{p-1}, \bar\epsilon_1\}
  \otimes P(t^{\pm1}, \mu_2)
$$
and converges to $V(1)_* THH(\ell/p)^{tC_{p^n}}$.  It is a module
spectral sequence over the algebra spectral sequence $\hat E^*(C_{p^n},
\ell)$ converging to $V(1)_* THH(\ell)^{tC_{p^n}}$.

There is an initial $d^2$-differential generated by
$$
d^2(\epsilon_0 \mu_0^{i-1}) = t \mu_0^i
$$
for $0<i<p$.
Next, there are $2n$ families of even length differentials generated by
$$
d^{2\rho(2k-1)}(t^{p^{2k-1}-p^{2k}+i} \cdot \bar\epsilon_1)
	= (t\mu_2)^{\rho(2k-3)} \cdot t^i
$$
for $v_p(i) = 2k-2$, for each $k = 1, \dots, n$, and
$$
d^{2\rho(2k)}(t^{p^{2k-1}-p^{2k}})
	= \lambda_2 \cdot t^{p^{2k-1}} \cdot (t\mu_2)^{\rho(2k-2)}
$$
for each $k = 1, \dots, n$.
Finally, there is a differential of odd length generated by
$$
d^{2\rho(2n)+1}(u_n \cdot t^{-p^{2n}}) = (t\mu_2)^{\rho(2n-2)+1} \,.
$$
\endproclaim

We shall prove Theorem~7.1 by induction on $n$.  The base case $n=1$ is
covered by Proposition~6.8.  We can therefore assume that Theorem~7.1
holds for some fixed $n\ge1$.  First we make the following deduction.

\proclaim{Corollary 7.2}
The initial differential in the $C_{p^n}$-Tate spectral sequence in
$V(1)$-homotopy for $THH(\ell/p)$ leaves
$$
\hat E^3_{**}(C_{p^n}, \ell/p) = E(u_n, \bar\epsilon_1, \lambda_2)
	\otimes P(t^{\pm1}, t\mu_2) \,.
$$
The next $2n$ families of differentials leave the intermediate terms
$$
\align
\hat E^{2\rho(1)+1}_{**}(C_{p^n}&, \ell/p) = E(u_n, \lambda_2)
	\otimes \F_p\{t^{-i} \mid 0<i<p\} \otimes P(t^{\pm p}) \\
&\qquad \oplus E(u_n, \bar\epsilon_1, \lambda_2) \otimes P(t^{\pm p},
	t\mu_2)
\endalign
$$
(for $m=1$),
$$
\align
\hat E^{2\rho(2m-1)+1}_{**}(C_{p^n}&, \ell/p) = E(u_n, \lambda_2)
	\otimes \F_p\{t^{-i} \mid 0<i<p\} \otimes P(t^{\pm p^2}) \\
&\quad \oplus \bigoplus_{k=2}^m E(u_n, \lambda_2) \otimes
	\F_p\{t^j \mid v_p(j) = 2k-2\} \otimes P_{\rho(2k-3)}(t\mu_2) \\
&\quad \oplus \bigoplus_{k=2}^{m-1} E(u_n, \bar\epsilon_1) \otimes
	\F_p\{t^j \lambda_2 \mid v_p(j) = 2k-1\} \otimes
	P_{\rho(2k-2)}(t\mu_2) \\
&\qquad \oplus E(u_n, \bar\epsilon_1, \lambda_2) \otimes P(t^{\pm
	p^{2m-1}}, t\mu_2)
\endalign
$$
for $m=2, \dots, n$, and
$$
\align
\hat E^{2\rho(2m)+1}_{**}(C_{p^n}, \ell/p) &= E(u_n, \lambda_2) \otimes
	\F_p\{t^{-i} \mid 0<i<p\} \otimes P(t^{\pm p^2}) \\
&\quad \oplus \bigoplus_{k=2}^m E(u_n, \lambda_2) \otimes
	\F_p\{t^j \mid v_p(j) = 2k-2\} \otimes P_{\rho(2k-3)}(t\mu_2) \\
&\quad \oplus \bigoplus_{k=2}^m E(u_n, \bar\epsilon_1) \otimes
	\F_p\{t^j \lambda_2 \mid v_p(j) = 2k-1\} \otimes
	P_{\rho(2k-2)}(t\mu_2) \\
&\qquad \oplus E(u_n, \bar\epsilon_1, \lambda_2) \otimes P(t^{\pm
	p^{2m}}, t\mu_2)
\endalign
$$
for $m = 1, \dots, n$.  The final differential leaves the
$E^{2\rho(2n)+2} = E^\infty$-term, equal to
$$
\align
\hat E^\infty_{**}(C_{p^n}, \ell/p) &= E(u_n, \lambda_2) \otimes
	\F_p\{t^{-i} \mid 0<i<p\} \otimes P(t^{\pm p^2})\\
&\quad \oplus \bigoplus_{k=2}^n E(u_n, \lambda_2) \otimes
	\F_p\{t^j \mid v_p(j) = 2k-2\} \otimes P_{\rho(2k-3)}(t\mu_2) \\
&\quad \oplus \bigoplus_{k=2}^n E(u_n, \bar\epsilon_1) \otimes
	\F_p\{t^j \lambda_2 \mid v_p(j) = 2k-1\} \otimes
	P_{\rho(2k-2)}(t\mu_2) \\
&\qquad \oplus E(\bar\epsilon_1, \lambda_2) \otimes P(t^{\pm
	p^{2n}}) \otimes P_{\rho(2n-2)+1}(t\mu_2) \,.
\endalign
$$
\endproclaim

\demo{Proof}
The statements about the $E^3$-, $E^{2\rho(1)+1}$- and
$E^{2\rho(2)+1}$-terms are clear from Proposition~6.8.  For each $m =
2, \dots, n$ we proceed by a secondary induction.  The differential
$$
d^{2\rho(2m-1)}(t^{p^{2m-1}-p^{2m}+i} \cdot \bar\epsilon_1)
	= (t\mu_2)^{\rho(2m-3)} \cdot t^i
$$
for $v_p(i) = 2m-2$ is non-trivial only on the summand
$$
E(u_n, \bar\epsilon_1, \lambda_2) \otimes P(t^{\pm p^{2m-2}}, t\mu_2)
$$
of the $E^{2\rho(2m-2)+1} = E^{2\rho(2m-1)}$-term, with homology
$$
\align
&E(u_n, \lambda_2) \otimes \F_p\{t^j \mid v_p(j) = 2m-2\} \otimes
	P_{\rho(2m-3)}(t\mu_2) \\
&\quad \oplus E(u_n, \bar\epsilon_1, \lambda_2) \otimes P(t^{\pm
	p^{2m-1}}, t\mu_2) \,.
\endalign
$$
This gives the stated $E^{2\rho(2m-1)+1}$-term.  Similarly, the differential
$$
d^{2\rho(2m)}(t^{p^{2m-1}-p^{2m}})
	= \lambda_2 \cdot t^{p^{2m-1}} \cdot (t\mu_2)^{\rho(2m-2)}
$$
is non-trivial only on the summand
$$
E(u_n, \bar\epsilon_1, \lambda_2) \otimes P(t^{\pm p^{2m-1}}, t\mu_2)
$$
of the $E^{2\rho(2m-1)+1} = E^{2\rho(2m)}$-term, with homology
$$
\align
&E(u_n, \bar\epsilon_1) \otimes \F_p\{t^j \lambda_2 \mid v_p(j) = 2m-1\}
	\otimes P_{\rho(2m-2)}(t\mu_2) \\
&\quad \oplus E(u_n, \bar\epsilon_1, \lambda_2) \otimes P(t^{\pm
	p^{2m}}, t\mu_2)
\,.
\endalign
$$
This gives the stated $E^{2\rho(2m)+1}$-term.  The final differential
$$
d^{2\rho(2n)+1}(u_n \cdot t^{-p^{2n}}) = (t\mu_2)^{\rho(2n-2)+1}
$$
is non-trivial only on the summand
$$
E(u_n, \bar\epsilon_1, \lambda_2) \otimes P(t^{\pm p^{2n}}, t\mu_2)
$$
of the $E^{2\rho(2n)+1}$-term, with homology
$$
E(\bar\epsilon_1, \lambda_2) \otimes P(t^{\pm p^{2n}}) \otimes
	P_{\rho(2n-2)+1}(t\mu_2)
\,.
$$
This gives the stated $E^{2\rho(2n)+2}$-term.  At this stage there is
no room for any further differentials, since the spectral sequence is
concentrated in a narrower horizontal band than the vertical height of
the following differentials.
\qed
\enddemo

Next we compare the $C_{p^n}$-Tate spectral sequence with the
$C_{p^n}$-homotopy spectral sequence obtained by restricting the
$E^2$-term to the second quadrant ($s\le0$, $t\ge0$).  It is
algebraically easier to handle the latter after inverting $\mu_2$,
which can be interpreted as comparing $THH(\ell/p)$ with its
$C_p$-Tate construction.

In general, there is a commutative diagram
$$
\xymatrix{
THH(B)^{C_{p^n}} \ar[r]^-R \ar[d]^{\Gamma_n}
  & THH(B)^{C_{p^{n-1}}} \ar[r]^-{\Gamma_{n-1}} \ar[d]^{\hat\Gamma_n}
  & THH(B)^{hC_{p^{n-1}}} \ar[d]^{\hat\Gamma_1^{hC_{p^{n-1}}}} \\
THH(B)^{hC_{p^n}} \ar[r]^-{R^h}
  & THH(B)^{tC_{p^n}} \ar[r]^-{G_{n-1}}
  & (THH(B)^{tC_p})^{hC_{p^{n-1}}}
}
\tag 7.3
$$
where $G_{n-1}$ is the comparison map from the $C_{p^{n-1}}$-fixed points
to the $C_{p^{n-1}}$-homotopy fixed points of $THH(B)^{tC_p}$, in view
of the identification
$$
(THH(B)^{tC_p})^{C_{p^{n-1}}} = THH(B)^{tC_{p^n}} \,.
$$

We are of course considering the case $B = \ell/p$.  In $V(1)$-homotopy
all four maps with labels containing $\Gamma$ are $(2p-2)$-coconnected, by
Corollary~6.10, so $G_{n-1}$ is at least $(2p-1)$-coconnected.  (We shall
see in Lemma~7.11 that $V(1)_* G_{n-1}$ is an isomorphism in all degrees.)
By Proposition~6.9 the map $\hat\Gamma_1$ precisely inverts $\mu_2$, so
the $E_2$-term of the $C_{p^n}$-homotopy fixed point spectral sequence in
$V(1)$-homotopy for $THH(\ell/p)^{tC_p}$ is obtained by inverting $\mu_2$
in $E^2_{**}(C_{p^n}, \ell/p)$.  We denote it by $\mu_2^{-1} E^*(C_{p^n},
\ell/p)$, even though in later terms only a power of $\mu_2$ is present.

\proclaim{Theorem 7.4}
The $C_{p^n}$-homotopy fixed point spectral sequence $\mu_2^{-1}
E^*(C_{p^n}, \ell/p)$ in $V(1)$-homotopy for $THH(\ell/p)^{tC_p}$
begins
$$
\mu_2^{-1} E^2_{**}(C_{p^n}, \ell/p) = E(u_n, \lambda_2) \otimes \F_p\{1,
  \epsilon_0, \mu_0, \epsilon_0 \mu_0, \dots, \mu_0^{p-1}, \bar\epsilon_1\}
  \otimes P(t, \mu_2^{\pm1})
$$
and converges to $V(1)_* (THH(\ell/p)^{tC_p})^{hC_{p^n}}$, which
receives a $(2p-2)$-coconnected map $(\hat\Gamma_1)^{hC_{p^n}}$ from
$V(1)_* THH(\ell/p)^{hC_{p^n}}$.
There is an initial $d^2$-differential generated by
$$
d^2(\epsilon_0 \mu_0^{i-1}) = t \mu_0^i
$$
for $0<i<p$.
Next, there are $2n$ families of even length differentials generated by
$$
d^{2\rho(2k-1)}(\mu_2^{p^{2k}-p^{2k-1}+j} \cdot \bar\epsilon_1)
	= (t\mu_2)^{\rho(2k-1)} \cdot \mu_2^j
$$
for $v_p(j) = 2k-2$, for each $k = 1, \dots, n$, and
$$
d^{2\rho(2k)}(\mu_2^{p^{2k}-p^{2k-1}})
	= \lambda_2 \cdot \mu_2^{-p^{2k-1}} \cdot (t\mu_2)^{\rho(2k)}
$$
for each $k = 1, \dots, n$.
Finally, there is a differential of odd length generated by
$$
d^{2\rho(2n)+1}(u_n \cdot \mu_2^{p^{2n}}) = (t\mu_2)^{\rho(2n)+1} \,.
$$
\endproclaim

\demo{Proof}
The differential pattern follows from Theorem~7.1 by naturality with
respect to the maps of spectral sequences
$$
\mu_2^{-1} E^*(C_{p^n}, \ell/p)
@<\hat\Gamma_1^{hC_{p^n}}<<
E^*(C_{p^n}, \ell/p)
@>R^h>>
\hat E^*(C_{p^n}, \ell/p)
$$
induced by $\hat\Gamma_1^{hC_{p^n}}$ and $R^h$.  The first inverts
$\mu_2$ and the second inverts $t$, at the level of $E^2$-terms.
We are also using that $t\mu_2$, the image of $v_2$, multiplies as an
infinite cycle in all of these spectral sequences.
\qed
\enddemo

\proclaim{Corollary 7.5}
The initial differential in the $C_{p^n}$-homotopy fixed point spectral
sequence in $V(1)$-homotopy for $THH(\ell/p)^{tC_p}$ leaves
$$
\align
\mu_2^{-1} E^3_{**}(C_{p^n}, \ell/p) &= E(u_n, \lambda_2) \otimes
	\F_p\{\mu_0^i \mid 0<i<p\} \otimes P(\mu_2^{\pm1}) \\
&\qquad \oplus E(u_n, \bar\epsilon_1, \lambda_2)
        \otimes P(\mu_2^{\pm1}, t\mu_2) \,.
\endalign
$$
The next $2n$ families of differentials leave the intermediate terms
$$
\align
\mu_2^{-1} E^{2\rho(2m-1)+1}_{**}(C_{p^n}&, \ell/p) = E(u_n, \lambda_2)
	\otimes \F_p\{\mu_0^i \mid 0<i<p\} \otimes P(\mu_2^{\pm1}) \\
&\oplus \bigoplus_{k=1}^m E(u_n, \lambda_2) \otimes
        \F_p\{\mu_2^j \mid v_p(j) = 2k-2\} \otimes P_{\rho(2k-1)}(t\mu_2) \\
&\oplus \bigoplus_{k=1}^{m-1} E(u_n, \bar\epsilon_1) \otimes
        \F_p\{\lambda_2 \mu_2^j \mid v_p(j) = 2k-1\} \otimes
        P_{\rho(2k)}(t\mu_2) \\
&\quad \oplus E(u_n, \bar\epsilon_1, \lambda_2) \otimes P(\mu_2^{\pm
        p^{2m-1}}, t\mu_2)
\endalign
$$
and
$$
\align
\mu_2^{-1} E^{2\rho(2m)+1}_{**}(C_{p^n}, \ell/p) &= E(u_n, \lambda_2)
	\otimes \F_p\{\mu_0^i \mid 0<i<p\} \otimes P(\mu_2^{\pm1}) \\
&\oplus \bigoplus_{k=1}^m E(u_n, \lambda_2) \otimes
        \F_p\{\mu_2^j \mid v_p(j) = 2k-2\} \otimes P_{\rho(2k-1)}(t\mu_2) \\
&\oplus \bigoplus_{k=1}^m E(u_n, \bar\epsilon_1) \otimes
	\F_p\{\lambda_2 \mu_2^j \mid v_p(j) = 2k-1\} \otimes
	P_{\rho(2k)}(t\mu_2) \\
&\quad \oplus E(u_n, \bar\epsilon_1, \lambda_2) \otimes P(\mu_2^{\pm
        p^{2m}}, t\mu_2)
\endalign
$$
for $m = 1, \dots, n$.  The final differential leaves the
$E^{2\rho(2n)+2} = E^\infty$-term, equal to
$$
\align
\mu_2^{-1} E^\infty_{**}(C_{p^n}, \ell/p) &= E(u_n, \lambda_2) \otimes
	\F_p\{\mu_0^i \mid 0<i<p\} \otimes P(\mu_2^{\pm1}) \\
&\quad \oplus \bigoplus_{k=1}^n E(u_n, \lambda_2) \otimes
	\F_p\{\mu_2^j \mid v_p(j) = 2k-2\} \otimes P_{\rho(2k-1)}(t\mu_2) \\
&\quad \oplus \bigoplus_{k=1}^n E(u_n, \bar\epsilon_1) \otimes
	\F_p\{\lambda_2 \mu_2^j \mid v_p(j) = 2k-1\} \otimes
	P_{\rho(2k)}(t\mu_2) \\
&\qquad \oplus E(\bar\epsilon_1, \lambda_2) \otimes P(\mu_2^{\pm
	p^{2n}}) \otimes P_{\rho(2n)+1}(t\mu_2) \,.
\endalign
$$
\endproclaim

\demo{Proof}
The computation of the $E^3$-term from the $E^2$-term is straightforward.
The rest of the proof goes by a secondary induction on $m = 1, \dots,
n$, very much like the proof of Corollary~7.2.  The differential
$$
d^{2\rho(2m-1)}(\mu_2^{p^{2m}-p^{2m-1}+j} \cdot \bar\epsilon_1)
	= (t\mu_2)^{\rho(2m-1)} \cdot \mu_2^j
$$
for $v_p(j) = 2m-2$ is non-trivial only on the summand
$$
E(u_n, \bar\epsilon_1, \lambda_2) \otimes P(\mu_2^{\pm p^{2m-2}}, t\mu_2)
$$
of the $E^3 = E^{2\rho(1)}$-term (for $m=1$), resp.~the
$E^{2\rho(2m-2)+1} = E^{2\rho(2m-1)}$-term (for $m = 2, \dots, n$).
Its homology is
$$
\align
&E(u_n, \lambda_2) \otimes \F_p\{ \mu_2^j \mid v_p(j) = 2m-2 \} \otimes
	P_{\rho(2m-1)}(t\mu_2) \\
&\quad\oplus E(u_n, \bar\epsilon_1, \lambda_2) \otimes P(\mu_2^{\pm
	p^{2m-1}}, t\mu_2) \,,
\endalign
$$
which gives the stated $E^{2\rho(2m-1)+1}$-term.
The differential
$$
d^{2\rho(2m)}(\mu_2^{p^{2m}-p^{2m-1}}) = \lambda_2 \cdot \mu_2^{-p^{2m-1}}
	\cdot (t\mu_2)^{\rho(2m)}
$$
is non-trivial only on the summand
$$
E(u_n, \bar\epsilon_1, \lambda_2) \otimes P(\mu_2^{\pm p^{2m-1}}, t\mu_2)
$$
of the $E^{2\rho(2m-1)+1} = E^{2\rho(2m)}$-term, leaving
$$
\align
&E(u_n, \bar\epsilon_1) \otimes \F_p\{ \lambda_2 \mu_2^j
        \mid v_p(j) = 2m-1 \} \otimes P_{\rho(2m)}(t\mu_2) \\
&\quad\oplus E(u_n, \bar\epsilon_1, \lambda_2) \otimes P(\mu_2^{\pm p^{2m}},
        t\mu_2) \,.
\endalign
$$
This gives the stated $E^{2\rho(2m)+1}$-term.
The final differential
$$
d^{2\rho(2n)+1}(u_n \cdot \mu_2^{p^{2n}}) = (t\mu_2)^{\rho(2n)+1}
$$
is non-trivial only on the summand
$$
E(u_n, \bar\epsilon_1, \lambda_2) \otimes P(\mu_2^{\pm p^{2n}}, t\mu_2)
$$
of the $E^{2\rho(2n)+1}$-term, with homology
$$
E(\bar\epsilon_1, \lambda_2) \otimes P(\mu_2^{\pm p^{2n}}) \otimes
        P_{\rho(2n)+1}(t\mu_2)
\,.
$$
This gives the stated $E^{2\rho(2n)+2}$-term.  There is
no room for any further differentials, since the spectral sequence is
concentrated in a narrower vertical band than the horizontal width of
the following differentials, so $E^{2\rho(2n)+2} = E^\infty$.
\qed
\enddemo

\demo{Proof of Theorem~7.1}
To make the inductive step to $C_{p^{n+1}}$, we use that
the first $d^r$-differential of odd length in $\hat E^*(C_{p^n}, \ell/p)$
occurs for $r = r_0 = 2\rho(2n)+1$.  It follows from \cite{AuR02, 5.2}
that the terms $\hat E^r(C_{p^n}, \ell/p)$ and $\hat E^r(C_{p^{n+1}},
\ell/p)$ are isomorphic for $r \le 2\rho(2n)+1$, via the Frobenius
map (taking $t^i$ to $t^i$) in even columns and the Verschiebung map
(taking $u_n t^i$ to $u_{n+1} t^i$) in odd columns.  Furthermore, the
differential $d^{2\rho(2n)+1}$ is zero in the latter spectral sequence.
This proves the part of Theorem~7.1 for $n+1$ that concerns the
differentials leading up to the term
$$
\align
\hat E^{2\rho(2n)+2}(C_{p^{n+1}}&, \ell/p) = E(u_{n+1}, \lambda_2) \otimes
        \F_p\{t^{-i} \mid 0<i<p\} \otimes P(t^{\pm p^2}) \\
&\quad \oplus \bigoplus_{k=2}^n E(u_{n+1}, \lambda_2) \otimes
        \F_p\{t^j \mid v_p(j) = 2k-2\} \otimes P_{\rho(2k-3)}(t\mu_2) \\
&\quad \oplus \bigoplus_{k=2}^n E(u_{n+1}, \bar\epsilon_1) \otimes
        \F_p\{t^j \lambda_2 \mid v_p(j) = 2k-1\} \otimes
        P_{\rho(2k-2)}(t\mu_2)
\tag 7.6 \\
&\qquad \oplus E(u_{n+1}, \bar\epsilon_1, \lambda_2) \otimes P(t^{\pm
        p^{2n}}, t\mu_2) \,.
\endalign
$$

Next we use the following commutative diagram, where we abbreviate
$THH(B)$ to $T(B)$:
$$
\xymatrix{
(T(B)^{tC_p})^{hC_{p^n}} \ar[d]^F
& T(B)^{hC_{p^n}} \ar[l]_-{\hat\Gamma_1^{hC_{p^n}}} \ar[d]^F
& T(B)^{C_{p^n}} \ar[l]_{\Gamma_n} \ar[r]^-{\hat\Gamma_{n+1}} \ar[d]^F
& T(B)^{tC_{p^{n+1}}} \ar[d]^F \\
T(B)^{tC_p}
& T(B) \ar[l]_{\hat\Gamma_1}
& T(B) \ar@{=}[l] \ar[r]^{\hat\Gamma_1}
& T(B)^{tC_p}
}
\tag 7.7
$$
The horizontal maps all induce $(2p-2)$-coconnected maps in
$V(1)$-homotopy for $B = \ell/p$.  Here $F$ is the Frobenius map,
forgetting part of the equivariance.  Thus the map $\hat\Gamma_{n+1}$ to
the right induces an isomorphism of $E(\lambda_2) \otimes P(v_2)$-modules
in all degrees $> (2p-2)$ from $V(1)_* THH(\ell/p)^{C_{p^n}}$, implicitly
identified to the left with the abutment of $\mu_2^{-1} E^*(C_{p^n},
\ell/p)$, to $V(1)_* THH(\ell/p)^{tC_{p^{n+1}}}$, which is the abutment
of $\hat E^*(C_{p^{n+1}}, \ell/p)$.  The diagram above ensures that the
isomorphism induced by $\hat\Gamma_{n+1}$ is compatible with the one
induced by $\hat\Gamma_1$.  By Proposition~6.9 it takes $\bar\epsilon_1$,
$\lambda_2$ and $\mu_2$ to $\bar\epsilon_1$, $\lambda_2$ and $t^{-p^2}$,
respectively, and similarly for monomials in these classes.

We focus on the summand
$$
E(u_n, \lambda_2) \otimes \F_p\{\mu_2^j \mid v_p(j) = 2n-2\}
\otimes P_{\rho(2n-1)}(t\mu_2)
$$
in $\mu_2^{-1} E^\infty_{**}(C_{p^n}, \ell/p)$, abutting to $V(1)_*
THH(\ell/p)^{C_{p^n}}$ in degrees $> (2p-2)$.  In the $P(v_2)$-module
structure on the abutment, each class $\mu_2^j$ with $v_p(j) = 2n-2$,
$j>0$, generates a copy of $P_{\rho(2n-1)}(v_2)$, since there are no
permanent cycles in the same total degree as $y = (t\mu_2)^{\rho(2n-1)}
\cdot \mu_2^j$ that have lower (= more negative) homotopy fixed point
filtration.  See Lemma~7.8 below for the elementary verification.  The
$P(v_2)$-module isomorphism induced by $\hat\Gamma_{n+1}$ must take this
to a copy of $P_{\rho(2n-1)}(v_2)$ in $V(1)_* THH(\ell/p)^{tC_{p^{n+1}}}$,
generated by $t^{-p^2j}$.

Writing $i = -p^2j$, we deduce that for $v_p(i) = 2n$, $i<0$, the
infinite cycle $z = (t\mu_2)^{\rho(2n-1)} \cdot t^i$ must represent zero
in the abutment, and must therefore be hit by a differential $z = d^r(x)$
in the $C_{p^{n+1}}$-Tate spectral sequence.  Here $r \ge 2\rho(2n)+2$.

Since $z$ generates a free copy of $P(t\mu_2)$ in the
$E^{2\rho(2n)+2}$-term displayed in~(7.6), and $d^r$ is
$P(t\mu_2)$-linear, the class $x$ cannot be annihilated by any power
of $t\mu_2$.  This means that $x$ must be contained in the summand
$$
E(u_{n+1}, \bar\epsilon_1, \lambda_2) \otimes P(t^{\pm p^{2n}}, t\mu_2)
$$
of $\hat E^{2\rho(2n)+2}_{**}(C_{p^{n+1}}, \ell/p)$.  By an elementary
check of bidegrees, see Lemma~7.9 below, the only possibility is that $x$
has vertical degree $(2p-1)$, so that we have differentials
$$
d^{2\rho(2n+1)}(t^{p^{2n+1}-p^{2n+2}+i} \cdot \bar\epsilon_1)
	= (t\mu_2)^{\rho(2n-1)} \cdot t^i
$$
for all $i<0$ with $v_p(i) = 2n$.  The cases $i>0$ follow by the module
structure over the $C_{p^{n+1}}$-Tate spectral sequence for $\ell$.
The remaining two differentials,
$$
d^{2\rho(2n+2)}(t^{p^{2n+1}-p^{2n+2}})
        = \lambda_2 \cdot t^{p^{2n+1}} \cdot (t\mu_2)^{\rho(2n)}
$$
and
$$
d^{2\rho(2n+2)+1}(u_{n+1} \cdot t^{-p^{2n+2}}) = (t\mu_2)^{\rho(2n)+1}
$$
are also present in the $C_{p^{n+1}}$-Tate spectral sequence for
$\ell$, see \cite{AuR02, 6.1}, hence follow in the present case by the
module structure.  With this we have established the complete
differential pattern asserted by Theorem~7.1.
\qed
\enddemo

\proclaim{Lemma 7.8}
For $v_p(j) = 2n-2$, $n\ge1$, there are no classes in $\mu_2^{-1}
E^\infty_{**}(C_{p^n}, \ell/p)$ in the same total degree as
$y = (t\mu_2)^{\rho(2n-1)} \cdot \mu_2^j$ that have lower homotopy fixed
point filtration.
\endproclaim

\demo{Proof}
The total degree of $y$ is $2(p^{2n+2} - p^{2n+1} + p - 1) + 2p^2 j \equiv
(2p-2) \mod 2p^{2n}$, which is even.

Looking at the formula for $\mu_2^{-1} E^\infty_{**}(C_{p^n}, \ell/p)$
in Corollary~7.5, the classes of lower filtration than $y$ all lie in
the terms
$$
E(u_n, \bar\epsilon_1) \otimes \F_p\{ \lambda_2 \mu_2^i \mid v_p(i)=2n-1 \}
  \otimes P_{\rho(2n)}(t\mu_2)
$$
and
$$
E(\bar\epsilon_1, \lambda_2) \otimes P(\mu_2^{\pm p^{2n}}) \otimes
  P_{\rho(2n)+1}(t\mu_2) \,.
$$
Those in even total degree and of lower filtration than $y$ are
$$
u_n \lambda_2 \cdot \mu_2^i (t\mu_2)^e, \quad
\bar\epsilon_1 \lambda_2 \cdot \mu_2^i (t\mu_2)^e
$$
with $v_p(i) = 2n-1$, $\rho(2n-1) < e < \rho(2n)$,
and
$$
\mu_2^i (t\mu_2)^e, \quad
\bar\epsilon_1 \lambda_2 \cdot \mu_2^i (t\mu_2)^e
$$
with $v_p(i) \ge 2n$, $\rho(2n-1) < e \le \rho(2n)$.

The total degree of $u_n \lambda_2 \cdot \mu_2^i (t\mu_2)^e$ for $v_p(i)
= 2n-1$ is $(-1) + (2p^2-1) + 2p^2 i + (2p^2-2)e \equiv (2p^2-2)(e+1)
\mod 2p^{2n}$.  For this to agree with the total degree of $y$, we must
have $(2p-2) \equiv (2p^2-2)(e+1) \mod 2p^{2n}$, so $(e+1) \equiv 1/(1+p)
\mod p^{2n}$ and $e \equiv \rho(2n-1) - 1 \mod p^{2n}$.  There is no
such $e$ with $\rho(2n-1) < e < \rho(2n)$.

The total degree of $\bar\epsilon_1 \lambda_2 \cdot \mu_2^i (t\mu_2)^e$
for $v_p(i) = 2n-1$ is $(2p-1) + (2p^2-1) + 2p^2 i + (2p^2-2)e \equiv 2p
+ (2p^2-2)(e+1) \mod 2p^{2n}$.  To agree with that
of $y$, we must have $(2p-2) \equiv 2p + (2p^2-2)(e+1) \mod 2p^{2n}$, so
$(e+1) \equiv 1/(1-p^2) \mod p^{2n}$ and $e \equiv \rho(2n) \mod p^{2n}$.
There is no such $e$ with $\rho(2n-1) < e < \rho(2n)$.

The total degree of $\mu_2^i (t\mu_2)^e$ for $v_p(i) \ge 2n$ is $2p^2 i +
(2p^2-2)e \equiv (2p^2-2)e \mod 2p^{2n}$.  To agree with that of $y$, we
must have $(2p-2) \equiv (2p^2-2)e \mod 2p^{2n}$, so $e \equiv 1/(1+p)
\equiv \rho(2n-1) \mod p^{2n}$.  There is no such $e$ with $\rho(2n-1)
< e \le \rho(2n)$.

The total degree of $\bar\epsilon_1 \lambda_2 \cdot \mu_2^i (t\mu_2)^e$
for $v_p(i) \ge 2n$ is $(2p-1) + (2p^2-1) + 2p^2 i + (2p^2-2)e$.
To agree modulo~$2p^{2n}$ with that of $y$, we must have $e \equiv
\rho(2n) \mod p^{2n}$.  The only such $e$ with $\rho(2n-1) < e \le
\rho(2n)$ is $e = \rho(2n)$.  But in that case, the total degree of
$\bar\epsilon_1 \lambda_2 \cdot \mu_2^i (t\mu_2)^e$ is $2p + 2p^2 i +
(2p^2-2)(\rho(2n)+1) = 2(p^{2n+2} + p - 1) + 2p^2 i$.  To be equal to that
of $y$, we must have $2p^2 i + 2p^{2n+1} = 2p^2 j$, which is impossible
for $v_p(i) \ge 2n$ and $v_p(j) = 2n-2$.
\qed
\enddemo

\proclaim{Lemma 7.9}
For $v_p(i) = 2n$, $n\ge1$ and $z = (t\mu_2)^{\rho(2n-1)} \cdot t^i$,
the only class in
$$
E(u_{n+1}, \bar\epsilon_1, \lambda_2) \otimes P(t^{\pm p^{2n}}, t\mu_2)
$$
that can support a differential $d^r(x) = z$
for $r \ge 2\rho(2n)+2$
is (a unit times)
$$
x = t^{p^{2n+1} - p^{2n+2} + i} \cdot \bar\epsilon_1 \,.
$$
\endproclaim

\demo{Proof}
The class $z$ has total degree $(2p^2-2)\rho(2n-1) - 2i = 2p^{2n+2}
- 2p^{2n+1} + 2p - 2 - 2i \equiv (2p-2) \mod 2p^{2n}$, which is even,
and vertical degree $2p^2 \rho(2n-1)$.  Hence $x$ has odd total degree,
and vertical degree at most $2p^2 \rho(2n-1) - 2\rho(2n) - 1 = 2p^{2n+2}
- 2p^{2n+1} - \dots - 2p^3 - 1$.  This leaves the possibilities
$$
u_{n+1} \cdot t^j (t\mu_2)^e, \quad
\bar\epsilon_1 \cdot t^j (t\mu_2)^e, \quad
\lambda_2 \cdot t^j (t\mu_2)^e
$$
with $v_p(j) \ge 2n$ and $0 \le e < p^{2n} - p^{2n-1} - \dots - p
= \rho(2n-1) - \rho(2n-2) - 1$,
and
$$
u_{n+1} \bar\epsilon_1 \lambda_2 \cdot t^j (t\mu_2)^e
$$
with $v_p(j) \ge 2n$ and $0 \le e < p^{2n} - p^{2n-1} - \dots - p - 1
= \rho(2n-1) - \rho(2n-2) - 2$.

The total degree of $x$ must be one more than the total degree of $z$,
hence is congruent to $(2p-1)$ modulo~$2p^{2n}$.

The total degree of $u_{n+1} \cdot t^j (t\mu_2)^e$ is $-1 - 2j +
(2p^2-2)e \equiv -1 + (2p^2-2)e \mod 2p^{2n}$.  To have $(2p-1) \equiv -1
+ (2p^2-2)e \mod 2p^{2n}$ we must have $e \equiv -p/(1-p^2) \equiv p^{2n}
- p^{2n-1} - \dots - p \mod p^{2n}$, which does not happen for $e$
in the allowable range.

The total degree of $\lambda_2 \cdot t^j (t\mu_2)^e$ is $(2p^2-1) - 2j +
(2p^2-2)e \equiv (2p^2-1) + (2p^2-2)e \mod 2p^{2n}$.  To have $(2p-1)
\equiv (2p^2-1) + (2p^2-2)e \mod 2p^{2n}$ we must have $e \equiv -p/(1+p)
\equiv \rho(2n-1) - 1 \mod p^{2n}$, which does not happen.

The total degree of $u_{n+1} \bar\epsilon_1 \lambda_2 \cdot t^j
(t\mu_2)^e$ is $-1 + (2p-1) + (2p^2-1) - 2j + (2p^2-2)e \equiv (2p -
1) + (2p^2-2)(e+1) \mod 2p^{2n}$.  To have $(2p-1) \equiv (2p - 1) +
(2p^2-2)(e+1) \mod 2p^{2n}$ we must have $(e+1) \equiv 0 \mod p^{2n}$,
so $e \equiv p^{2n} - 1 \mod p^{2n}$, which does not happen.

The total degree of $\bar\epsilon_1 \cdot t^j (t\mu_2)^e$ is $(2p-1) -
2j + (2p^2-2)e \equiv (2p-1) + (2p^2-2)e \mod 2p^{2n}$.  To have $(2p-1)
\equiv (2p-1) + (2p^2-2)e \mod 2p^{2n}$, we must have $e \equiv 0 \mod
p^{2n}$, so $e = 0$ is the only possibility in the allowable range.
In that case, a check of total degrees shows that we must have $j =
p^{2n+1} - p^{2n+2} + i$.
\qed
\enddemo

\proclaim{Corollary 7.10}
$V(1)_* THH(\ell/p)^{C_{p^n}}$ is finite in each degree.
\endproclaim

\demo{Proof}
This is clear by inspection of the $E^\infty$-term in
Corollary~7.2.
\qed
\enddemo

% ((Deduce what $\hat\Gamma_{n+1}$ does on $u_n$.))

\proclaim{Lemma 7.11}
The map $G_n$ induces an isomorphism
$$
V(1)_* THH(\ell/p)^{tC_{p^{n+1}}} @>\cong>>
V(1)_* (THH(\ell/p)^{tC_p})^{hC_{p^n}}
$$
in all degrees.  In the limit over the Frobenius maps $F$,
there is a map $G$ inducing an isomorphism
$$
V(1)_* THH(\ell/p)^{tS^1} @>\cong>>
V(1)_* (THH(\ell/p)^{tC_p})^{hS^1} \,.
$$
\endproclaim

\demo{Proof}
As remarked after diagram~(7.3), $G_n$ induces an isomorphism
in $V(1)$-homotopy above degree $(2p-2)$.  The permanent cycle
$t^{-p^{2n+2}}$ in $\hat E^\infty_{**}(C_{p^{n+1}}, \ell)$ acts
invertibly on $\hat E^\infty_{**}(C_{p^{n+1}}, \ell/p)$, and
its image $G_n(t^{-p^{2n+2}}) = \mu_2^{p^{2n}}$ in $\mu_2^{-1}
E^\infty_{**}(C_{p^n}, \ell)$ acts invertibly on $\mu_2^{-1}
E^\infty_{**}(C_{p^n}, \ell/p)$.  Therefore the module action derived
from the $\ell$-algebra structure on $\ell/p$ ensures that $G_n$ induces
isomorphisms in $V(1)$-homotopy in all degrees.
\qed
\enddemo

\proclaim{Theorem 7.12}
(a)
The associated graded of $V(1)_* THH(\ell/p)^{tS^1}$ for the
$S^1$-Tate spectral sequence is
$$
\align
\hat E^\infty_{**}(S^1, \ell/p) &= E(\lambda_2) \otimes
        \F_p\{t^{-i} \mid 0<i<p\} \otimes P(t^{\pm p^2})\\
&\quad \oplus \bigoplus_{k\ge2} E(\lambda_2) \otimes
        \F_p\{t^j \mid v_p(j) = 2k-2\} \otimes P_{\rho(2k-3)}(t\mu_2) \\
&\quad \oplus \bigoplus_{k\ge2} E(\bar\epsilon_1) \otimes
        \F_p\{t^j \lambda_2 \mid v_p(j) = 2k-1\} \otimes
        P_{\rho(2k-2)}(t\mu_2) \\
&\qquad \oplus E(\bar\epsilon_1, \lambda_2) \otimes P(t\mu_2) \,.
\endalign
$$

(b)
The associated graded of $V(1)_* THH(\ell/p)^{hS^1}$ for the
$S^1$-homotopy fixed point spectral sequence maps by a
$(2p-2)$-coconnected map to
$$
\align
\mu_2^{-1} E^\infty_{**}(S^1, \ell/p) &= E(\lambda_2) \otimes
        \F_p\{\mu_0^i \mid 0<i<p\} \otimes P(\mu_2^{\pm1}) \\
&\quad \oplus \bigoplus_{k\ge1} E(\lambda_2) \otimes
        \F_p\{\mu_2^j \mid v_p(j) = 2k-2\} \otimes P_{\rho(2k-1)}(t\mu_2) \\
&\quad \oplus \bigoplus_{k\ge1} E(\bar\epsilon_1) \otimes
        \F_p\{\lambda_2 \mu_2^j \mid v_p(j) = 2k-1\} \otimes
        P_{\rho(2k)}(t\mu_2) \\
&\qquad \oplus E(\bar\epsilon_1, \lambda_2) \otimes P(t\mu_2) \,.
\endalign
$$

(c)
The isomorphism from~(a) to~(b) induced by $G$ takes $t^{-i}$ to
$\mu_0^i$ for $0<i<p$ and $t^i$ to $\mu_2^j$ for $i + p^2j = 0$.
Furthermore, it takes multiples by $\bar\epsilon_1$, $\lambda_2$ or
$t\mu_2$ in the source to the same multiples in the target.
\endproclaim

\demo{Proof}
Claims~(a) and~(b) follow by passage to the limit over $n$ from
Corollaries~7.2 and~7.5.  Claim~(c) follows by passage to the same limit
from the formulas for the isomorphism induced by $\hat\Gamma_{n+1}$,
which were given below diagram~(7.7).
\qed
\enddemo

\head 8. Topological cyclic homology \endhead

By definition, there is a fiber sequence
$$
TC(B) @>\pi>> TF(B) @>R-1>> TF(B)
$$
inducing a long exact sequence
$$
\dots @>\partial>> V(1)_* TC(B) @>\pi>> V(1)_* TF(B) @>R-1>>
V(1)_* TF(B) @>\partial>> \dots
\tag 8.1
$$
in $V(1)$-homotopy.  By Corollary~6.10, there are $(2p-2)$-coconnected
maps $\Gamma$ and $\hat\Gamma$ from $V(1)_* TF(\ell/p)$ to $V(1)_*
THH(\ell/p)^{hS^1}$ and $V(1)_* THH(\ell/p)^{tS^1}$, respectively.  We
model $V(1)_* TF(\ell/p)$ in degrees $> (2p-2)$ by the map $\hat\Gamma$
to the $S^1$-Tate construction.  Then, by diagram~(7.3), $R$ is modeled
in the same range of degrees by the chain of maps below.
$$
\xymatrix{
V(1)_* THH(B)^{tS^1} \ar[dr]^G &
V(1)_* THH(B)^{hS^1} \ar[d]^{(\hat\Gamma_1)^{hS^1}} \ar[r]^{R^h} &
V(1)_* THH(B)^{tS^1} \\
& V(1)_* (THH(B)^{tC_p})^{hS^1}
}
$$
Here $R^h$ induces a map of spectral sequences
$$
E^*(R^h) \: E^*(S^1, B) \to \hat E^*(S^1, B) \,,
$$
which at the $E^2$-term equals the inclusion that algebraically inverts
$t$.  When $B = \ell/p$, the left hand map $G$ is an isomorphism by
Lemma~7.11, and the middle (wrong-way) map is $(2p-2)$-coconnected.

\proclaim{Proposition 8.2}
In degrees $> (2p-2)$, the homomorphism
$$
E^\infty(R^h) \: E^\infty(S^1, \ell/p) \to \hat E^\infty(S^1, \ell/p)
$$
maps

(a)
$E(\bar\epsilon_1, \lambda_2) \otimes P(t\mu_2)$ identically to the
same expression;

(b)
$E(\lambda_2) \otimes \F_p\{\mu_2^{-j}\} \otimes P_{\rho(2k-1)}(t\mu_2)$
surjectively onto
$$
E(\lambda_2) \otimes \F_p\{t^j\} \otimes P_{\rho(2k-3)}(t\mu_2)
$$
for each $k\ge2$, $j = dp^{2k-2}$, $0<d<p^2-p$ and $p\nmid d$;

(c)
$E(\bar\epsilon_1) \otimes \F_p\{\lambda_2 \mu_2^{-j}\} \otimes
P_{\rho(2k)}(t\mu_2)$
surjectively onto
$$
E(\bar\epsilon_1) \otimes \F_p\{t^j \lambda_2\} \otimes
P_{\rho(2k-2)}(t\mu_2)
$$
for each $k\ge2$, $j = dp^{2k-1}$ and $0<d<p$;

(d)
the remaining terms to zero.
\endproclaim

\demo{Proof}
Consider the summands of $E^\infty(S^1, \ell/p)$ and $\hat
E^\infty(S^1, \ell/p)$, as given in Theorem~7.12.  Clearly, the first
term $E(\lambda_2) \otimes \F_p\{\mu_0^i \mid 0<i<p\} \otimes P(\mu_2)$
goes to zero (these classes are hit by $d^2$-differentials), and the
last term $E(\bar\epsilon_1, \lambda_2) \otimes P(t\mu_2)$ maps
identically to the same term.  This proves~(a) and part of~(d).

For each $k\ge1$ and $j = dp^{2k-2}$ with $p\nmid d$, the term
$E(\lambda_2) \otimes \F_p\{\mu_2^{-j}\} \otimes
P_{\rho(2k-1)}(t\mu_2)$ maps to the term $E(\lambda_2) \otimes
\F_p\{t^j\} \otimes P_{\rho(2k-3)}(t\mu_2)$, except that the target is
zero for $k=1$.  In symbols, the element $\lambda_2^\delta \mu_2^{-j}
(t\mu_2)^i$ maps to the element $\lambda_2^\delta t^j (t\mu_2)^{i-j}$.
If $d<0$, then the $t$-exponent in the target is bounded above by
$dp^{2k-2} + \rho(2k-3) < 0$, so the target lives in the right
half-plane and is essentially not hit by the source, which lives in the
left half-plane.  If $d>p^2-p$, then the total degree in the source is
bounded above by $(2p^2-1) - 2dp^{2k} + \rho(2k-1)(2p^2-2) < 2p-2$, so
the source lives in total degree $< (2p-2)$ and will be disregarded.
If $0<d<p^2-p$, then $\rho(2k-1)-dp^{2k-2} > \rho(2k-3)$ and
$-dp^{2k-2} < 0$, so the source surjects onto the target.  This
proves~(b) and part of~(d).

Lastly, for each $k\ge1$ and $j = dp^{2k-1}$ with $p\nmid d$, the term
$E(\bar\epsilon_1) \otimes \F_p\{\lambda_2 \mu_2^{-j}\} \otimes
P_{\rho(2k)}(t\mu_2)$ maps to the term $E(\bar\epsilon_1) \otimes
\F_p\{t^j \lambda_2\} \otimes P_{\rho(2k-2)}(t\mu_2)$.  The target is
zero for $k=1$.  If $d<0$, then $dp^{2k-1} + \rho(2k-2) < 0$ so the
target lives in the right half-plane.  If $d>p$, then $(2p-1) +
(2p^2-1) - 2dp^{2k+1} + \rho(2k)(2p^2-2) < 2p-2$, so the source lives
in total degree $< (2p-2)$.  If $0<d<p$, then $\rho(2k)-dp^{2k-1} >
\rho(2k-2)$ and $-dp^{2k-1} < 0$, so the source surjects onto the
target.  This proves~(c) and the remaining part of~(d).
\qed
\enddemo

\definition{Definition 8.3}
Let
$$
\align
A &= E(\bar\epsilon_1, \lambda_2) \otimes P(t\mu_2) \\
B_k &= E(\lambda_2) \otimes \F_p\{t^{dp^{2k-2}} \mid 0<d<p^2-p,
	p\nmid d\} \otimes P_{\rho(2k-3)}(t\mu_2) \\
C_k &= E(\bar\epsilon_1) \otimes \F_p\{t^{dp^{2k-1}} \lambda_2 \mid
	0<d<p\} \otimes P_{\rho(2k-2)}(t\mu_2)
\endalign
$$
for $k\ge2$ and let $D$ be the span of the remaining monomials in $\hat
E^\infty(S^1, \ell/p)$.  Let $B = \bigoplus_{k\ge2} B_k$ and $C =
\bigoplus_{k\ge2} C_k$.  Then $\hat E^\infty(S^1, \ell/p) = A \oplus B
\oplus C \oplus D$.
\enddefinition

\proclaim{Proposition 8.4}
In degrees $> (2p-2)$, there are closed subgroups $\widetilde A =
E(\bar\epsilon_1, \lambda_2) \otimes P(v_2)$, $\widetilde B_k$, $\widetilde
C_k$ and $\widetilde D$ in $V(1)_* TF(\ell/p)$, represented by $A$, $B_k$,
$C_k$ and $D$ in $\hat E^\infty(S^1, \ell/p)$, respectively, such that
the homomorphism induced by the restriction map $R$

(a)
is the identity on $\widetilde A$;

(b)
maps $\widetilde B_{k+1}$ surjectively onto $\widetilde B_k$ for all $k\ge2$;

(c)
maps $\widetilde C_{k+1}$ surjectively onto $\widetilde C_k$ for all $k\ge2$;

(d)
is zero on $\widetilde B_2$, $\widetilde C_2$ and $\widetilde D$.

\noindent
In these degrees, $V(1)_* TF(\ell/p) \cong \widetilde A \oplus \widetilde B
\oplus \widetilde C \oplus \widetilde D$, where $\widetilde B = \prod_{k\ge2}
\widetilde B_k$ and $\widetilde C = \prod_{k\ge2} \widetilde C_k$.
\endproclaim

\demo{Proof}
In terms of the model $THH(\ell/p)^{tS^1}$ for $TF(\ell/p)$, the
restriction map $R$ is given in these degrees as the composite of
the isomorphism $G$, computed in Theorem~7.12(c), and the map $\hat
E^\infty(R^h)$, computed in Proposition~8.2.  This gives the desired
formulas at the level of $E^\infty$-terms.  The rest of the argument is
the same as that for Theorem~7.7 of \cite{AuR02}, using Corollary~7.10
to control the topologies, and will be omitted.
\qed
\enddemo

\remark{Remark 8.5}
Here we have followed the basic computational strategy of \cite{BM94},
\cite{BM95} and \cite{AuR02}.  It would be desirable to have a more
concrete construction of the lifts $\widetilde B_k$, $\widetilde C_k$ and
$\widetilde D$, in terms of de\,Rham--Witt operators $R$, $F$, $V$ and
$d = \sigma$, like in the algebraic case of \cite{HM97} and
\cite{HM03}.
\endremark

\proclaim{Proposition 8.6}
In degrees $>(2p-2)$ there are isomorphisms
$$
\align
\ker(R-1) &\cong \widetilde A \oplus \lim_k \widetilde B_k
	\oplus \lim_k \widetilde C_k \\
	&\cong E(\bar\epsilon_1, \lambda_2) \otimes P(v_2) \\
&\qquad \oplus E(\lambda_2) \otimes \F_p\{t^d \mid 0<d<p^2-p,
	p\nmid d\} \otimes P(v_2) \\
&\qquad \oplus E(\bar\epsilon_1) \otimes \F_p\{t^{dp} \lambda_2 \mid
	0<d<p\} \otimes P(v_2)
\endalign
$$
and $\cok(R-1) \cong \widetilde A = E(\bar\epsilon_1, \lambda_2)
\otimes P(v_2)$.  Hence there is an isomorphism
$$
\align
V(1)_* TC(\ell/p) &\cong
	E(\partial, \bar\epsilon_1, \lambda_2) \otimes P(v_2) \\
&\qquad \oplus E(\lambda_2) \otimes \F_p\{t^d \mid 0<d<p^2-p,
	p\nmid d\} \otimes P(v_2) \\
&\qquad \oplus E(\bar\epsilon_1) \otimes \F_p\{t^{dp} \lambda_2 \mid
	0<d<p\} \otimes P(v_2)
\endalign
$$
in these degrees, where $\partial$ has degree $-1$ and represents the
image of $1$ under the connecting map $\partial$ in~$(8.1)$.
\endproclaim

\demo{Proof}
By Proposition~8.4, the homomorphism $R-1$ is zero on $\widetilde A$
and an isomorphism on $\widetilde D$.  Furthermore, there is an exact
sequence
$$
0 \to \lim_k \widetilde B_k \to \prod_{k\ge2} \widetilde B_k
@>R-1>> \prod_{k\ge2} \widetilde B_k \to \Rlim_k \widetilde B_k \to 0
$$
and similarly for the $C$'s.  The derived limit on the right vanishes
since each $\widetilde B_{k+1}$ surjects onto $\widetilde B_k$.

Multiplication by $t\mu_2$ in each $B_k$ is realized by multiplication by
$v_2$ in $\widetilde B_k$.  Each $\widetilde B_k$ is a sum of $2(p-1)^2$
cyclic $P(v_2)$-modules, and since $\rho(2k-3)$ grows to infinity with $k$
their limit is a free $P(v_2)$-module of the same rank, with the indicated
generators $t^d$ and $t^d \lambda_2$ for $0 < d < p^2-p$, $p \nmid d$.
The argument for the $C$'s is practically the same.

The long exact sequence~(8.1) yields the short exact sequence
$$
0 \to \Sigma^{-1} \cok(R-1) @>\partial>> V(1)_* TC(\ell/p) @>\pi>>
\ker(R-1) \to 0 \,,
$$
from which the formula for the middle term follows.
\qed
\enddemo

\remark{Remark 8.7}
A more obvious set of $E(\lambda_2) \otimes P(v_2)$-module generators
for $\lim_k \widetilde B_k$ would be the classes $t^{dp^2}$ in $B_2
\cong \widetilde B_2$, for $0 < d < p^2-p$, $p \nmid d$.  Under the
canonical map $TF(\ell/p) \to THH(\ell/p)^{C_p}$, modeled here by
$THH(\ell/p)^{tS^1} \to (THH(\ell/p)^{tC_p})^{hC_p}$, these map to
the classes $\mu_2^{-d}$.  Since we are only concerned with degrees
$> (2p-2)$ we may equally well use their $v_2$-power multiplies
$(t\mu_2)^d \cdot \mu_2^{-d} = t^d$ as generators, with the advantage
that these are in the image of the localization map $THH(\ell/p)^{hC_p}
\to (THH(\ell/p)^{tC_p})^{hC_p}$.  Hence the class denoted $t^d$ in
$\lim_k \widetilde B_k$ is chosen so as to map under $TF(\ell/p) \to
THH(\ell/p)^{hC_p}$ to $t^d$ in $E^\infty_{**}(C_p; \ell/p)$.  Similarly,
the class denoted $t^{dp} \lambda_2$ in $\lim_k \widetilde C_k$ is chosen
so as to map to $t^{dp} \lambda_2$ in $E^\infty_{**}(C_p; \ell/p)$.
\endremark

\medskip

The map $\pi \: \ell/p \to \Z/p$ is $(2p-2)$-connected, hence induces
$(2p-1)$-connected maps $\pi^* \: K(\ell/p) \to K(\Z/p)$ and $\pi^* \:
V(1)_* TC(\ell/p) \to V(1)_* TC(\Z/p)$, by \cite{BM94, 10.9} and
\cite{Du97}.  Here $TC(\Z/p) \simeq H\Z_p \vee \Sigma^{-1} H\Z_p$ and
$V(1)_* TC(\Z/p) \cong E(\partial, \bar\epsilon_1)$, so we can recover
$V(1)_* TC(\ell/p)$ in degrees $\le (2p-2)$ from this map.

\proclaim{Theorem 8.8}
There is an isomorphism of $E(\lambda_1, \lambda_2) \otimes P(v_2)$-modules
$$
\align
V(1)_* TC(\ell/p) &\cong P(v_2) \otimes E(\partial, \bar\epsilon_1,
	\lambda_2) \\
&\qquad \oplus P(v_2) \otimes E(\dlog v_1) \otimes \F_p\{t^d
	v_2 \mid 0<d<p^2-p, p\nmid d\} \\
&\qquad \oplus P(v_2) \otimes E(\bar\epsilon_1) \otimes \F_p\{t^{dp}
	\lambda_2 \mid 0<d<p\}
\endalign
$$
where $v_2 \cdot \dlog v_1 = \lambda_2$.  The degrees are $|\partial| =
-1$, $|\bar\epsilon_1| = |\lambda_1| = 2p-1$, $|\lambda_2| = 2p^2-1$
and $|v_2| = 2p^2-2$.  The formal multipliers have degrees $|t| = -2$
and $|\dlog v_1| = 1$.
\endproclaim

The notation $\dlog v_1$ for the multiplier $v_2^{-1} \lambda_2$ is
suggested by the relation $v_1 \cdot \dlog p = \lambda_1$ in $V(0)_*
TC(\Z_{(p)}|\Q)$.

% ((Discuss $\lambda_1$-action.))

\demo{Proof}
Only the additive generators $t^d$ for $0<d<p^2-p$, $p\nmid d$ from
Proposition~8.6 do not appear in $V(1)_* TC(\ell/p)$, but their
multiples by $\lambda_2$ and positive powers of $v_2$ do.  This leads
to the given formula, where $\dlog v_1 \cdot t^d v_2$ must be read as
$t^d\lambda_2$.
\qed
\enddemo

By \cite{HM97} the cyclotomic trace map of \cite{BHM93} induces
cofiber sequences
$$
K(B_p)_p @>trc>> TC(B)_p @>g>> \Sigma^{-1} H\Z_p
\tag 8.9
$$
for each connective $S$-algebra $B$ with $\pi_0(B_p) = \Z_p$ or $\Z/p$,
and thus long exact sequences
$$
\dots \to V(1)_* K(B_p) @>trc>> V(1)_* TC(B) @>g>> \Sigma^{-1}
E(\bar\epsilon_1) \to \dots \,.
$$
This uses the identifications $W(\Z_p)_F \cong W(\Z/p)_F \cong \Z_p$ of
Frobenius coinvariants of Witt rings, and applies in particular for $B
= H\Z_{(p)}$, $H\Z/p$, $\ell$ and $\ell/p$.

\proclaim{Theorem 8.10}
There is an isomorphism of $E(\lambda_1, \lambda_2) \otimes P(v_2)$-modules
$$
\align
V(1)_* K(\ell/p) &\cong P(v_2) \otimes E(\bar\epsilon_1) \otimes
	\F_p\{1, \partial\lambda_2, \lambda_2, \partial v_2\} \\
&\qquad \oplus P(v_2) \otimes E(\dlog v_1) \otimes \F_p\{t^d
	v_2 \mid 0<d<p^2-p, p\nmid d\} \\
&\qquad \oplus P(v_2) \otimes E(\bar\epsilon_1) \otimes \F_p\{t^{dp}
	\lambda_2 \mid 0<d<p\} \,.
\endalign
$$
This is a free $P(v_2)$-module of rank~$(2p^2-2p+8)$ and of zero Euler
characteristic, where $p\ge5$ is assumed.
\endproclaim

\demo{Proof}
In the case $B = \Z/p$, $K(\Z/p)_p \simeq H\Z_p$ and the map $g$ is
split surjective up to homotopy.  So the induced homomorphism to $V(1)_*
\Sigma^{-1} H\Z_p = \Sigma^{-1} E(\bar\epsilon_1)$ is surjective.  Since
$\pi \: \ell/p \to \Z/p$ induces a $(2p-1)$-connected map in topological
cyclic homology, and $\Sigma^{-1} E(\bar\epsilon_1)$ is concentrated
in degrees $\le (2p-2)$, it follows by naturality that also in the
case $B = \ell/p$ the map $g$ induces a surjection in $V(1)$-homotopy.
The kernel of the surjection $P(v_2) \otimes E(\partial, \bar\epsilon_1,
\lambda_2) \to \Sigma^{-1} E(\bar\epsilon_1)$ gives the first row in
the asserted formula.
\qed
\enddemo

% ((Explain/recall the origin of the classes?))

\head 9. The fraction field \endhead

We wish to determine the effect on algebraic $K$-theory of excising the
closed subspaces $\Spec(\Z_p)$ and $\Spnc(\ell/p)$ from
$\Spec(\ell_p)$, meeting at the closed point $\Spec(\Z/p)$.  It will be
a little cleaner to express the computation in topological cyclic
homology.  Therefore we consider another $3 \times 3$ diagram of
cofiber sequences
$$
\xymatrix{
TC(\Z/p) \ar[r]^{i_*} \ar[d]^{\pi_*} &
TC(\Z_{(p)}) \ar[r]^{j^*} \ar[d]^{\pi_*} &
TC(\Z_{(p)}|\Q) \ar[d]^{\pi_*} \\
TC(\ell/p) \ar[r]^{i_*} \ar[d]^{\rho^*} &
TC(\ell) \ar[r]^{j^*} \ar[d]^{\rho^*} &
TC(\ell|p^{-1}\ell) \ar[d]^{\rho^*} \\
TC(\ell/p|L/p) \ar[r]^{i_*} &
TC(\ell|L) \ar[r]^{j^*} &
TC(\ff(\ell)) \,,
}
\tag 9.1
$$
generated as usual by the upper left hand square.  Like in the $THH$-case
displayed in diagram~(5.5), the upper horizontal transfer map $i_*$ was
defined in \cite{HM03}, and the remaining transfer maps can be defined
using \cite{BM:loc}.

In the top row, the transfer map $i_*$ from
$$
V(1)_* TC(\Z/p) = E(\partial, \bar\epsilon_1)
$$
to
$$
V(1)_* TC(\Z_{(p)}) = E(\partial, \lambda_1)
	\oplus \F_p\{t^d \lambda_1 \mid 0<d<p\}
$$
is a module map over the target.  We claim that $i_*(1) = 0$ and
$i_*(\bar\epsilon_1) = \lambda_1$, up to a unit in $\F_p$, so
$$
V(1)_* TC(\Z_{(p)}|\Q) = E(\partial, \dlog p)
	\oplus \F_p\{t^d \lambda_1 \mid 0<d<p\} \,,
$$
where the connecting map from~(9.1) takes the degree~$1$ class $\dlog p$
to $1$.  To prove the claim, recall from \cite{HM03} that $v_1 \cdot
\dlog p = \lambda_1$ in the $V(0)$-homotopy of $TC(\Z_{(p)}|\Q)$.
Hence $\lambda_1$ in $V(1)_* TC(\Z_{(p)})$ maps to zero under $j^*$,
and must be in the image of $i_*$.  The only class that can hit it is
$\bar\epsilon_1$, up to a unit.

In the middle row, the transfer map $i_*$ from
$$
\align
V(1)_* TC(\ell/p) &= P(v_2) \otimes E(\partial, \bar\epsilon_1, \lambda_2) \\
&\qquad \oplus P(v_2) \otimes E(\dlog v_1) \otimes
	\F_p\{t^d v_2 \mid 0<d<p^2-p, p\nmid d\} \\
&\qquad \oplus P(v_2) \otimes E(\bar\epsilon_1) \otimes
	\F_p\{t^{dp} \lambda_2 \mid 0<d<p\}
\endalign
$$
(computed in Theorem~8.8), to
$$
\align
V(1)_* TC(\ell) &= P(v_2) \otimes E(\partial, \lambda_1, \lambda_2) \\
&\qquad \oplus P(v_2) \otimes E(\lambda_2) \otimes
	\F_p\{t^d \lambda_1 \mid 0<d<p\} \\
&\qquad \oplus P(v_2) \otimes E(\lambda_1) \otimes
	\F_p\{t^{dp} \lambda_2 \mid 0<d<p\}
\endalign
$$
(computed in \cite{AuR02, 8.4}), is also a module map over the target.
By naturality in the diagram
$$
\xymatrix{
TC(\ell/p) \ar[r]^{i_*} \ar[d]^{\pi^*} & TC(\ell) \ar[d]^{\pi^*} \ar[r]^{j^*}
	& TC(\ell|p^{-1}\ell) \ar[d]^{\pi^*} \\
TC(\Z/p) \ar[r]^{i_*} & TC(\Z_{(p)}) \ar[r]^{j^*} & TC(\Z_{(p)}|\Q)
}
$$
we find that $i_*(1) = 0$ and $i_*(\bar\epsilon_1) = \lambda_1$ in $V(1)_*
TC(\ell)$, since the vertical maps $\pi^*$ are all $(2p-1)$-connected.

To be precise, $i_*(\bar\epsilon_1)$ will be a unit times $\lambda_1$
plus a multiple of $t^{p^2-p}\lambda_2$, but $\lambda_1$ was only
defined modulo such indeterminacy in \cite{AuR02, 1.3}, and this is
a good occasion to make a more definite choice.  In $p$-adic homotopy,
$\pi^* \: \pi_{2p-1} TC(\ell)_p \to \pi_{2p-1} TC(\Z_{(p)})_p
\cong \Z_p\{\lambda_1\}$ is surjective with kernel isomorphic to
$\pi_{2p-2}(\ell)_p \cong \Z_p$, so
$$
\pi_{2p-1} TC(\ell)_p \to V(1)_{2p-1} TC(\ell) = \F_p\{\lambda_1,
  t^{(p-1)p}\lambda_2\}
$$
is surjective, and this ensures that we may, indeed, choose $\lambda_1^K
\in K_{2p-1}(\ell_p)_p$ so that its cyclotomic trace image $\lambda_1
\in \pi_{2p-1} TC(\ell)_p$ reduces mod~$p$ and~$v_1$ to a unit times
$i_*(\bar\epsilon_1)$.

Furthermore, $i_*(t^d v_2) = 0$ for $0<d<p^2-p$, $p \nmid d$, since the
target is zero in these degrees, namely the even degrees strictly between
$2p-2$ and $2p^2-2$ that are not congruent to $-2 \mod 2p$.  Similarly,
$i_*(t^d v_2 \dlog v_1) = 0$ for the same $d$.  Therefore the
cofiber of $i_*$ is
$$
\align
V(1)_* TC(\ell|p^{-1}\ell) &= P(v_2) \otimes E(\partial, \dlog p, \lambda_2) \\
&\quad \oplus P(v_2) \otimes E(\lambda_2) \otimes
	\F_p\{ t^d \lambda_1 \mid 0<d<p\} \\
&\quad \oplus P(v_2) \otimes E(\dlog v_1) \otimes
	\F_p\{ t^d v_2 \dlog p \mid 0<d<p^2-p, p\nmid d\} \\
&\quad \oplus P(v_2) \otimes E(\dlog p) \otimes
	\F_p\{t^{dp}\lambda_2 \mid 0<d<p\}
\endalign
$$
% ((Used \quad in place of \qquad to avoid overfull hbox.))
where the connecting map from~(9.1) takes $\dlog p$ to $1$.

In the middle column the transfer map $\pi_* \: V(1)_* TC(\Z_{(p)}) \to
V(1)_* TC(\ell)$ is zero.  For the natural map $\pi^*$ is
$(2p-1)$-connected and thus surjective, and the composite $\pi_* \circ
\pi^*$ is multiplication by the cyclotomic trace class of $H\Z_{(p)}$
in $\pi_0 TC(\ell)$, realized as the mapping cone of $v_1 \:
\Sigma^{2p-2} \ell \to \ell$, which represents zero already in $\pi_0
K(\ell)$.

We introduce a degree~$1$ class $\dlog v_1$ in the cofiber $V(1)_*
TC(\ell|L)$, and suppose that this notation is compatible with that of
Theorem~8.8, in the sense that $v_2 \cdot \dlog v_1 = \lambda_2$ as
elements of $V(1)_* TC(\ell|L)$, not just formally as in $V(1)_*
TC(\ell/p)$.  This is again to be expected by analogy with the case of
$V(0)_* TC(\Z_{(p)}|\Q)$.   Then the terms $E(\partial, \lambda_1) \otimes
\F_p\{\lambda_2\}$ and $\F_p\{t^d \lambda_1\lambda_2 \mid 0<d<p\}$ in
$V(1)_* TC(\ell)$ become once divisible by $v_2$ in the cofiber, which
we can write as
$$
\align
V(1)_* TC(\ell|L) &= P(v_2) \otimes E(\partial, \lambda_1, \dlog v_1) \\
&\qquad \oplus P(v_2) \otimes E(\dlog v_1) \otimes
	\F_p\{t^d \lambda_1 \mid 0<d<p\} \\
&\qquad \oplus P(v_2) \otimes E(\lambda_1) \otimes
	\F_p\{t^{dp} v_2 \dlog v_1 \mid 0<d<p\} \,.
\endalign
$$
In the left hand column, a very similar argument gives $\pi_* = 0$ and
$$
\align
V(1)_* TC(\ell/p|L/p) &= P(v_2) \otimes
	E(\partial, \bar\epsilon_1, \dlog v_1) \\
&\qquad \oplus P(v_2) \otimes E(\dlog v_1) \otimes
	\F_p\{t^d v_2 \mid 0<d<p^2-p, p\nmid d\} \\
&\qquad \oplus P(v_2) \otimes E(\bar\epsilon_1) \otimes
	\F_p\{t^{dp} v_2 \dlog v_1 \mid 0<d<p\} \,.
\endalign
$$

We are most interested in the right hand column, where again $\pi^*$ is
surjective and $\pi_* \circ \pi^* = 0$, so $\pi_* = 0$.

\proclaim{Proposition 9.2}
There is an exact sequence
$$
0 \to V(1)_* TC(\ell|p^{-1}\ell) @>\rho^*>>
V(1)_* TC(\ff(\ell)) \to V(1)_{*-1} TC(\Z_{(p)}|\Q) \to 0 \,.
$$
If $v_2 \cdot \dlog v_1 = \lambda_2$, where the (right hand) connecting
map takes $\dlog v_1$ to $1$, then there is an isomorphism
$$
\align
V(1)_* TC(\ff(\ell)) &= P(v_2) \otimes E(\partial, \dlog p, \dlog v_1) \\
&\qquad \oplus P(v_2) \otimes E(\dlog v_1) \otimes
	\F_p\{ t^d \lambda_1 \mid 0<d<p\} \\
&\qquad \oplus P(v_2) \otimes E(\dlog v_1) \otimes
	\F_p\{ t^d v_2 \dlog p \mid 0<d<p^2-p, p\nmid d\} \\
&\qquad \oplus P(v_2) \otimes E(\dlog p) \otimes
	\F_p\{t^{dp} v_2 \dlog v_1 \mid 0<d<p\} \,.
\endalign
$$
This is a free module over $P(v_2)$ of rank~$(2p^2+6)$ and of zero
Euler characteristic.
\endproclaim

\demo{Proof}
The assumption in the second clause is that the terms $E(\partial,
\dlog p) \otimes \F_p\{\lambda_2\}$ and $\F_p\{t^d \lambda_1\lambda_2
\mid 0<d<p\}$ in $V(1)_* TC(\ell|p^{-1}\ell)$ become once divisible by $v_2$
in $V(1)_* TC(\ff(\ell))$, and the formula then follows
from those for $V(1)_* TC(\Z_{(p)}|\Q)$ and $V(1)_* TC(\ell|p^{-1}\ell)$.
\qed
\enddemo

\remark{Remark 9.3}
We expect to verify the relation $v_2 \cdot \dlog v_1 = \lambda_2$ in
$V(1)_* TC(\ell|L)$ and $V(1)_* TC(\ff(\ell))$ by means of a logarithmic
geometric model for $TC(\ell|L)$, to be discussed in \cite{Rog:ltc}.
\endremark

\proclaim{Theorem 9.4}
There is an exact sequence
$$
\multline
0 \to \Sigma^{-2} E(\dlog p, \dlog v_1) \otimes \F_p\{\bar\epsilon_1\} \to
	V(1)_* K(\ff(\ell_p)) @>trc>> \\
@>trc>> V(1)_* TC(\ff(\ell)) \to \Sigma^{-1} E(\dlog p, \dlog v_1) \to 0
	\,.
\endmultline
$$
Thus, if $v_2 \cdot \dlog v_1 = \lambda_2$, then
there is an isomorphism of $P(v_2)$-modules
$$
V(1)_* K(\ff(\ell_p)) \equiv P(v_2) \otimes \Lambda_*
$$
modulo the kernel $\Sigma^{-2} E(\dlog p, \dlog v_1) \otimes
\F_p\{\bar\epsilon_1\}$ of $trc$, where
$$
\align
\Lambda_* &= E(\partial v_2, \dlog p, \dlog v_1) \\
&\qquad \oplus E(\dlog v_1) \otimes \F_p\{ t^d \lambda_1 \mid 0<d<p\}
	\\
&\qquad \oplus E(\dlog v_1) \otimes \F_p\{ t^d v_2 \dlog p \mid
	0<d<p^2-p, p\nmid d\} \\
&\qquad \oplus E(\dlog p) \otimes \F_p\{t^{dp} v_2 \dlog v_1 \mid
	0<d<p\} \,.
\endalign
$$
\endproclaim

See Figure~10.3 below for a chart of $\Lambda_*$.

\demo{Proof}
To pass to algebraic $K$-theory, we consider the map from diagram~(3.10)
to diagram~(9.1) induced by the cyclotomic trace map in~(8.9).  The
cofiber is the $3 \times 3$ diagram of cofiber sequences with a copy of
$\Sigma^{-1} H\Z_p$ in each corner of the upper left hand square, and
$V(1)_* \Sigma^{-1} H\Z_p = \Sigma^{-1} E(\bar\epsilon_1)$.  The transfer
maps $i_*$ and $\pi_*$ in this diagram are all trivial, since the
natural maps $i^*$ induce the isomorphisms $W(\Z_p)_F \cong W(\Z/p)_F$
and the maps $\pi^*$ induce identity maps.  So there is a long exact
sequence
$$
\dots \to V(1)_* K(\ff(\ell_p)) @>trc>> V(1)_* TC(\ff(\ell))
@>g>> \Sigma^{-1} E(\bar\epsilon_1, \dlog p, \dlog v_1) \to \dots \,.
$$
We claim that the image of $g$ equals $\Sigma^{-1} E(\dlog p, \dlog
v_1)$.  It is clear that the term $E(\dlog p, \dlog v_1) \otimes
\F_p\{\partial\}$ in $V(1)_* TC(\ff(\ell))$ has this image.  We must argue
that $\Sigma^{-1} E(\dlog p, \dlog v_1) \otimes \F_p\{\bar\epsilon_1\}$ is
not in the image of $g$.

In degree~$(2p-2)$ the class $\Sigma^{-1} \bar\epsilon_1$ could only be hit
by $t\lambda_1 \cdot \dlog v_1$, but $t\lambda_1$ is the image of
$\alpha_1$ in $\pi_{2p-3}(S)$ and $g$ takes $\dlog v_1$ to zero, so $g$
also takes $t\lambda_1 \cdot \dlog v_1$ to zero, since it is a map of
$S$-modules.  In degree~$(2p-1)$ the class $t^{p^2-p} v_2 \dlog v_1$
comes from $t^{p^2-p}\lambda_2$ in $V(1)_* TC(\ell)$, so by naturality
it cannot map to $\Sigma^{-1} \bar\epsilon_1 \dlog p$ or $\Sigma^{-1}
\bar\epsilon_1 \dlog v_1$.  In degree~$2p$ the class $t^{p^2-p} v_2 \dlog p
\dlog v_1$ comes from $t^{p^2-p}\lambda_2 \dlog p$ in $V(1)_*
TC(\ell|p^{-1}\ell)$, so it cannot hit $\Sigma^{-1} \bar\epsilon_1 \dlog p \dlog
v_1$.  This exhausts all possibilities, and concludes the proof.
\qed
\enddemo

The low-degree $v_2$-divisibility assumptions needed for
Proposition~9.2 and Theorem~9.4 become irrelevant upon inverting $v_2$,
since the $V(1)$-homotopy of both $TC(\Z_{(p)}|\Q)$ and $\Sigma^{-1} H\Z_p$ is
$v_2$-torsion.  Hence we have the following unconditional result, for
primes $p\ge5$.

\proclaim{Theorem 9.5}
In $v_2$-periodic homotopy there are isomorphisms of
$P(v_2^{\pm1})$-modules
$$
\align
V(1)_* TC(\ell|p^{-1}\ell) [v_2^{-1}] &\cong V(1)_* TC(\ff(\ell)) [v_2^{-1}] \cong
	V(1)_* K(p^{-1}\ell_p) [v_2^{-1}] \\
&\cong V(1)_* K(\ff(\ell_p)) [v_2^{-1}] \cong P(v_2^{\pm1}) \otimes \Lambda_*
	\,,
\endalign
$$
where $\Lambda_*$ is as in Theorem~9.4.
\endproclaim

\head 10. Galois cohomology \endhead

For motivation, let $F \to E$ be a $G$-Galois extension of local or
global number fields, i.e., finite extensions of $\Q$ or $\Q_p$ for
some prime $p$.  There is an induced $G$-action on $K(E)$ and a natural
map $K(F) \to K(E)^{hG}$.  In their simplest form, the
Lichtenbaum--Quillen conjectures predict that this map induces an
isomorphism in mod~$p$ homotopy, in sufficiently high degrees.  In
other words, the abutment of the homotopy fixed point spectral sequence
$$
E^2_{s,t} = H^{-s}_{gp}(G; V(0)_t K(E))
\Longrightarrow V(0)_{s+t} K(E)^{hG}
$$
is conjectured to agree with $V(0)_* K(F)$ in sufficiently high degrees.
Thomason \cite{Th85, 0.1, 4.1} proved these conjectures for algebraic
$K$-theory with the Bott element inverted (in much greater generality
than the present one), or equivalently, for the $v_1$-periodic
algebraic $K$-theory functor $V(0)_* K(-) [v_1^{-1}]$.  Hence there is a
spectral sequence
$$
E^2_{s,t} = H^{-s}_{gp}(G; V(0)_t K(E) [v_1^{-1}])
\Longrightarrow V(0)_{s+t} K(F) [v_1^{-1}]
$$
for each Galois extension $E$ of $F$, with $G = \Gal(E/F)$.  Passing to
the colimit over all such Galois extensions contained in a separable
closure $\bar F$ of $F$ still gives a spectral sequence
$$
E^2_{s,t} = H^{-s}_{cont}(G_F; V(0)_t K(\bar F) [v_1^{-1}])
\Longrightarrow V(0)_{s+t} K(F) [v_1^{-1}] \,.
$$
Here $H^r_{cont}(G_F; M)$ equals the continuous group cohomology of the
absolute Galois group $G_F = \Gal(\bar F/F)$, for any discrete
$G_F$-module $M$, which is also denoted $H^r_{Gal}(F; M)$ and known as
Galois cohomology.  By Suslin's theorem \cite{Su84}, $K(\bar F)_p
\simeq ku_p$, so $V(0)_t K(\bar F) [v_1^{-1}] \cong V(0)_t KU =
\F_p(t/2)$ equals the $G_F$-module $\mu_p^{\otimes (t/2)}$ for $t$
even, and $0$ for $t$ odd.  Here $\F_p(1) = V(0)_2 K(\bar F) \cong
\mu_p(\bar F) = \mu_p$ is the group of $p$-th roots of unity in $\bar
F$, considered as the $p$-torsion in $K_1(\bar F) = \bar F^\times$.
Thus the spectral sequence above can be rewritten as the Galois descent
spectral sequence
$$
E^2_{s,t} = H^{-s}_{Gal}(F; \F_p(t/2))
\Longrightarrow V(0)_{s+t} K(F) [v_1^{-1}] \,.
\tag 10.1
$$

We will conjecture below that there are similar results for $G$-Galois
extensions $\ff(A) \to \ff(B)$ of the fraction fields of a class of
$K(1)$-local commutative $S$-algebras, containing $L_p$ and $KU_p$, when
considering the $v_2$-periodic algebraic $K$-theory functor $V(1)_* K(-)
[v_2^{-1}]$.  Recall from~(3.5) and~(3.6) that this functor does not
distinguish between the fraction fields of $L$ and $L_p$, or between
the fraction fields of $KU$, $KU_{(p)}$ or $KU_p$, since it vanishes
on finite, global and local fields in the ordinary sense.  Note that we
here shift our focus from connective spectra like $\ell$ and $\ell_p$,
which are convenient for topological cyclic homology calculations, to
$K(1)$-local spectra like $L_p$, which are convenient for Galois theory.
As discussed in Remark~2.6, $\ff(\ell_p) = \ff(L_p)$ and
$\ff(ku_p) = \ff(KU_p)$, so the shift is only one of perspective.

More precisely, we might consider a $K(1)$-local pro-Galois extension
$L_{K(1)}S \to C$ of commutative $S$-algebras in the sense of
\cite{Rog08, 4.1}, with profinite Galois group $H$, and a pair $K
\subset N$ of closed subgroups of $H$, with $K$ normal in $N$.  Letting
$A = C^{hN}$, $B = C^{hK}$ and $G = N/K$ we obtain a diagram
$$
J_p = L_{K(1)}S \to A @>G>> B \to C \,.
$$
Here $J_p$ is the $p$-complete image-of-$J$ spectrum, which plays an
analogous initial role in the $K(1)$-local category as $H\Q = L_0S$
does in the rational category.  There is a $K(1)$-local
pro-$\Gamma$-Galois extension $J_p \to L_p$, where $\Gamma = 1 + p\Z_p
\subset \Z_p^\times$ acts on $L_p$ through $p$-adic Adams operations,
and $\Spec(J_p) = [\Spec(L_p)/\Gamma]$ is essentially the derived orbit
stack for the corresponding $\Gamma$-action on the derived scheme
$\Spec(L_p)$.  Similar remarks apply in a non-commutative sense to $J/p
\to L/p$, so $\Spec(\ff(J_p))$ may be interpreted as the derived orbit
stack for $\Gamma$ acting on $\Spec(\ff(L_p))$.  In both cases, some
care in interpretation is needed, since $\Gamma$ is a profinite group.

For such a $K(1)$-local $G$-Galois extension $A \to B$ as above, we
expect that the map $K(A) \to K(B)^{hG}$ induces an isomorphism in
$V(1)$-homotopy, in sufficiently high degrees, so that there is a weak
equivalence
$$
V(1)_* K(A) [v_2^{-1}] @>\simeq>> V(1)_* K(B)^{hG} [v_2^{-1}]
$$
and a homotopy fixed point spectral sequence
$$
E^2_{s,t} = H^{-s}_{gp}(G; V(1)_t K(B) [v_2^{-1}])
\Longrightarrow V(1)_{s+t} K(A) [v_2^{-1}] \,.
$$
In the special case of the $\Delta$-Galois extension $L_p \to KU_p$,
this conjecture is somewhat trivially compatible with the predicted
formula~(3.8), since $H^r_{gp}(\Delta; P(b^{\pm1}))$ equals
$P(v_2^{\pm1})$ for $r=0$ and is $0$ otherwise.

Considering $A$ and $B$ as valuation rings in their fraction fields, when
defined, these extensions correspond to unramified Galois extensions of
the fraction fields.  It is difficult to use the obstruction theories
of \cite{GH04} or \cite{Rob03} to construct ramified extensions
in commutative $S$-algebras, since the ramification creates highly
nontrivial obstruction groups of Andr{\'e}--Quillen type.  In fact,
there are only a few unramified Galois extensions of $L_p$ and $KU_p$.
By \cite{BR08, 1.1, 7.9}, the extended Lubin--Tate spectrum $KU^{nr}_p$,
with $\pi_* KU^{nr}_p = \W(\bar\F_p)[u^{\pm1}]$, is the maximal pro-Galois
extension of $KU_p$, both as commutative $S$-algebras and as $K(1)$-local
commutative $S$-algebras.

More generally, we suppose that there is a notion of $K(1)$-local
Galois extensions in an extended framework that contains the fraction
fields of commutative $S$-algebras discussed in Section~3, but which
also allows for ramified Galois extensions.  Then let $\Omega_1$ be
a maximal connected $K(1)$-local pro-Galois extension of $\ff(L_p)$,
or equivalently of $\ff(J_p)$.  In other words, $\Omega_1$ is to be a
separable closure of the fraction field of $L_p$ (or $J_p$).  Then for
each intermediate such $K(1)$-local $G$-Galois extension
$$
\ff(J_p) = \ff(L_{K(1)}S) \to \ff(A) @>G>> \ff(B) \to \Omega_1
$$
(where $\ff(B)$ or $K(\ff(B))$ might only exist in the extended
framework), we expect to have a weak equivalence
$$
V(1)_* K(\ff(A)) [v_2^{-1}] @>\simeq>> V(1)_* K(\ff(B))^{hG} [v_2^{-1}]
$$
and a homotopy fixed point spectral sequence
$$
E^2_{s,t} = H^{-s}_{gp}(G; V(1)_t K(\ff(B)) [v_2^{-1}])
\Longrightarrow V(1)_{s+t} K(\ff(A)) [v_2^{-1}] \,,
$$
as above.  Fixing $\ff(A)$ and passing to the colimit over all such
Galois extensions contained in $\Omega_1$, we are then led to a spectral
sequence
$$
E^2_{s,t} = H^{-s}_{cont}(G_{\ff(A)}; V(1)_t K(\Omega_1) [v_2^{-1}])
\Longrightarrow V(1)_{s+t} K(\ff(A)) [v_2^{-1}] \,,
$$
where $G_{\ff(A)} = \Gal(\Omega_1/\ff(A))$ is the absolute Galois
group, defined as the limit of the Galois groups $G$ above.

In the case $A = L_p$ for $p\ge5$, we achieved a concrete calculation
of the abutment of this spectral sequence in Theorem~9.5, by methods
completely independent of the conjectural Galois descent properties.
We are therefore entitled to make an educated guess at the structure of
the $E^2$-term of this spectral sequence, which in turn makes strong
suggestions about the form of the Galois module $V(1)_*
K(\Omega_1) [v_2^{-1}]$, and about the absolute Galois group $G_{\ff(L_p)}$.

We therefore study the spectral sequence
$$
E^2_{s,t} = H^{-s}_{cont}(G_{\ff(L_p)}; V(1)_t K(\Omega_1) [v_2^{-1}])
\Longrightarrow V(1)_{s+t} K(\ff(L_p)) [v_2^{-1}] \,,
\tag 10.2
$$
which we assume to be an algebra spectral sequence, and one that
collapses at the $E^2$-term for $p\ge5$, as is the case for the Galois
descent spectral sequence~(10.1) for any local number field, as well
as for any global number field when $p\ge3$.  The abutment $V(1)_*
K(\ff(L_p)) [v_2^{-1}]$ is isomorphic to $P(v_2^{\pm1})$ tensored with
$\Lambda_*$, as given in Theorem~9.4.  Hence the $E^2 = E^\infty$-term
of~(10.2) is the associated graded for a ``Galois'' filtration of this
$P(v_2^{\pm1})$-module.  Since $v_2$ is invertible, it must be represented
in Galois filtration $s=0$.  By consideration of localization sequences
in motivic cohomology, and the relation of motivic cohomology to Galois
cohomology, we can partially justify why the other algebra generators
$\partial v_2$, $\dlog p$, $\dlog v_1$, $t^d \lambda_1$, $t^d v_2 \dlog p$
and $t^{dp} v_2 \dlog v_1$ must all be represented in Galois filtration
$s=-1$.  It follows that the entire spectral sequence is concentrated in
the four columns $-3 \le s \le 0$, with the product $\partial v_2 \cdot
\dlog p \cdot \dlog v_1$ as the only generator in filtration $s=-3$.

We display the placement in the $E^2$-term of the
$P(v_2^{\pm1})$-module generators for $V(1)_* K(\ff(L_p)) [v_2^{-1}]$, or
equivalently the $\F_p$-basis for $\Lambda_*$, in Figure~10.3 below.
To make room on the page, we write $\dlp$ and $\dlv$ for $\dlog p$ and
$\dlog v_1$, respectively.  Multiple generators in the same bidegree
are separated by colons.  The symbol $(v_2)$ indicates the placement of
this $v_2$-multiple of $1$.  Also to save space, the chart is drawn for
$p=3$, even if we are really assuming $p\ge5$.  
\midinsert
$$
\xymatrix@=0.5pc{
\partial v_2 \dlp \dlv & . & . & . & 2p^2+2 \\ 
. & \partial v_2 \dlp : \partial v_2 \dlv : t v_2 \dlp \dlv & . & . & 2p^2 \\
. & t^2 v_2 \dlp \dlv & \partial v_2 : t v_2 \dlp & (v_2) & 2p^2-2 \\
. & t^p v_2 \dlp \dlv & t^2 v_2 \dlp & . \\
. & t^{p+1} v_2 \dlp \dlv & t^p v_2 \dlv & . \\
. & t^{p+2} v_2 \dlp \dlv & t^{p+1} v_2 \dlp & . \\
. & t^{2p} v_2 \dlp \dlv & t^{p+2} v_2 \dlp & . \\
. & t \lambda_1 \dlv & t^{2p} v_2 \dlv & . & 2p \\
. & \dlp \dlv : t^2 \lambda_1 \dlv & t \lambda_1 & . & 2p-2 \\
. & . & \dlp : \dlv : t^2 \lambda_1 & . & 2\\
. & . & . & \qquad 1 \qquad \vrule height12pt depth6pt width0pt & 0 \\
-3 & -2 & -1 & 0 & s \backslash t
}
$$
\botcaption
{Figure 10.3: $P(v_2^{\pm1})$-basis for $E^2_{s,t} \Longrightarrow
V(1)_{s+t} K(\ff(L_p)) [v_2^{-1}]$}
\endcaption
\endinsert

From these considerations, it is apparent that $V(1)_t K(\Omega_1)
[v_2^{-1}]$ will be nontrivial for each even $t$, and zero for each odd
$t$.  Let us briefly write $M_t$ for this Galois module.  The presence
of nonzero products by $\dlog p$ or $\dlog v_1$ from $E^2_{s,t}$ to
$E^2_{s-1,t+2}$, for each even $t$ and some $s$, implies that the pairing
$$
M_t \otimes M_2 \to M_{t+2}
$$
(tensor product over $M_0$, here and below) is also nontrivial for each
even $t$, since the product in the spectral sequence should be induced
by the group cohomology cup product and this pairing of the
coefficients.  This can most simply be the case if these pairings
induce an isomorphism
$$
M_2^{\otimes (t/2)} @>\cong>> M_t
$$
for all even $t$, which we now assume.  Since $E^2_{0,t}=0$ for
$0<t<2p^2-2$, the Galois action on $M_2^{\otimes (t/2)}$ should have no
invariants for these $t$, but $M_2^{\otimes (p^2-1)}$ should have
trivial action.  This indicates that the $M_0$-linear automorphism
group of $M_2$ should have exponent $(p^2-1)$, which excludes $\F_p$,
but strongly suggests the minimal example $M_2 = \F_{p^2}(1)$, with
some Galois automorphism of $\Omega_1$ over $\ff(L_p)$ acting by
multiplication by a generator of the group of units $\F_{p^2}^\times
\cong \Z/(p^2-1)$.  Then $M_0 = \F_{p^2}(0)$ with the trivial Galois
action, and letting $u \in M_2 = V(1)_2 K(\Omega_1) [v_2^{-1}]$ be a
generator for $M_2$ as an $M_0$-module, we see that $u^{p^2-1}$
generates $M_{2p^2-2} = M_0\{v_2\}$, hence equals a unit multiple of
$v_2$.  Thus we deduce that the coefficient module in~(10.2) can most
simply be
$$
V(1)_* K(\Omega_1) [v_2^{-1}] = M_* = \F_{p^2}[u^{\pm1}] = \pi_*(K_2)
$$
where $K_2$ is the $2$-periodic form of the Morava $K$-theory spectrum
$K(2)$, related to the (elliptic) Lubin--Tate spectrum $E_2$ with
$\pi_* E_2 = \W(\F_{p^2})[[u_1]][u^{\pm1}]$ by $K_2 \simeq V(1) \wedge
E_2$.  Here $K(2) \to K_2$ takes $v_2$ to $u^{p^2-1}$, and $v_1$ acts
on $V(0) \wedge E_2$ by multiplication by $u_1 u^{p-1}$.  In other
words, $\pi_*(K_2) \cong \pi_*(E_2)/(p, u_1) = \pi_*(E_2)/(p, v_1)$.
We are therefore led to the formula
$$
V(1)_* K(\Omega_1) [v_2^{-1}] \cong V(1)_* E_2 \,.
$$
By \cite{HoSt99, 7.2}, the left hand side is the $V(1)$-homotopy
of $L_{K(2)} K(\Omega_1)$, and $E_2$ is $K(2)$-local.

\proclaim{Conjecture 10.4}
$K(\Omega_1)$ is a connective form of $E_2$, in the sense that
there is an equivalence $L_{K(2)} K(\Omega_1) \simeq E_2$.
\endproclaim

This would generalize Suslin's theorem, which can be formulated as
an equivalence $L_{K(1)} K(\bar F) \simeq KU_p = E_1$.  Assuming
this conjecture, we write $\F_{p^2}(t/2)$ for the Galois module
$V(1)_t K(\Omega_1) [v_2^{-1}] = V(1)_t E_2 = \pi_t(K_2)$, which is
$\F_{p^2}\{u^{t/2}\}$ for $t$ even, and zero for $t$ odd.  Then the
Galois descent spectral sequence takes the following form.

\proclaim{Conjecture 10.5}
For $K(1)$-local fraction fields $\ff(A)$ contained in the separable
closure $\Omega_1$ of $\ff(J_p) = \ff(L_{K(1)}S)$ there is a natural
algebra spectral sequence
$$
E^2_{s,t} = H^{-s}_{Gal}(\ff(A); \F_{p^2}(t/2))
\Longrightarrow V(1)_{s+t} K(\ff(A)) [v_2^{-1}] \,,
$$
where $H^r_{Gal}(\ff(A); M) = H^r_{cont}(G_{\ff(A)}; M)$ is the
continuous group cohomology of the absolute Galois group of
$\ff(A)$.
\endproclaim

The apparent self-duality of the $\F_p$-algebra $\Lambda_*$ displayed
in Figure~10.3 suggests that there is such a self-duality in the Galois
cohomology of most or all $K(1)$-local fraction fields contained in
$\Omega_1$.  Recall \cite{Se97, \S II.5} that for each $p$-local number field
$F$ there is a canonical isomorphism $H^2_{Gal}(F, \F_p(1))
\cong \F_p$, and that by the Tate--Poitou arithmetic duality theorem
the cup product
$$
H^r_{Gal}(F; \F_p(i)) \otimes H^{2-r}_{Gal}(F; \F_p(1-i)) @>\cup>>
H^2_{Gal}(F; \F_p(1)) \cong \F_p
$$
is a perfect pairing for each $r$ and $i$.

\proclaim{Conjecture 10.6}
For each finite $K(1)$-local Galois extension $\ff(A)$ of $\ff(L_p)$
there is a canonical isomorphism
$$
H^3_{Gal}(\ff(A); \F_{p^2}(2)) \cong \F_p
$$
and the cup product
$$
H^r_{Gal}(\ff(A); \F_{p^2}(i)) \otimes H^{3-r}_{Gal}(\ff(A); \F_{p^2}(2-i))
	@>\cup>> H^3_{Gal}(\ff(A); \F_{p^2}(2)) \cong \F_p
$$
is a perfect pairing for each $r$ and $i$.
\endproclaim

\remark{Remark 10.7}
Assuming these conjectures, it would be interesting to interpret the
Galois module $\F_{p^2}(1) \cong  V(1)_2 K(\Omega_1) [v_2^{-1}]$ directly
in terms of $K(\Omega_1)$, perhaps as suitable torsion points of a
derived algebraic or formal group over a Galois extension $\ff(A)$ of
$\ff(L_p)$.  In particular, it would be illuminating to find a (minimal)
such extension for which
$$
V(1)_2 K(\ff(A)) [v_2^{-1}] \to V(1)_2 K(\Omega_1) [v_2^{-1}]
$$
is surjective.  This would play the role of the $p$-th
cyclotomic field $\Q(\zeta_p)$ in the case of mod~$p$ algebraic $K$-theory
of number fields, since there is a Bott element $\beta \in V(0)_2
K(\Q(\zeta_p)) [v_1^{-1}]$ mapping to $u \in V(0)_2 K(\bar\Q) [v_1^{-1}]
\cong \F_p(1)$, up to a unit in $\F_p$.  We recall that the choice of
a $p$-th root of unity $\zeta_p$ determines a group homomorphism
$C_p \to GL_1(\Q(\zeta_p))$,
a map of $E_\infty$ spaces
$$
BC_p \to BGL_1(\Q(\zeta_p)) \to \Omega^\infty_\otimes K(\Q(\zeta_p)) \,,
$$
and a map of commutative $S$-algebras $S[BC_p] \to K(\Q(\zeta_p))$.
The Bott class $\beta \in V(0)_2(BC_p)$ satisfies $\beta^p = v_1 \beta$,
as can be seen by mapping along $V(0)_*(BC_p) \to K(1)_*(BC_p)$, and its
image in $V(0)_2 K(\Q(\zeta_p))$ satisfies $\beta^{p-1} = v_1 = u^{p-1}$.

The calculations of \cite{Au:tcku, 1.1} show that a generator in the
$(p+1)$-st tensor power $\F_{p^2}(p+1) \cong V(1)_{2p+2} K(\Omega_1)
[v_2^{-1}]$ of the Galois module $\F_{p^2}(1)$ is realized by the
unramified $\Delta$-Galois extension $L_p \to KU_p$, and its fraction
field analogue $\ff(L_p) \to \ff(KU_p)$.  For there is a higher Bott
element $b \in V(1)_{2p+2} K(ku)$ that satisfies $b^{p-1} = -v_2$, whose
image in $V(1)_{2p+2} K(\Omega_1) [v_2^{-1}] \cong \F_{p^2}\{u^{p+1}\}$
must have the form $w u^{p+1}$, for some element $w \in \F_{p^2}$ with
$w^{p-1} = -1$.  Hence $ku_p \to \Omega_1$ realizes the $\F_p$-subalgebra
generated by $wu^{p+1}$ and $u^{\pm(p^2-1)}$ in $\F_{p^2}[u^{\pm1}]$.
The higher Bott element is constructed from the map
$$
K(\Z, 2) \simeq BU(1) \to BU_\otimes \to GL_1(ku)
$$
that interprets the left hand abelian group as the $v_1$-torsion points
in the homotopy units of $ku$.  More precisely, $BU_\otimes \simeq
BU(1) \times BSU_\otimes$, and $p$-locally $BSU \simeq BSU_\otimes$
\cite{AP76} is $v_1$-torsion free, so $K(\Z, 2)$ is the $v_1$-torsion
in $BU_\otimes$.  The delooped $E_\infty$ map $K(\Z, 3) \to BGL_1(ku) \to
\Omega^\infty_\otimes(ku)$ is adjoint to a map $S[K(\Z, 3)] \to K(ku)$ of
commutative $S$-algebras.  The Bott element $\beta \in V(1)_{2p+2}(K(\Z,
3))$ satisfies $\beta^p = -v_2 \beta$, as can be detected in $K(2)_*(K(\Z,
3))$ using \cite{RaWi80, 9.2, 12.1}, and its image $b$ in $V(1)_{2p+2} K(ku)$
satisfies $b^{p-1} = -v_2 = -u^{p^2-1}$.  A generator $\delta$ of
$\Delta$ multiplies $b$ by a generator of $\F_p^\times \cong \Z/(p-1)$,
so an extension of $\delta$ to a Galois automorphism of $\Omega_1$
must multiply $u$ by a generator of $\F_{p^2}^\times \cong \Z/(p^2-1)$.
Hence the absolute Galois group of $\ff(L_p)$ will act faithfully on
$\F_{p^2}(1)$.

Now suppose that there is a homotopy-commutative $H$-space $X$ with
$$
K(2)_*(X) \cong K(2)_*[x]/(x^{p^2} = v_2^j x) \,,
$$
realizing one of the finite irreducible commutative Morava--Hopf algebras
\cite{SW98, 1.4}.  Here $|x| = 2j$, and we assume that $(j, p+1) = 1$.
Let $\Omega X$ be the loop group of~$X$, and let $S[\Omega X]$ be the
spherical monoid ring.  As usual, there will be a homotopy commutative
$S$-algebra map
$$
S[X] \to K(S[\Omega X]) = A(X) \,,
$$
where $A(X)$ is as in \cite{Wa85}.
The unit map $S \to E(2)$ induces a map $V(1)_*(X) \to K(2)_*(X)$,
and we assume that $x$ lifts to a class in $V(1)_{2j}(X)$, with image
$\xi \in V(1)_{2j} K(S[\Omega X])$, still satisfying $\xi^{p^2} =
v_2^j \xi$.  If one can extend $S[\Omega X]$ to a $K(1)$-local Galois
extension $\ff(A)$ of $\ff(J_p)$ or $\ff(L_p)$, then the image $\xi \in
V(1)_{2j} K(\ff(A))$ will map to $u^j \in V(1)_{2j} K(\Omega_1)[v_2^{-1}]
\cong \F_{p^2}(j)$, up to a unit in $\F_{p^2}$.  If we also arrange
that $\ff(A)$ contains $\ff(KU_p)$, then $\ff(A)$ will realize all of
$\F_{p^2}[u^{\pm1}]$, since this is generated as an $\F_p$-algebra by
the classes $b \mapsto wu^{p+1}$, $\xi \mapsto u^j$ and $v_2^{\pm1}
\mapsto u^{\pm(p^2-1)}$.
\endremark

%\remark{Remark 10.?}
%At odd primes, the $E(1)$-local Picard group only contains the usual
%spheres $L_1 S^n$ for $n \in \Z$, by \cite{HoSa99, Thm.~A}, but the
%$K(1)$-local Picard group is $\Z/(2p-2) \times \Z_p$, by \cite{HMS94},
%with the embedding $k \mapsto L_{K(1)} S^{2(p-1)k}$ for $k \in \Z$
%extending over $\Z \subset \Z_p$.  For example, with $k = 1/(p^2-1)
%\in \Z_p$ there is an invertible $K(1)$-local spectrum $X = X_{1/(p+1)}
%= L_{K(1)} S^{2/(p+1)}$ with a $K(1)$-equivalence $X^{\wedge (p+1)} \to
%S^2$.  So in the $K(1)$-local category, we might form
%$$
%A = \bigvee_{k=0}^p KU_p \wedge X^{\wedge k}
%$$
%as an associative $KU_p$-algebra.  The challenge remains to make $A$ a
%commutative $KU_p$-algebra, or a $\Z/(p+1)$-Galois extension of $KU_p$.
%By \cite{BR08}, this can only be possible in the putative extended
%framework of $K(1)$-local fraction fields, yielding $\ff(KU_p) \to
%\ff(A)$.
%\endremark

\remark{Remark 10.8}
The optimistic reader can now extend these conjectures for all
$n\ge1$ to Galois extensions in the $K(n)$-local category, as seen
by the $v_{n+1}$-periodic algebraic $K$-theory functor $V(n)_* K(-)
[v_{n+1}^{-1}]$, where the Smith--Toda complex $V(n)$ might be replaced
with any fixed type $(n+1)$ finite complex (to ensure that it exists).
The extension $J_p \to L_p$ is best replaced with the $K(n)$-local
pro-Galois extension $L_{K(n)}S \to E_n$ with Galois group $\G_n$
constructed in \cite{DH04}, see \cite{Rog08, 5.4} and \cite{AuR08}.
We suggest writing $\Omega_n$ for the $K(n)$-local separable closure of
the fraction field of $L_{K(n)}S$ or $E_n$.  Then Conjecture~10.4 may
be extended to the prediction
$$
L_{K(n+1)} K(\Omega_n) \simeq E_{n+1} \,,
$$
and similarly for Conjecture~10.5.  If correct, then $e_{n+1} =
K(\Omega_n)_p$ is a good connective form of the $p$-complete
commutative $S$-algebra $E_{n+1}$.
\endremark

\enddocument